\documentstyle[11pt]{article}

\newtheorem{theo}{{\bf Theorem}}

\newtheorem{lemma}{{\bf Lemma}}[section]

\newtheorem{coro}[lemma]{{\bf Corollary}}
\newtheorem{defn}{{\bf Definition}}
\newtheorem{remark}{{\bf Remark}}

\font\bbb=msbm10 scaled\magstep1

\newcommand{\CC}{\mbox{\bbb C}}

\newcommand{\RR}{\mbox{\bbb R}}
\newcommand{\ZZ}{\mbox{\bbb Z}}

\newcommand{\PPP}{P^{\!\!\!^{^{\circ}}}}
\newcommand{\intalpha}{\alpha^{\!\!\!^{{\circ}}}}
\newcommand{\intbeta}{\beta^{\!\!\!^{^{\circ}}}}

\def\C{{\rm \kern.24em
 \vrule width.02em
   height1.5ex depth-.05ex
 \kern-.26em C}}

 \def\CCC{{\rm \kern.24em
 \vrule width.08em
   height1.5ex depth-.08ex
 \kern-.36em C}}

\begin{document}

\title{\bf A triangulation of \boldmath{$\CCC {\rm P}^{3}$} as
symmetric cube of \boldmath{${\rm S}^{2}$}}
\author{{\bf Bhaskar Bagchi}$^{\rm a}$, {\bf Basudeb Datta}$^{\rm b,
1}$ }

\date{}

\maketitle

\vspace{-3mm}

\noindent {\small $^{\rm a}$Theoretical Statistics and Mathematics
Unit, Indian Statistical Institute,  Bangalore 560\,059, India

\smallskip

\noindent $^{\rm b}$Department of Mathematics, Indian Institute of
Science, Bangalore 560\,012,  India}

\footnotetext[1]{Supported by UGC-SAP/DSA-IV.

\vspace{2mm}

{\em E-mail addresses:} bbagchi@isibang.ac.in (B. Bagchi),
dattab@math.iisc.ernet.in (B. Datta). }

\begin{center}

\date{December 15, 2010}

\end{center}

\smallskip

\hrule

\bigskip

\noindent {\bf Abstract.} {The symmetric group $S_3$ acts on
$S^{\,2} \times S^{\,2} \times S^{\,2}$ by coordinate permutation,
and the quotient space $(S^{\,2}\times S^{\,2}\times S^{\,2})/S_3$
is homeomorphic to the complex projective space $\CC P^{\,3}$. In
this paper, we construct an 124-vertex simplicial subdivision
$(S^{\,2}\times S^{\,2} \times S^{\,2})_{124}$ of the 64-vertex
standard cellulation $S^{\,2}_4 \times S^{\,2}_4 \times S^{\,2}_4$
of $S^{\,2} \times S^{\,2} \times S^{\,2}$, such that the
$S_3$-action on this cellulation naturally extends to an action on
$(S^{\,2}\times S^{\,2} \times S^{\,2})_{124}$. Further, the
$S_3$-action on $(S^{\,2}\times S^{\,2} \times S^{\,2})_{124}$ is
``good", so that the quotient simplicial complex $(S^{\,2} \times
S^{\,2} \times S^{\,2})_{124}/S_3$ is a 30-vertex triangulation
$\CC P^{\,3}_{30}$ of $\CC P^{\,3}$. In other words, we construct
a simplicial realization $(S^{\,2}\times S^{\,2} \times
S^{\,2})_{124} \to \CC P^{\,3}_{30}$ of the branched covering
$S^{\,2} \times S^{\,2} \times S^{\,2} \to \CC P^{\,3}$. Finally,
we apply the BISTELLAR program of Lutz on $\CC P^{\,3}_{30}$,
resulting in an 18-vertex 2-neighbourly triangulation $\CC
P^{\,3}_{18}$ of $\CC P^{\,3}$. The automorphism group of $\CC
P^{\,3}_{18}$ is trivial. It may be recalled that, by a result of
Arnoux and Marin, any triangulation of $\CC P^{\,3}$ requires at
least 17 vertices. So, $\CC P^{\,3}_{18}$ is close to
vertex-minimal, if not actually vertex-minimal. Moreover, no
explicit triangulation of $\CC P^{\,3}$ was known so far.}

\bigskip

{\small

\noindent  {\em MSC 2010:} 57Q15, 57R05, 57M60.

\smallskip

\noindent {\em Keywords:} Triangulated manifolds; Complex
projective space; Symmetric power; Product of 2-spheres.
} \bigskip

\hrule

\section{Main Results}

\noindent {\bf Description of \boldmath{$(S^{\,2}\times S^{\,2}
\times S^{\,2})_{124}$}\,:} Its vertex-set is $V_{124} = \{x_{ijk}
~ : 1\leq i, j, k \leq 4\} \cup \{v_{pqr} ~ : ~ 1\leq p, q, r \leq
4, \, p, q, r$ are distinct$\} \cup \{u_{ij} ~ : ~ 1\leq i, j \leq
6\}$. The group $S_3 \times A_4$ acts as an automorphism group,
where the actions of $S_3$ and $A_4$ on the vertices are as
follows. On the vertices other than $u_{ij}$, $S_3$ acts by
permuting the positions of the three subscripts while $A_4$ acts
by permuting the values of these subscripts (which are elements of
$\{1, 2, 3, 4\}$). The action on the vertices $u_{ij}$ of the
generators $\alpha = (1, 2, 3)$, $\beta = (1, 2, 4)$ of $A_4$ and
$\gamma = (1, 2, 3)$, $\delta = (1, 2)$ of $S_3$ is as follows\,:
\begin{eqnarray*}
\alpha = \prod_{i=1}^6(u_{i1}, u_{i2}, u_{i4})(u_{i6}, u_{i5},
u_{i+3,3}), && \beta = \prod_{i=1}^6(u_{i1}, u_{i3}, u_{i+3,2})
(u_{i4}, u_{i5}, u_{i+3,6}), \\
\gamma = \prod_{i=1}^6(u_{1i}, u_{2i}, u_{3i})(u_{4i}, u_{5i},
u_{6i}), && \delta = \prod_{i=1}^6(u_{1i}, u_{6i}) (u_{2i},
u_{5i}) (u_{3i}, u_{4i}).
\end{eqnarray*}
(Here, summation in the subscripts is modulo 6.)

Modulo this group the facets (maximal simplices) are the following\,:
 \begin{eqnarray*}
  x_{111}x_{133}x_{131}x_{433}x_{233}u_{22}u_{11},
& x_{111}x_{133}x_{131}x_{233}v_{312}u_{22}u_{11},
& x_{111}x_{133}x_{131}v_{312}x_{121}u_{22}u_{11},   \\
  x_{111}x_{133}x_{223}x_{433}x_{233}u_{22}u_{11},
& x_{111}x_{133}x_{223}x_{233}v_{312}u_{22}u_{11},
& x_{111}x_{133}x_{223}v_{312}x_{121}u_{22}u_{11},   \\
  x_{111}x_{221}x_{131}x_{433}x_{233}u_{22}u_{11},
& x_{111}x_{221}x_{131}x_{233}v_{312}u_{22}u_{11},
& x_{111}x_{221}x_{131}v_{312}x_{121}u_{22}u_{11},   \\
  x_{111}x_{221}x_{223}x_{433}x_{233}u_{22}u_{11},
& x_{111}x_{221}x_{223}x_{233}v_{312}u_{22}u_{11},
& x_{111}x_{221}x_{223}v_{312}x_{121}u_{22}u_{11},   \\
  x_{111}x_{433}x_{233}x_{221}x_{223}u_{65}u_{11},
& x_{111}x_{433}x_{233}x_{223}v_{342}u_{65}u_{11},
& x_{111}x_{433}x_{233}v_{342}x_{434}u_{65}u_{11},   \\
  x_{111}x_{433}x_{424}x_{221}x_{223}u_{65}u_{11},
& x_{111}x_{433}x_{424}x_{223}v_{342}u_{65}u_{11},
& x_{111}x_{433}x_{424}v_{342}x_{434}u_{65}u_{11},   \\
  x_{111}x_{224}x_{233}x_{221}x_{223}u_{65}u_{11},
& x_{111}x_{224}x_{233}x_{223}v_{342}u_{65}u_{11},
& x_{111}x_{224}x_{233}v_{342}x_{434}u_{65}u_{11},   \\
  x_{111}x_{224}x_{424}x_{221}x_{223}u_{65}u_{11},
& x_{111}x_{224}x_{424}x_{223}v_{342}u_{65}u_{11},
& x_{111}x_{224}x_{424}v_{342}x_{434}u_{65}u_{11},   \\
  x_{111}x_{133}x_{131}x_{424}x_{121}u_{22}u_{11},
& x_{111}x_{133}x_{131}x_{424}x_{433}u_{22}u_{11},
& x_{111}x_{133}x_{131}x_{424}x_{121}x_{124}u_{11},  \\
  x_{111}x_{133}x_{131}x_{424}x_{124}x_{134}u_{11},
& x_{111}x_{133}x_{131}x_{424}x_{134}x_{434}u_{11},
& x_{111}x_{133}x_{131}x_{424}x_{434}x_{433}u_{11},  \\
  x_{111}x_{131}x_{224}x_{424}x_{221}u_{65}u_{11},
& x_{111}x_{131}x_{224}x_{424}x_{434}u_{65}u_{11},
 & x_{111}x_{131}x_{224}x_{424}x_{221}x_{121}u_{11}, \\
 x_{111}x_{131}x_{224}x_{424}x_{121}x_{124}u_{11},
& x_{111}x_{131}x_{224}x_{424}x_{124}x_{134}u_{11},
& x_{111}x_{131}x_{224}x_{424}x_{134}x_{434}u_{11}, \\
 x_{111}x_{131}x_{133}x_{233}x_{434}x_{134}u_{11}, &
 x_{111}x_{131}x_{133}x_{233}x_{434}x_{433}u_{11}, &
 x_{111}x_{121}x_{223}x_{224}x_{424}x_{124}u_{11},\\
 x_{111}x_{121}x_{223}x_{224}x_{424}x_{221}u_{11}, &
 x_{111}x_{131}x_{224}x_{233}x_{434}x_{134}u_{11},&
 x_{111}x_{121}x_{133}x_{223}x_{424}x_{124}u_{11}, \\
 x_{111}x_{131}x_{233}x_{224}x_{221}u_{65}u_{11}, &
 x_{111}x_{131}x_{233}x_{224}x_{434}u_{65}u_{11}, &
 x_{111}x_{131}x_{233}x_{433}x_{221}u_{65}u_{11}, \\
 x_{111}x_{131}x_{233}x_{433}x_{434}u_{65}u_{11}, &
 x_{111}x_{223}x_{424}x_{133}x_{121}u_{22}u_{11}, &
 x_{111}x_{223}x_{424}x_{133}x_{433}u_{22}u_{11}, \\
 x_{111}x_{223}x_{424}x_{221}x_{121}u_{22}u_{11}, &
 x_{111}x_{223}x_{424}x_{221}x_{433}u_{22}u_{11}, &
 x_{111}x_{131}x_{424}x_{433}x_{221}u_{65}u_{11}, \\
 x_{111}x_{131}x_{424}x_{433}x_{434}u_{65}u_{11}, &
 x_{111}x_{131}x_{221}x_{424}x_{121}u_{22}u_{11}, &
 x_{111}x_{131}x_{221}x_{424}x_{433}u_{22}u_{11}, \\
 x_{111}x_{221}x_{421}x_{424}x_{121}x_{131}u_{22}, &
 x_{111}x_{221}x_{421}x_{424}x_{121}x_{423}u_{22}, &
 x_{111}x_{221}x_{421}x_{424}x_{433}x_{131}u_{22}, \\
 x_{111}x_{221}x_{421}x_{424}x_{433}x_{423}u_{22}, &
 x_{111}x_{223}x_{423}x_{424}x_{133}x_{121}u_{22}, &
 x_{111}x_{223}x_{423}x_{424}x_{133}x_{433}u_{22}, \\
 x_{111}x_{223}x_{423}x_{424}x_{221}x_{121}u_{22}, &
 x_{111}x_{223}x_{423}x_{424}x_{221}x_{433}u_{22}, &
 x_{111}x_{131}x_{431}x_{433}x_{221}x_{231}u_{65}, \\
 x_{111}x_{131}x_{431}x_{433}x_{221}x_{424}u_{65}, &
 x_{111}x_{131}x_{431}x_{433}x_{434}x_{231}u_{65}, &
 x_{111}x_{131}x_{431}x_{433}x_{434}x_{424}u_{65}, \\
 x_{111}x_{131}x_{231}x_{233}x_{224}x_{221}u_{65}, &
 x_{111}x_{131}x_{231}x_{233}x_{224}x_{434}u_{65}, &
 x_{111}x_{131}x_{231}x_{233}x_{433}x_{221}u_{65},  \\
 x_{111}x_{131}x_{231}x_{233}x_{433}x_{434}u_{65}, &
 x_{111}x_{131}x_{221}x_{421}x_{424}x_{431}x_{433}, &
 x_{111}x_{131}x_{224}x_{233}x_{234}x_{434}x_{134}, \\
 x_{111}x_{131}x_{224}x_{233}x_{234}x_{434}x_{231}, &
 x_{111}x_{121}x_{123}x_{133}x_{223}x_{424}x_{124}, &
 x_{111}x_{121}x_{123}x_{133}x_{223}x_{424}x_{423}, \\
 x_{111}x_{333}x_{114}x_{224}x_{334}x_{124}x_{134}, &
 x_{111}x_{333}x_{114}x_{224}x_{334}x_{134}x_{234}, &
 x_{111}x_{333}x_{114}x_{224}x_{334}x_{234}x_{214}, \\
 x_{111}x_{333}x_{114}x_{123}x_{124}x_{133}x_{223}, &
 x_{111}x_{333}x_{114}x_{134}x_{224}x_{233}x_{234}, &
 x_{111}x_{333}x_{114}x_{133}x_{124}x_{134}u_{11}, \\
 x_{111}x_{333}x_{114}x_{133}x_{124}x_{223}u_{11}, &
 x_{111}x_{333}x_{114}x_{133}x_{233}x_{134}u_{11}, &
 x_{111}x_{333}x_{114}x_{133}x_{233}x_{223}u_{11}, \\
 x_{111}x_{333}x_{114}x_{224}x_{124}x_{134}u_{11}, &
 x_{111}x_{333}x_{114}x_{224}x_{124}x_{223}u_{11}, &
 x_{111}x_{333}x_{114}x_{224}x_{233}x_{134}u_{11}, \\
 x_{111}x_{333}x_{114}x_{224}x_{233}x_{223}u_{11}, &
 x_{111}x_{333}x_{113}x_{114}x_{223}x_{123}x_{133}, &
 x_{111}x_{333}x_{113}x_{114}x_{223}x_{133}x_{233}, \\
 x_{111}x_{333}x_{113}x_{114}x_{223}x_{233}x_{213}, &
 x_{111}x_{333}x_{114}x_{214}x_{233}x_{213}x_{223}, &
 x_{111}x_{333}x_{114}x_{214}x_{233}x_{223}x_{224}, \\
 x_{111}x_{333}x_{114}x_{214}x_{233}x_{224}x_{234}, &
 x_{111}x_{333}x_{211}x_{233}x_{221}x_{223}x_{224}, &
 x_{111}x_{333}x_{211}x_{233}x_{214}x_{223}x_{224}, \\
 x_{111}x_{333}x_{211}x_{233}x_{213}x_{214}x_{223}, &
 x_{111}x_{333}x_{211}x_{233}x_{214}x_{224}x_{234}, &
 x_{111}x_{333}x_{211}x_{233}x_{224}x_{231}x_{234}, \\
 x_{111}x_{333}x_{211}x_{233}x_{221}x_{224}x_{231}, &
x_{111}x_{333}x_{223}x_{213}x_{211}x_{214}x_{313}, &
x_{111}x_{333}x_{221}x_{224}x_{323}x_{321}x_{311}, \\
x_{111}x_{333}x_{224}x_{323}x_{321}x_{324}x_{311}, &
x_{111}x_{333}x_{224}x_{323}x_{324}x_{311}x_{314}, &
x_{111}x_{333}x_{223}x_{221}x_{224}x_{323}u_{61}, \\
x_{111}x_{333}x_{223}x_{221}x_{224}x_{211}u_{61}, &
x_{111}x_{333}x_{223}x_{224}x_{211}x_{214}u_{61}, &
x_{111}x_{333}x_{223}x_{211}x_{214}x_{313}u_{61}, \\
x_{111}x_{333}x_{221}x_{224}x_{323}x_{311}u_{61}, &
x_{111}x_{333}x_{221}x_{224}x_{211}x_{311}u_{61}, &
x_{111}x_{333}x_{224}x_{211}x_{214}x_{311}u_{61}, \\
x_{111}x_{333}x_{224}x_{214}x_{311}x_{314}u_{61}, &
x_{111}x_{333}x_{224}x_{323}x_{311}x_{314}u_{61}, &
x_{111}x_{333}x_{211}x_{214}x_{313}x_{311}u_{61}, \\
x_{111}x_{333}x_{214}x_{313}x_{311}x_{314}u_{61}, &
x_{111}x_{333}x_{323}x_{313}x_{311}x_{314}u_{61}, &
x_{111}x_{333}x_{223}x_{221}x_{211}v_{132}u_{61}, \\
x_{111}x_{333}x_{223}x_{221}x_{323}v_{132}u_{61}, &
x_{111}x_{333}x_{223}x_{313}x_{211}v_{132}u_{61}, &
x_{111}x_{333}x_{223}x_{313}x_{323}v_{132}u_{61}, \\
x_{111}x_{333}x_{311}x_{221}x_{211}v_{132}u_{61}, &
x_{111}x_{333}x_{311}x_{221}x_{323}v_{132}u_{61}, &
x_{111}x_{333}x_{311}x_{313}x_{211}v_{132}u_{61}, \\
x_{111}x_{333}x_{311}x_{313}x_{323}v_{132}u_{61}, &
 x_{111}x_{333}x_{224}x_{231}x_{234}x_{211}x_{331}, &
 x_{111}x_{333}x_{224}x_{214}x_{334}x_{311}x_{314}, \\
 x_{111}x_{333}x_{224}x_{211}x_{234}x_{214}u_{12}, &
 x_{111}x_{333}x_{224}x_{211}x_{234}x_{331}u_{12}, &
 x_{111}x_{333}x_{224}x_{211}x_{311}x_{214}u_{12}, \\
 x_{111}x_{333}x_{224}x_{211}x_{311}x_{331}u_{12}, &
 x_{111}x_{333}x_{224}x_{334}x_{234}x_{214}u_{12}, &
 x_{111}x_{333}x_{224}x_{334}x_{234}x_{331}u_{12}, \\
 x_{111}x_{333}x_{224}x_{334}x_{311}x_{214}u_{12}, &
 x_{111}x_{333}x_{224}x_{334}x_{311}x_{331}u_{12}, &
 x_{111}x_{333}x_{221}x_{224}x_{331}x_{231}x_{211}, \\
 x_{111}x_{333}x_{221}x_{224}x_{331}x_{211}x_{311}, &
 x_{111}x_{333}x_{221}x_{224}x_{331}x_{311}x_{321}, &
 x_{111}x_{333}x_{224}x_{324}x_{311}x_{321}x_{331}, \\
 x_{111}x_{333}x_{224}x_{324}x_{311}x_{331}x_{334}, &
 x_{111}x_{333}x_{224}x_{324}x_{311}x_{334}x_{314}, &
 x_{111}x_{222}x_{333}x_{122}x_{133}x_{112}x_{132}, \\
 x_{111}x_{222}x_{333}x_{122}x_{133}x_{112}x_{113}, &
 x_{111}x_{222}x_{333}x_{122}x_{133}x_{113}x_{123}, &
 x_{111}x_{222}x_{333}x_{112}x_{133}x_{113}v_{321}, \\
 x_{111}x_{222}x_{333}x_{112}x_{212}x_{113}v_{321}, &
 x_{111}x_{222}x_{333}x_{132}x_{112}x_{133}x_{232}, &
 x_{111}x_{222}x_{333}x_{213}x_{233}x_{212}x_{113},  \\
 x_{111}x_{222}x_{333}x_{112}x_{133}x_{232}v_{321}, &
 x_{111}x_{222}x_{333}x_{233}x_{212}x_{113}v_{321}. &
\end{eqnarray*}
The full list of facets of $(S^{\,2}\times S^{\,2} \times
S^{\,2})_{124}$  may be obtained from these 149 basic facets by
applying the group $S_3\times A_4$. Under this group, the first
145 basic facets form orbits of length 72 each, while each of the
last 4 forms an orbit of length 24, yielding a total of $145
\times 72 + 4 \times 24 = 10536$ facets. The face vector of
$(S^{\,2}\times S^{\,2} \times S^{\,2})_{124}$ is $(124, 1908,
11740, 34140, 50532, 36876, 10536)$. The group $S_3 \times A_4$
appears to be its full group of automorphisms. In Theorem
\ref{T1}, we shall show that this simplicial complex triangulates
$S^{\,2}\times S^{\,2} \times S^{\,2}$.

\bigskip

\noindent {\bf Description of \boldmath{$\CC P^{\,3}_{30}$}\,:}
Consider $(S^{\,2}\times S^{\,2} \times S^{\,2})_{124}$ with the
$S_3$ action given above. Quotienting $(S^{\,2}\times S^{\,2}
\times S^{\,2})_{124}$ by the group $S_3$, we get the
triangulation $\CC P^{\,3}_{30}$. Its vertex-set is $\{x_{ij} ~ :
1\leq i, j \leq 4\} \cup \{y_i ~ : 1\leq i \leq 4\} \cup \{v_i ~ :
1\leq i \leq 4\} \cup \{u_j ~ : 1\leq j \leq 6\}$. (Here $x_{rs} =
q(x_{rrs})$, $u_{t} = q (u_{st})$, $y_{i} = q(x_{jkl})$ and $v_{i}
= q(v_{jkl})$ where $\{i, j, k, l\} = \{1, 2, 3, 4\}$ and $q$ is
the quotient map.) Since the automorphism groups $A_4$ and $S_3$
of $(S^{\,2}\times S^{\,2} \times S^{\,2})_{124}$ commute, its
$A_4$ action induces an $A_4$ action on $\CC P^{\,3}_{30}$. More
explicitly, if $\alpha$, $\beta$ are the generators of the
alternating group $A_4$ given by $\alpha = (1, 2, 3)$, $\beta =
(1,2, 4)$, then $\alpha$, $\beta$ act on the vertices of $\CC
P^{\,3}_{30}$ by\,:
\begin{eqnarray*}
\alpha  & =  & (x_{11},x_{22},x_{33})(x_{12},x_{23},x_{31})
(x_{13}, x_{21},x_{32})(x_{14},x_{24},x_{34})(x_{41},x_{42},
x_{43}) \\
   && ~~~~ (y_{1},y_{2},y_{3})(v_{1},v_{2},v_{3})(u_{1},
   u_{2},u_{4})(u_{3},u_{6},u_{5}),    \\
\beta  &  =  & (x_{11},x_{22},x_{44})(x_{12},x_{24},x_{41}) (x_{21},x_{42},x_{14})
(x_{13},x_{23},x_{43})(x_{31},x_{32},x_{34}) \\
   && ~~~~ (y_{1},y_{2},y_{4})(v_{1},v_{2},v_{4})(u_{1},u_{3},u_{2})(u_{4},u_{5},u_{6}).
\end{eqnarray*}
The following are the basic facets of $\CC P^{\,3}_{30}$ modulo
$A_4 = \langle \alpha, \beta\rangle$\,:
\begin{eqnarray*}
&& \!\!\!\!\!\!\!  x_{11}x_{31}x_{13}x_{34}x_{32}u_{2}u_{1}, ~
  x_{11}x_{31}x_{13}x_{32}v_{4}u_{2}u_{1}, ~
  x_{11}x_{31}x_{13}v_{4}x_{12}u_{2}u_{1}, ~
  x_{11}x_{31}x_{23}x_{34}x_{32}u_{2}u_{1}, \\
&& \!\!\!\!\!\!\!  x_{11}x_{31}x_{23}x_{32}v_{4}u_{2}u_{1}, ~
  x_{11}x_{31}x_{23}v_{4}x_{12}u_{2}u_{1}, ~
  x_{11}x_{21}x_{13}x_{34}x_{32}u_{2}u_{1}, ~
  x_{11}x_{21}x_{13}x_{32}v_{4}u_{2}u_{1}, \\
&& \!\!\!\!\!\!\!  x_{11}x_{21}x_{13}v_{4}x_{12}u_{2}u_{1}, ~
  x_{11}x_{21}x_{23}x_{34}x_{32}u_{2}u_{1}, ~
  x_{11}x_{21}x_{23}x_{32}v_{4}u_{2}u_{1}, ~
  x_{11}x_{21}x_{23}v_{4}x_{12}u_{2}u_{1}, \\
&& \!\!\!\!\!\!\!  x_{11}x_{34}x_{32}x_{21}x_{23}u_{5}u_{1}, ~
  x_{11}x_{34}x_{32}x_{23}v_{1}u_{5}u_{1}, ~
  x_{11}x_{34}x_{32}v_{1}x_{43}u_{5}u_{1}, ~
  x_{11}x_{34}x_{42}x_{21}x_{23}u_{5}u_{1}, \\
&& \!\!\!\!\!\!\!  x_{11}x_{34}x_{42}x_{23}v_{1}u_{5}u_{1}, ~
  x_{11}x_{34}x_{42}v_{1}x_{43}u_{5}u_{1}, ~
  x_{11}x_{24}x_{32}x_{21}x_{23}u_{5}u_{1}, ~
  x_{11}x_{24}x_{32}x_{23}v_{1}u_{5}u_{1}, \\
&& \!\!\!\!\!\!\!  x_{11}x_{24}x_{32}v_{1}x_{43}u_{5}u_{1}, ~
  x_{11}x_{24}x_{42}x_{21}x_{23}u_{5}u_{1}, ~
  x_{11}x_{24}x_{42}x_{23}v_{1}u_{5}u_{1}, ~
  x_{11}x_{24}x_{42}v_{1}x_{43}u_{5}u_{1}, \\
&& \!\!\!\!\!\!\!  x_{11}x_{31}x_{13}x_{42}x_{12}u_{2}u_{1}, ~
  x_{11}x_{31}x_{13}x_{42}x_{34}u_{2}u_{1}, ~
  x_{11}x_{31}x_{13}x_{42}x_{12}y_{3}u_{1}, ~
  x_{11}x_{31}x_{13}x_{42}y_{3}y_{2}u_{1}, \\
&& \!\!\!\!\!\!\!  x_{11}x_{31}x_{13}x_{42}y_{2}x_{43}u_{1},
  x_{11}x_{31}x_{13}x_{42}x_{43}x_{34}u_{1},
  x_{11}x_{13}x_{24}x_{42}x_{21}u_{5}u_{1},
  x_{11}x_{13}x_{24}x_{42}x_{43}u_{5}u_{1}, \\
&& \!\!\!\!\!\!\!  x_{11}x_{13}x_{24}x_{42}x_{21}x_{12}u_{1}, ~
 x_{11}x_{13}x_{24}x_{42}x_{12}y_{3}u_{1}, ~
  x_{11}x_{13}x_{24}x_{42}y_{3}y_{2}u_{1}, ~
  x_{11}x_{13}x_{24}x_{42}y_{2}x_{43}u_{1}, \\
&& \!\!\!\!\!\!\! x_{11}x_{13}x_{31}x_{32}x_{43}y_{2}u_{1},
 x_{11}x_{13}x_{31}x_{32}x_{43}x_{34}u_{1},
 x_{11}x_{12}x_{23}x_{24}x_{42}y_{3}u_{1},
 x_{11}x_{12}x_{23}x_{24}x_{42}x_{21}u_{1}, \\
&& \!\!\!\!\!\!\! x_{11}x_{13}x_{24}x_{32}x_{43}y_{2}u_{1}, ~
  x_{11}x_{12}x_{31}x_{23}x_{42}y_{3}u_{1}, ~
 x_{11}x_{13}x_{32}x_{24}x_{21}u_{5}u_{1}, ~
 x_{11}x_{13}x_{32}x_{24}x_{43}u_{5}u_{1}, \\
 && \!\!\!\!\!\!\! x_{11}x_{13}x_{32}x_{34}x_{21}u_{5}u_{1}, ~
 x_{11}x_{13}x_{32}x_{34}x_{43}u_{5}u_{1}, ~
 x_{11}x_{23}x_{42}x_{31}x_{12}u_{2}u_{1}, ~
 x_{11}x_{23}x_{42}x_{31}x_{34}u_{2}u_{1}, ~ \\
 && \!\!\!\!\!\!\! x_{11}x_{23}x_{42}x_{21}x_{12}u_{2}u_{1}, ~
 x_{11}x_{23}x_{42}x_{21}x_{34}u_{2}u_{1}, ~
 x_{11}x_{13}x_{42}x_{34}x_{21}u_{5}u_{1}, ~
 x_{11}x_{13}x_{42}x_{34}x_{43}u_{5}u_{1}, \\
 && \!\!\!\!\!\!\! x_{11}x_{13}x_{21}x_{42}x_{12}u_{2}u_{1}, ~
 x_{11}x_{13}x_{21}x_{42}x_{34}u_{2}u_{1}, ~
 x_{11}x_{21}y_{3}x_{42}x_{12}x_{13}u_{2}, ~
 x_{11}x_{21}y_{3}x_{42}x_{12}y_{1}u_{2}, \\
 && \!\!\!\!\!\!\! x_{11}x_{21}y_{3}x_{42}x_{34}x_{13}u_{2}, ~
 x_{11}x_{21}y_{3}x_{42}x_{34}y_{1}u_{2}, ~
 x_{11}x_{23}y_{1}x_{42}x_{31}x_{12}u_{2}, ~
 x_{11}x_{23}y_{1}x_{42}x_{31}x_{34}u_{2}, \\
&& \!\!\!\!\!\!\!  x_{11}x_{23}y_{1}x_{42}x_{21}x_{12}u_{2}, ~
 x_{11}x_{23}y_{1}x_{42}x_{21}x_{34}u_{2}, ~
 x_{11}x_{13}y_{2}x_{34}x_{21}y_{4}u_{5}, ~
 x_{11}x_{13}y_{2}x_{34}x_{21}x_{42}u_{5}, \\
&& \!\!\!\!\!\!\! x_{11}x_{13}y_{2}x_{34}x_{43}y_{4}u_{5}, ~
 x_{11}x_{13}y_{2}x_{34}x_{43}x_{42}u_{5}, ~
 x_{11}x_{13}y_{4}x_{32}x_{24}x_{21}u_{5}, ~
 x_{11}x_{13}y_{4}x_{32}x_{24}x_{43}u_{5}, \\
&& \!\!\!\!\!\!\! x_{11}x_{13}y_{4}x_{32}x_{34}x_{21}u_{5}, ~
 x_{11}x_{13}y_{4}x_{32}x_{34}x_{43}u_{5}, ~
 x_{11}x_{13}x_{21}y_{3}x_{42}y_{2}x_{34}, ~
 x_{11}x_{13}x_{24}x_{32}y_{1}x_{43}y_{2}, \\
 && \!\!\!\!\!\!\!x_{11}x_{13}x_{24}x_{32}y_{1}x_{43}y_{4}, ~
 x_{11}x_{12}y_{4}x_{31}x_{23}x_{42}y_{3}, ~
 x_{11}x_{12}y_{4}x_{31}x_{23}x_{42}y_{1}, ~
 x_{11}x_{33}x_{14}x_{24}x_{34}y_{3}y_{2}, \\
&& \!\!\!\!\!\!\! x_{11}x_{33}x_{14}x_{24}x_{34}y_{2}y_{1}, ~
 x_{11}x_{33}x_{14}x_{24}x_{34}y_{1}y_{3}, ~
 x_{11}x_{33}x_{14}y_{4}y_{3}x_{31}x_{23}, ~
 x_{11}x_{33}x_{14}y_{2}x_{24}x_{32}y_{1}, \\
&& \!\!\!\!\!\!\! x_{11}x_{33}x_{14}x_{31}y_{3}y_{2}u_{1}, ~
 x_{11}x_{33}x_{14}x_{31}y_{3}x_{23}u_{1}, ~
 x_{11}x_{33}x_{14}x_{31}x_{32}y_{2}u_{1}, ~
 x_{11}x_{33}x_{14}x_{31}x_{32}x_{23}u_{1}, \\
&& \!\!\!\!\!\!\! x_{11}x_{33}x_{14}x_{24}y_{3}y_{2}u_{1}, ~
 x_{11}x_{33}x_{14}x_{24}y_{3}x_{23}u_{1}, ~
 x_{11}x_{33}x_{14}x_{24}x_{32}y_{2}u_{1}, ~
 x_{11}x_{33}x_{14}x_{24}x_{32}x_{23}u_{1}, \\
&& \!\!\!\!\!\!\! x_{11}x_{33}x_{13}x_{14}x_{23}y_{4}x_{31},
 x_{11}x_{33}x_{13}x_{14}x_{23}x_{31}x_{32},
 x_{11}x_{33}x_{13}x_{14}x_{23}x_{32}y_{4},
 x_{11}x_{33}x_{14}y_{3}x_{32}y_{4}x_{23}, \\
&& \!\!\!\!\!\!\! x_{11}x_{33}x_{14}y_{3}x_{32}x_{23}x_{24},
 x_{11}x_{33}x_{14}y_{3}x_{32}x_{24}y_{1},
 x_{11}x_{33}x_{12}x_{32}x_{21}x_{23}x_{24},
 x_{11}x_{33}x_{12}x_{32}y_{3}x_{23}x_{24}, \\
&& \!\!\!\!\!\!\! x_{11}x_{33}x_{12}x_{32}y_{4}y_{3}x_{23}, ~
 x_{11}x_{33}x_{12}x_{32}y_{3}x_{24}y_{1}, ~
 x_{11}x_{33}x_{12}x_{32}x_{24}y_{4}y_{1}, ~
 x_{11}x_{33}x_{12}x_{32}x_{21}x_{24}y_{4}, \\
 && \!\!\!\!\!\!\! x_{11}x_{33}x_{23}y_{4}x_{12}y_{3}x_{31}, ~
 x_{11}x_{33}x_{21}x_{24}x_{32}y_{4}x_{13}, ~
 x_{11}x_{33}x_{24}x_{32}y_{4}y_{1}x_{13}, ~
 x_{11}x_{33}x_{24}x_{32}y_{1}x_{13}y_{2}, \\
&& \!\!\!\!\!\!\! x_{11}x_{33}x_{23}x_{21}x_{24}x_{32}u_{1},
 x_{11}x_{33}x_{23}x_{21}x_{24}x_{12}u_{1},
x_{11}x_{33}x_{23}x_{24}x_{12}y_{3}u_{1}, ~
 x_{11}x_{33}x_{23}x_{12}y_{3}x_{31}u_{1}, \\
&& \!\!\!\!\!\!\! x_{11}x_{33}x_{21}x_{24}x_{32}x_{13}u_{1},
 x_{11}x_{33}x_{21}x_{24}x_{12}x_{13}u_{1}, ~
 x_{11}x_{33}x_{24}x_{12}y_{3}x_{13}u_{1}, ~
 x_{11}x_{33}x_{24}y_{3}x_{13}y_{2}u_{1}, \\
&& \!\!\!\!\!\!\! x_{11}x_{33}x_{24}x_{32}x_{13}y_{2}u_{1}, ~
 x_{11}x_{33}x_{12}y_{3}x_{31}x_{13}u_{1}, ~
 x_{11}x_{33}y_{3}x_{31}x_{13}y_{2}u_{1}, ~
 x_{11}x_{33}x_{32}x_{31}x_{13}y_{2}u_{1}, \\
&& \!\!\!\!\!\!\! x_{11}x_{33}x_{23}x_{21}x_{12}v_{4}u_{1}, ~
 x_{11}x_{33}x_{23}x_{21}x_{32}v_{4}u_{1}, ~
 x_{11}x_{33}x_{23}x_{31}x_{12}v_{4}u_{1}, ~
 x_{11}x_{33}x_{23}x_{31}x_{32}v_{4}u_{1}, \\
&& \!\!\!\!\!\!\! x_{11}x_{33}x_{13}x_{21}x_{12}v_{4}u_{1}, ~
 x_{11}x_{33}x_{13}x_{21}x_{32}v_{4}u_{1}, ~
 x_{11}x_{33}x_{13}x_{31}x_{12}v_{4}u_{1}, ~
 x_{11}x_{33}x_{13}x_{31}x_{32}v_{4}u_{1}, \\
&& \!\!\!\!\!\!\! x_{11}x_{33}x_{24}y_{4}y_{1}x_{12}x_{31}, ~
 x_{11}x_{33}x_{24}y_{3}x_{34}x_{13}y_{2}, ~
 x_{11}x_{33}x_{24}x_{12}y_{1}y_{3}u_{2}, ~
 x_{11}x_{33}x_{24}x_{12}y_{1}x_{31}u_{2}, \\
&& \!\!\!\!\!\!\! x_{11}x_{33}x_{24}x_{12}x_{13}y_{3}u_{2},
 x_{11}x_{33}x_{24}x_{12}x_{13}x_{31}u_{2}, ~
 x_{11}x_{33}x_{24}x_{34}y_{1}y_{3}u_{2}, ~
 x_{11}x_{33}x_{24}x_{34}y_{1}x_{31}u_{2}, \\
&& \!\!\!\!\!\!\! x_{11}x_{33}x_{24}x_{34}x_{13}y_{3}u_{2},
 x_{11}x_{33}x_{24}x_{34}x_{13}x_{31}u_{2},
 x_{11}x_{33}x_{21}x_{24}x_{31}y_{4}x_{12},
x_{11}x_{33}x_{21}x_{24}x_{31}x_{13}y_{4}, \\
&& \!\!\!\!\!\!\! x_{11}x_{33}x_{21}x_{24}x_{31}x_{12}x_{13},
   x_{11}x_{33}x_{24}y_{1}x_{13}y_{4}x_{31},
 x_{11}x_{33}x_{24}y_{1}x_{13}x_{31}x_{34},
 x_{11}x_{33}x_{24}y_{1}x_{13}x_{34}y_{2}, \\
&& \!\!\!\!\!\!\! x_{11}x_{22}x_{33}x_{21}x_{31}x_{12}y_{4}, ~
 x_{11}x_{22}x_{33}x_{21}x_{31}x_{12}x_{13}, ~
x_{11}x_{22}x_{33}x_{21}x_{31}x_{13}y_{4}, \\
&& \!\!\!\!\!\!\!
x_{11}x_{22}x_{33}x_{12}x_{31}x_{13}v_{4}, ~
 x_{11}x_{22}x_{33}x_{12}x_{21}x_{13}v_{4}, ~
x_{11}x_{22}x_{33}y_{4}x_{12}x_{31}x_{23},\\
 && \!\!\!\!\!\!\!
 x_{11}x_{22}x_{33}y_{4}x_{32}x_{21}x_{13}, ~
x_{11}x_{22}x_{33}x_{12}x_{31}x_{23}v_{4}, ~
 x_{11}x_{22}x_{33}x_{32}x_{21}x_{13}v_{4}.
\end{eqnarray*}

The full list of facets of $\CC P^{\,3}_{30}$ may be obtained from
these 149 basic facets by applying the group $A_4$. Under this
group, the first 145 basic facets form orbits of length 12 each,
while each of the last 4 forms an orbit of length 4, yielding a
total of $145 \times 12 + 4 \times 4 = 1756$ facets. The face
vector of $\CC P^{\,3}_{30}$ is $(30, 362, 2066, 5810, 8470, 6146,
1756)$. The group $A_4$ is the full group of automorphisms of this
simplicial complex (verified by the simpcomp program of
Effenberger and Spreer \cite{es}). We shall see in Theorem
\ref{T1} that it triangulates $\CC P^{\,3}$.

\medskip

Let $S^{\,2}_4 \times  S^{\,2}_4 \times S^{\,2}_4$ denote the
polytopal complex whose polytopes are $A \times B \times C$ as $A,
B, C$ range over the proper faces of a tetrahedron $T$. It is a
subcomplex of the boundary of the nine-dimensional polytope $T
\times T \times T$. Clearly, its geometric carrier is $S^{\,2}
\times  S^{\,2} \times S^{\,2}$. Our main result is\,:

\begin{theo}$\!\!\!${\bf .} \label{T1}
$(a)$ The simplicial complex $(S^{\,2}\times S^{\,2} \times
S^{\,2})_{124}$ is $($the abstract scheme of$\,)$ a simplicial
subdivision of $S^{\,2}_4 \times  S^{\,2}_4
\times S^{\,2}_4$. In consequence, it is a combinatorial
triangulation of $S^{\,2} \times  S^{\,2} \times S^{\,2}$.
\newline $(b)$ The action of $S_3$ on $(S^{\,2}\times S^{\,2}
\times S^{\,2})_{124}$ is good. In consequence, $\CC P^{\,3}_{30}
=  (S^{\,2}\times S^{\,2} \times S^{\,2})_{124}/S_3$ is a
triangulation of \, $\CC P^{\,3}$.
\end{theo}
(The notion of good action of a group on a simplicial complex
is introduced in Definition \ref{good}.)

Since $\CC P^{\,3}_{30}$ is a triangulation of \, $\CC P^{\,3}$,
it follows that $\CC P^{\,3}_{30}$ is a triangulated manifold. By
using the BISTELLAR program of Lutz (\cite{lu2}) we found that
$\CC P^{\,3}_{30}$ is a combinatorial manifold. We applied
BISTELLAR on $\CC P^{\,3}_{30}$ to reduce the number of vertices.
The final output was an 18-vertex simplicial complex which is
bistellar equivalent to $\CC P^{\,3}_{30}$. It is $2$-neighbourly
and its face vector is $(18, 153, 783, 2110, 3021, 2177, 622)$. It
is presented in the appendix. Thus we have\,:

\begin{theo}$\!\!\!${\bf .} \label{T2}
There exists an  $18$-vertex combinatorial triangulation of \,
$\CC P^{\,3}$.
\end{theo}

By \cite{am}, any triangulation of $\CC P^{\,3}$ requires at least
$(3+1)^2+1=17$ vertices. We still do not know if there is a
17-vertex triangulation of $\CC P^{\,3}$. But the output of
BISTELLAR seems to indicate that the face vector of the simplicial
complex $\CC P^{\,3}_{18}$ obtained here is the componentwise
minimum among all triangulations of $\CC P^{\,3}$. It is
noteworthy that while it has been known for the last twentyseven
years that $\CC P^{\,2}$ has a minimal triangulation (with 9
vertices, cf. \cite{am, bd1, kb2}), no explicit triangulation of
$\CC P^{\,3}$ was hitherto known.

\section{Preliminaries and some basic results}

\subsection{Polytopal complexes and their subdivisions}

All affine spaces considered here are finite dimensional. For a
set $A$ in an affine space, the smallest affine subspace
containing $A$ is called the {\em affine span} of $A$ and is
denoted by ${\rm Aff}(A)$. If $S$ is a finite set in an affine
space then the convex hull of $S$ (denoted by $\langle S \rangle$)
is called a {\em polytope}. If $P$ is a polytope in an affine
space $H$ and $E$ is a hyperplane of $H$ (i.e., an affine subspace
of dimension $\dim(H) -1$) such that $P$ is entirely contained in
one of the two closed halfspaces determined by $E$ then the
polytope $F: = P \cap E$ is called a {\em face} of $P$. If $F \neq
P$ then $F$ is called a {\em proper} face of $P$. The 0- and
1-dimensional faces of a polytope are called the {\em vertices} and
{\em edges} of the polytope, respectively. The dimension of a
polytope $P$ is the dimension of its affine span ${\rm Aff}(P)$.
If the dimension of a polytope $P$ is $d$ then we say that $P$ is
a {\em $d$-polytope}. The empty set is a polytope of dimension
$-1$. Clearly, a $d$-polytope has at least $d+1$ vertices. If it
has exactly $d+1$ vertices then it is called a (geometric) {\em
simplex}. If all the proper faces of a polytope $P$ are simplices
then $P$ is called a {\em simplicial} polytope. For a polytope
$P$, $\partial P$ denotes the topological boundary of $P$ and
$\PPP$ denotes the relative interior of $P$. So, $\partial P$ is
the union of all the proper faces of $P$ and $\PPP = P\setminus
\partial P$.

A {\em polytopal complex} $K$ is a finite collection of polytopes
in an affine space such that (i) if $P\in K$, then all the faces
of $P$ are also in $K$ and (ii) the intersection $P \cap Q$ of two
polytopes $P, Q \in K$ is a face of both $P$ and $Q$. The {\em
dimension} of $K$ (denoted by $\dim(K)$) is the largest dimension
of a polytope in $K$. If $K$ is a polytopal complex in an affine
space $H$ then the space $|K| := \cup\{P \, : \, P\in K\}$ (with
subspace topology of $H$) is called the {\em geometric carrier} of
$K$. If all the polytopes in a polytopal complex $K$ are simplices
then $K$ is called a ({\em geometric}) {\em simplicial complex}.

Two affine subspaces $E$ and $F$ in an affine space $H$ are said
to be {\em skew} if $\dim(E + F) = \dim(E) +\dim(F)+1$. For two
polytopes $P$ and $Q$, if ${\rm Aff}(P)$ and ${\rm Aff}(Q)$ are
skew then $\langle P\cup Q\rangle$ is a polytope (denoted by
$P\ast Q$) of dimension $\dim(P) + \dim(Q) +1$ and is called the
{\em join} of $P$ and $Q$. Let $K$ and $L$ be two polytopal
complexes. Suppose ${\rm Aff}(\alpha)$ and ${\rm Aff}(\beta)$ are
skew for all $\alpha \in K$ and $\beta \in L$. Then $K \ast L :=
\{\alpha\ast \beta \, : \, \alpha \in K, \beta \in L\}$ is a
polytopal complex (called the {\em join} of $K$ and $L$) of
dimension $\dim(K) + \dim(L) +1$. Clearly, the join of two
simplicial complexes is a simplicial complex.

If $P$ is a polytope, the collection $\overline{P}$ of all the
faces of $P$ is a polytopal complex, called the {\em face complex}
of $P$; its geometric carrier is $P$ itself. The collection
$\partial \overline{P}$ of all the proper faces of $P$ is another
polytopal complex, called the {\em boundary complex} of $P$; its
geometric carrier is $\partial P$. More generally, if $K$ is a
$d$-dimensional polytopal complex and $|K|$ is a $d$-ball then
consider the subcomplex $\partial K$ whose facets (i.e., maximal
faces) are the $(d-1)$-polytopes $A$ in $K$ such that $A$ is a
face of exactly one $d$-polytope in $K$. Then $|\partial K| =
\partial |K|$. The complex $\partial K$ is called the {\em
boundary} of $K$. If $K_1$, $K_2$ are two polytopal complexes,
then their product is the polytopal complex $K_1 \times K_2 :=
\{A_1 \times A_2 \, : \, A_1\in K_1, A_2\in K_2\}$. Clearly, we
have $|K_1 \times K_2| = |K_1| \times |K_2|$.

A ({\em polytopal}\,) {\em subdivision} of a polytope $P$ is a
polytopal complex $L$ whose geometric carrier is $P$. A {\em
subdivision} $K^{\prime}$ of a polytopal complex $K$ is a
polytopal complex such that $|K^{\prime}| = |K|$ and for each $P
\in K^{\prime}$ there exists $Q \in K$ such that $P \subseteq Q$.
Clearly, if $Q \in K$ then $K^{\prime}[Q] := \{P \in K^{\prime} \,
: \, P \subseteq Q\}$ is a subdivision of $Q$. We say that
$K^{\prime}[Q]$ is the subcomplex of $K^{\prime}$ induced on $Q$.
If a subdivision $K^{\prime}$ of a polytopal complex $K$ is a
geometric simplicial complex then $K^{\prime}$ is called a {\em
simplicial subdivision} of $K$ (cf. \cite[Chapter 5]{zi}). If $M$
is a topological space then a {\em cellulation} of $M$ is a
polytopal complex $K$ such that $|K|$ is homeomorphic to $M$. If,
further, $K$ is a simplicial complex, then $K$ is said to be a
{\em triangulation} of $M$. Clearly, any simplicial subdivision of
a cellulation of $M$ is a triangulation of $M$.

\begin{lemma}$\!\!\!${\bf .} \label{L2.1}
For $d\geq 1$, let $C_d = \sigma \times [0, 1]$, where $\sigma$
is a $(d - 1)$-simplex. Up to isomorphism, there exists a
unique $2d$-vertex simplicial subdivision $\widetilde{C}_d$ of
$C_d$. The facets in $\widetilde{C}_d$ are $a_1\cdots a_ib_i
\cdots b_d$, $1 \leq i \leq d$, where $\sigma = u_1\cdots u_d$,
$a_j = (u_j, 0)$ and $b_j = (u_j, 1)$ for $1\leq j \leq d$.
Moreover, the facets in $\widetilde{C}_d$ are precisely the
maximal cliques in the edge graph of $\widetilde{C}_d$.
\end{lemma}

\noindent {\bf Proof.} We prove the result by induction on $d$. If
$d=1$ then the result is trivial. So, assume that $d >1$ and the
result is true for the $(d-1)$-polytope $C_{d-1}$. Let
$\widetilde{C}_d$ be a $2d$-vertex simplicial subdivision of
$C_d$. Then $\widetilde{C}_d$ induces a simplicial subdivision
$\widetilde{C}_{d-1}$ on $C_{d-1} = u_1 \cdots u_{d-1} \times [0,
1]$. Consider the $(d-1)$-face (simplex) $\alpha = a_1\cdots a_d$
of $C_d$. It is a face of a unique facet $\beta$ in
$\widetilde{C}_d$. Clearly, $\beta$ must be of the form $\alpha
\cup \{b_i\}$ for some $i$. Assume, without loss, that $\beta =
a_1\cdots a_db_d$. Then $C_d$ is the union of $\beta$ and the cone
$C_{d -1}\ast b_d$. (This is actually a special case of Lemma
\ref{L2.2} below.) Since there is no extra vertex, any facet in
$C_{d -1}\ast b_d$ must contain $b_d$. Thus, the set of facets in
$\widetilde{C}_d$ is $\{a_1\cdots a_db_d, \gamma\cup \{b_d\} \, :
\, \gamma$ a facet in $\widetilde{C}_{d- 1}\}$. Since, by
induction hypothesis, the facets in $\widetilde{C}_{d -1}$ are
$a_1\cdots a_ib_i \cdots b_{d-1}$, $1 \leq i \leq d-1$, this
proves the first statement. The second statement is easy to
verify. \hfill $\Box$

\begin{remark}$\!\!\!${\bf .} \label{remark1}
{\rm There is a natural way to obtain a simplicial subdivision of
the product of two (or more) simplices. Namely, if $\Delta_1$,
$\Delta_2$ are simplices, say with vertex-sets $V_1$ and $V_2$,
then impose arbitrary linear orders on $V_1$ and $V_2$, and take
the product partial order $\leq$ on $V := V_1 \times V_2$ (i.e.,
$(x_1, x_2) \leq (y_1, y_2)$ if and only if $x_1 \leq y_1$ and
$x_2 \leq y_2$). Then the chain complex of the poset $(V, \leq)$
is a simplicial subdivision of $\Delta_1 \times \Delta_2$. Lemma
\ref{L2.1} says that - in case one of $\Delta_1$, $\Delta_2$ is
one-dimensional - this is the only way to simplicially subdivide
$\Delta_1 \times \Delta_2$ without adding new vertices. }
\end{remark}

\begin{defn}$\!\!\!${\bf .} \label{antistar}
{\rm Let $P$ be a polytope and $x$ be a point of $P$ ($x$ may or
may not be a vertex of $P$). Let $C$ be the smallest face of $P$
containing $x$. Then the polytopal complex ${\rm Ast}_P(x)$,
consisting of all the faces of $P$ not containing $C$, is called
the {\em antistar} of $x$ in $P$. }
\end{defn}

\begin{lemma} {\rm \bf (The antistar Lemma).} \label{L2.2}
Let $P$ be a $d$-polytope and $x$ be a point in $P$.  Then the
collection $\widetilde{P} := \{x\ast D \, : D \in {\rm
Ast}_P(x)\}$ is a polytopal subdivision of $P$. Moreover, if ${\rm
Ast}_P(x)^{\prime}$ is a simplicial subdivision of ${\rm
Ast}_P(x)$ then $\{x\ast \sigma \, : \, \sigma \in {\rm
Ast}_P(x)^{\prime}\}$ is a simplicial subdivision of $P$.
\end{lemma}

\noindent {\bf Proof.} Clearly, $\widetilde{P}$ is a polytopal
complex. Let $|\widetilde{P}|$ be the geometric carrier of
$\widetilde{P}$. Since any polytope in $\widetilde{P}$ is of the
form $x\ast D \subseteq P$ for some polytope $D$ in ${\rm
Ast}_P(x)$, it follows that $|\widetilde{P}| \subseteq P$.

We now prove $P\subseteq |\widetilde{P}|$. If $d\leq 2$ the result
is obvious. So, assume $d\geq 3$ and the result is true for
smaller values of $d$. Let $C$ be the smallest face of $P$ not
containing $x$. Let $x\neq y \in P$. If $y \in |{\rm Ast}_P(x)|$
then $y\in |\widetilde{P}|$. If $y\in \partial P \setminus |{\rm
Ast}_P(x)|$ then $y \in F$ for some $(d-1)$-face $F$ of $P$ which
contains $C$. By induction hypothesis, $y \in |\widetilde{F}|
\subseteq |\widetilde{P}|$. So, assume that $y\in P \setminus
\partial P$. Let $R$ be the ray containing $y$ and initial point
$x$. Then $R \cap P$ is a line segment of the form $[x, z]$ for
some $z\in \partial P$. If $z \in D$ for some proper face $D$ of
$P$ containing $x$ then $[x, z] \subseteq D$ and hence $y\in D
\subseteq \partial P$, a contradiction. So, $z \in E$ for some $E
\in {\rm Ast}_P(x)$. This implies that $y \in x \ast E \subseteq
|\widetilde{P}|$. This completes the proof of the first statement.
The second statement is now obvious. \hfill $\Box$

\bigskip

We need the following technical lemma in the proofs of Lemmas \ref{R5}
and \ref{R7}.

\begin{lemma}$\!\!\!${\bf .} \label{L2.3}
Let $\alpha$ $\beta$ and $\gamma$ be three simplices $($in an
affine space$)$ of dimensions $i\geq 1$, $j\geq 1$ and $k\geq -1$
respectively. Suppose the relative interiors of $\alpha$ and
$\beta$ have a non-empty intersection, and there is a $(j-1)$-face
$\beta_0$ of $\beta$ for which the convex hull $\langle \gamma
\cup \alpha \cup \beta_0\rangle$ is an $(i+j+k+1)$-simplex. Then,
both $\overline{\gamma}\ast\overline{\alpha}\ast\partial
\overline{\beta}$ and $\overline{\gamma} \ast \partial
\overline{\alpha} \ast \overline{\beta}$ are  simplicial
subdivisions of the convex hull $\langle \gamma \cup \alpha \cup
\beta \rangle$.
\end{lemma}

\noindent {\bf Proof.} The hypothesis on dimension means that the
affine span of each of the three polytopes $\alpha$, $\beta_0$,
$\gamma$ is skew from the affine span of the union of the other
two. In particular, ${\rm Aff}(\alpha)$ and ${\rm Aff}(\beta_0)$
are skew. Since $\beta_0$ is of co-dimension 1 in $\beta$, it
follows that the intersection of ${\rm Aff}(\alpha)$ and ${\rm
Aff}(\beta)$ is of dimension (at most, hence exactly) 0, i.e.,
this intersection is a point $a$. Therefore, the intersection of
the relative interiors $\intalpha$ and $\intbeta$ is $\{a\}$. Thus
${\rm Aff}(\gamma \cup \alpha)$ is disjoint from $\partial \beta$,
and hence from each facet $\beta_1$ of $\partial
\overline{\beta}$. Thus the hypothesis of the lemma holds for each
facet $\beta_1$ of $\partial \overline{\beta}$ (in place of
$\beta_0$). Since the join of two simplices is a simplex, it
follows that $\overline{\gamma} \ast \overline{\alpha} \ast
\partial \overline{\beta}$ is a geometric simplicial complex.
Clearly, $\langle\gamma \cup \alpha \cup \beta_1 \rangle \subseteq
\langle \gamma \cup \alpha \cup \beta \rangle$ for each facet
$\beta_1$ of $\partial \overline{\beta}$ and hence
$|\overline{\gamma} \ast \overline{\alpha} \ast \partial
\overline{\beta}| \subseteq \langle \gamma \cup \alpha \cup \beta
\rangle$. For the opposite inclusion, let $x \in \langle \gamma
\cup \alpha \cup \beta \rangle$. Then $x = t_1x_1 + t_2x_2 +
t_3x_3$ for some $x_1 \in \gamma$, $x_2\in \alpha$, $x_3 \in
\beta$, $0\leq t_1, t_2, t_3 \leq 1$ and $t_1+t_2+t_3 = 1$. Since
$a \in \intbeta$, $\beta = \cup\{a \ast \lambda \, : \, \lambda$
facet of $\partial \overline{\beta}\}$. Thus, there exists a facet
$\lambda$ of $\partial \overline{\beta}$ such that $x_3 \in a \ast
\lambda$. Accordingly, let $x_3 = (1-s)a + sy_3$ for some $y_3 \in
\lambda$ and $0\leq s \leq 1$. Then $y_2 := \frac{t_2}{t_2+
t_3(1-s)}x_2 + \frac{t_3(1 -s)}{t_2+t_3(1-s)}a \in \alpha$ and $x
= t_1x_1 + (t_2+t_3(1-s))y_2 + (st_3)y_3 =  t_1x_1 + (1 - t_1 -
st_3)y_2 + (st_3)y_3 \in \gamma \ast \alpha \ast \lambda \subseteq
|\overline{\gamma} \ast \overline{\alpha} \ast \partial
\overline{\beta}|$. Thus, $|\overline{\gamma} \ast
\overline{\alpha} \ast \partial \overline{\beta}| = \langle \gamma
\cup \alpha \cup \beta\rangle$.

Since ${\rm Aff}(\alpha) \cap {\rm Aff}(\beta) = \{a\} = \intalpha
\cap \intbeta$ and $\dim({\rm Aff}(\alpha\cup\beta\cup\gamma)) =
i+j+1$, it follows that $\dim({\rm Aff}(\alpha_1\cup\beta\cup
\gamma)) = i+j+1$ for any facet $\alpha_1$ of $\partial
\overline{\alpha}$. Then, by the same argument, $\overline{\gamma}
\ast \partial \overline{\alpha} \ast \overline{\beta}$ is a
geometric simplicial complex and $|\overline{\gamma} \ast \partial
\overline{\alpha} \ast \overline{\beta}| = \langle \gamma \cup
\alpha \cup \beta \rangle$. This completes the proof. \hfill
$\Box$

\subsection{Simplicial complexes and their geometric realization}

An {\em abstract simplicial complex} (or simply a {\em simplicial
complex}) is a finite collection of finite sets such that every
subset of an element is also an element. For $i \geq 0$, an
element of size $i + 1$ is called a {\em face} of dimension $i$
(or an {\em $i$-face)} of the complex. By convention, the empty
set is a face of dimension $-1$. The {\em dimension} of a
simplicial complex $X$ (denoted by $\dim(X)$) is the maximum of
the dimensions of its faces. The 0- and 1-dimensional faces of a
complex are also called the {\em vertices} and {\em edges} of the
complex respectively. The set $V(X)$ of vertices of a simplicial
complex $X$ is called the {\em vertex-set} of $X$. A maximal face
in a simplicial complex is also called a {\em facet} of the
complex. Since a simplicial complex is uniquely determined by the
set of facets, we sometimes identify a simplicial complex with its
set of facets.

If $K$, $L$ are two simplicial complexes, then a {\em simplicial
isomorphism} from $K$ to $L$ is a bijection $\pi : V(K) \to V(L)$
such that for $\sigma \subseteq V(K)$, $\sigma$ is a face of $K$
if and only if $\pi (\sigma)$ is a face of $L$. The complexes $K$,
$L$ are called {\em isomorphic} when such an isomorphism exists.
We identify two simplicial complexes if they are isomorphic.

The {\em face vector} of a $d$-dimensional simplicial complex is
the vector $(f_0, f_1, \dots, f_d)$, where $f_i$ is the number of
$i$-faces in the complex. It is called {\em $k$-neighbourly} if
$f_{k-1} = {f_0 \choose k}$, i.e., if all the possible
$(k-1)$-faces occur in the complex.

If $K$ is a geometric simplicial complex and $V$ is the set of
vertices of $K$ then ${\cal K} := \{A \subseteq V \, : \langle A
\rangle \in K\}$ is a simplicial complex, called the {\em abstract
scheme} of $K$.

If $X$ is a $d$-dimensional simplicial complex, then let us
identify the vertex-set $V$ of $X$ with a set $V$ of points in
$\RR^{\,2d+1}$ such that any subset of $V$ of size at most $2d+2$
is affinely independent. Then the geometric simplicial complex
${\cal X}:=\{ \langle \sigma \rangle \, : \, \sigma\in X\}$ is
called a {\em geometric realization} of $X$. The geometric carrier
$|{\cal X}|$ of ${\cal X}$ is also called the {\em geometric
carrier} of $X$ and is denoted by $|X|$. Clearly, $X$ is
isomorphic to the abstract scheme of ${\cal X}$. If a topological
space $M$ is homeomorphic to $|X|$ then say $X$ {\em triangulates}
$M$. It is trivial that isomorphic finite simplicial complexes
have homeomorphic geometric carriers. If $|X|$ is a topological
manifold (respectively, $d$-sphere) then $X$ is called a {\em
triangulated manifold} (resp. {\em triangulated $d$-sphere}). If
$|X|$ is a pl manifold (with the pl structure induced by $X$) then
$X$ is called a {\em combinatorial manifold}.

For a finite set $V$ with $d+2$ ($d\geq -1$) elements, the
collection $\overline{V}$ of all the subsets of $V$ is a
simplicial complex which triangulates the $(d+1)$-ball. The
complex $\overline{V}$ is called the {\em standard} $(d+1)$-ball
and is also denoted by $D^{d+1}_{d+2}(V)$ (or simply by
$D^{d+1}_{d+2}$). It is the abstract scheme of the face complex of
a $(d+1)$-simplex. The collection $\partial \overline{V}$ of all
the proper subsets of $V$ is a simplicial complex (a subcomplex of
$\overline{V}$) which triangulates the $d$-sphere $S^{\,d}$. The
complex $\partial \overline{V}$ is called the {\em standard
$d$-sphere} and is also denoted by $S^{\,d}_{d+2}(V)$ (or simply
by $S^{\,d}_{d+2}$). (Generally, we write $X = X^d_n$ to indicate
that $X$ has $n$ vertices and dimension $d$.) Note that
$S^{\,d}_{d+2}$ is the abstract scheme of the boundary complex of
a $(d+1)$-simplex.

If $\sigma$ is a face of a simplicial complex $X$ then the {\em
link} of $\sigma$ in $X$, denoted by ${\rm lk}_X(\sigma)$, is the
simplicial complex whose faces are the faces $\tau$ of $X$ such
that $\tau \cap \sigma = \emptyset$ and $\sigma\cup\tau$ is a face
of $X$. The number of vertices in the link of $\sigma$ is called
the {\em degree} of $\sigma$.  For $1\leq d\leq 4$, $X$ is a
combinatorial $d$-manifold if and only if the vertex links are
triangulated $(d-1)$-spheres. The {\em star} of $\sigma$, denoted
by ${\rm star}_X(\sigma)$ or ${\rm star}(\sigma)$, is the
subcomplex $\bar{\sigma} \ast {\rm lk}_X(\sigma)$ of $X$. The {\em
anti-star} of $\sigma$ in $X$ is the simplicial complex ${\rm
ast}_X(\sigma) =\{\tau \in X ~ : ~ \tau \cap \sigma =
\emptyset\}$.

Let $X$ be a simplicial complex whose maximal faces are
$d$-dimensional. Let $A$ be a $(d-i)$-face ($0 \leq i \leq d$) of
$X$ whose link in $X$ is the standard $(i-1)$-sphere $\partial
\overline{B}$. Suppose also that $B \not\in X$. Then the
simplicial complex $\widehat{X} := (X \setminus (\overline{A} \ast
\partial \overline{B})) \cup (\partial \overline{A} \ast
\overline{B})$ and $X$ triangulate the same topological space. The
complex $\widehat{X}$ is said to be obtained from $X$ by a {\em
bistellar $i$-move}. Two simplicial complexes $X$ and $Y$ are
called {\em bistellar equivalent} if there exists a sequence $X_1,
\dots, X_k$ of simplicial complexes such that $X=X_1$, $Y=X_k$ and
$X_{i+1}$ is obtained from $X_i$ by a bistellar move for $1\leq i
\leq k-1$.

\subsection{Group actions and quotients}

Consider $(S^{\,2})^d$, the Cartesian product of $d$ copies of the
2-sphere $S^{\,2}$. The symmetric group $S_d$ acts on
$(S^{\,2})^d$, by co-ordinate permutations, as a group of
homeomorphisms\,:
$$
\pi : (x_1, \dots, x_d) \mapsto (x_{\pi(1)}, \dots, x_{\pi(d)}),
\mbox{ for } \pi\in S_d.
$$
This action is far from free. Indeed, the diagonal $S^{\,2}$ is
pointwise fixed under $S_d$. Yet, miraculously, the quotient space
$(S^{\,2})^d/S_d$ (with the quotient topology) is a manifold.
(This result is well known to algebraic geometers.)

\begin{lemma}$\!\!\!${\bf .} \label{L2.4}
The quotient space $(S^{\,2})^d/S_d$ is homeomorphic to the
$d$-dimensional complex projective space $\CC P^{\,d}$.
\end{lemma}

\noindent {\bf Proof.} We use the usual identification of
$S^{\,2}$ with $\CC P^{1}$ (``the Riemann sphere") via
stereographic projection. So, we need to show that $(\CC
P^{1})^d/S_d \cong \CC P^{\,d}$. Recall the usual description of
$\CC P^{\,d}$ by homogeneous co-ordinates\,: $\CC P^{\,d} =
(\CC^{d+1}\setminus\{(0, \dots, 0)\})/\!\sim$, where the binary
relation $\sim$ is defined by $(y_0, \dots, y_{d}) \sim (x_0,
\dots, x_{d})$ if and only if $y_i = \lambda x_i$, $0\leq i\leq
d$, for some $\lambda \in \CC\setminus\{0\}$. As usual, $[x_0,
x_1,  \dots, x_d]$ denotes the $\sim$-class containing $(x_0,
\dots, x_{d})$.  Now consider the map $\varphi \colon (\CC
P^{\,1})^d \to \CC P^{\,d}$ defined as
$$
\varphi([z_1, w_1], \dots, [z_d, w_d]) = [\alpha_0, \alpha_1,
\dots, \alpha_d],
$$
where the complex numbers $\alpha_0, \alpha_1, \dots, \alpha_d$
are determined by the identity
$$
\sum_{j=0}^{d} \alpha_j X^j Y^{d-j} = \prod_{k=1}^d (z_kX + w_kY).
$$
Since multiplication in $\CC$ is commutative, $\varphi$ is
well-defined. It is clearly continuous. Since the field $\CC$ is
algebraically closed, each homogeneous polynomial of degree $d$
($\neq 0$) in $\CC[X, Y]$ has exactly $d$ roots in $\CC P^{1}$
(counting with multiplicity). Also, modulo multiplication by
non-zero scalars, such a polynomial is uniquely determined by its
zeros, and in turn it determines its zero-set except for
permutations. Thus, $\varphi$ is onto and its fibres are precisely
the $S_d$-orbits in its domain. Therefore, $\varphi$ induces a
continuous bijection $\widehat{\varphi} \colon (\CC P^{1})^d/S_d
\to \CC P^{\,d}$. Since both domain and range are compact
Hausdorff spaces, $\widehat{\varphi}$ is a homeomorphism. \hfill
$\Box$

\medskip

One says that $\CC P^{\,d}$ is the $d$-th symmetric power of
$S^{\,2} = \CC P^1$. In \cite{bd11}, we used the case $d=2$ of
this lemma to construct a
10-vertex triangulation of $\CC P^{\,2}$. Here we use its $d=3$
case to obtain an explicit triangulation of $\CC P^{\,3}$.

\begin{defn}$\!\!\!${\bf .} \label{pure}
{\rm Let $G$ be a group of simplicial automorphisms of a
simplicial complex $X$. We shall say that the action of $G$ on $X$
is {\em pure} if it satisfies\,: $(a)$ whenever $u$, $v$ are
distinct vertices from the same $G$-orbit, $uv$ is a non-edge of
$X$, and $(b)$ for each $G$-orbit $\theta \subseteq V(X)$ and each
$\alpha \in X$, the stabiliser $G_{\alpha}$ of $\alpha$ in $G$
acts transitively on $\theta \cap V({\rm lk}_X(\alpha))$.}
\end{defn}

The ``if\," part of the following lemma is from \cite{bd11}. This
lemma clarifies the significance of the notion of purity.

\begin{lemma}  {\rm \bf (The purity Lemma).} \label{L2.5}
Let $G$ be a group of simplicial automorphisms of a simplicial
complex $X$. Let $q \colon V(X) \to V(X)/G$ denote the quotient
map, and $X/G := \{q(\alpha) \, : \, \alpha\in X\}$. Extend the
action of $\,G$ on $V(X)$ to an action of $\,G$ on $|X|$ piecewise
linearly, i.e., affinely on the geometric carrier of each face.
Then there is a homeomorphism from $|X|/G$ onto $|X/G|$ induced by
$|q| \colon |X| \to |X/G|$ if and only if the action of $G$ on $X$
is pure. In consequence, when the $G$-action on $X$ is pure, we
have $|X/G| = |X|/G$.
\end{lemma}

\noindent {\bf Proof.} First suppose that the $G$-action is pure.
The condition $(a)$ ensures that the quotient map $q$ is one-one
on each face of $X$. The simplicial map $q \colon X \to X/G$
induces a piecewise linear continuous map $|q|$ from $|X|$ onto
$|X/G|$.

\smallskip

\noindent {\sf Claim.} {\em The fibres of $\,q \colon X \to X/G$
are precisely the $G$-orbits on simplices of $X$ $($that is, if
$\alpha$, $\alpha^{\,\prime} \in X$ are such that $q(\alpha) =
q(\alpha^{\,\prime})$ then there exists $g \in G$ such that
$g(\alpha) = \alpha^{\,\prime})$.}

\smallskip

We prove the claim by induction on $k = \dim(\alpha) =
\dim(\alpha^{\,\prime})$. The claim is trivial for $k = -1$. So,
assume $k \geq 0$ and the claim is true for all smaller
dimensions. Choose a face $\beta \subseteq \alpha$ of dimension
$k-1$, and let $\beta^{\,\prime}\subseteq \alpha^{\,\prime}$ be
such that $q(\beta^{\,\prime}) = q(\beta)$. By induction
hypothesis, $\beta^{\,\prime}$ and $\beta$ are in the same
$G$-orbit. Therefore, applying a suitable element of $G$, we may
assume, without loss of generality, that $\beta^{\,\prime} =
\beta$. Let $\alpha = \beta \cup \{x\}$, $\alpha^{\,\prime} =
\beta \cup \{x^{\,\prime}\}$. Then $q(x) = q(x^{\,\prime})$, i.e.,
$x$ and $x^{\,\prime}$ are in the same $G$-orbit. Now, by
assumption $(b)$, there is a $g\in G_{\beta}$ such that $g(x)=
x^{\,\prime}$. Then $g(\alpha) =\alpha^{\,\prime}$. This proves
the claim.

\smallskip

In the presence of condition $(a)$, the claim ensures that the
fibres of $|q|$ are precisely the $G$-orbits on points of $|X|$.
Hence $|q|$ induces the required homeomorphism between $|X|/G$ and
$|X/G|$.

Now assume that the map between $|X|/G$ and $|X/G|$ induced by
$|q|$ is a homeomorphism. Clearly, condition $(a)$ holds since the
quotient map preserves dimension. So, $(a)$ holds, and the fibres
of $q$ are the $G$-orbits in $X$. Let $\alpha$, $\theta$ be as in
$(b)$, and let $x\neq y$ be vertices in $\theta \cap {\rm
lk}(\alpha)$. Then $q(\alpha \cup \{x\}) = q(\alpha \cup \{y\})$
and hence there is an element $g$ of $G$ such that $\alpha
\cup\{y\} = g(\alpha \cup \{x\})$. Since $y \in \theta$, $(a)$
ensures that $\theta$ is disjoint from the vertex-set of $\alpha$,
so that we must have $y = g(x)$. Hence $g \in G_{\alpha}$, and
condition $(b)$ holds. Thus the $G$-action on $X$ is pure. \hfill
$\Box$

\begin{defn}$\!\!\!${\bf .} \label{good}
{\rm Let $G$ be a group of simplicial automorphisms of a
simplicial complex $X$. We shall say that the group $G$ is {\em
good} (or, the action of $G$ on $X$ is {\em good}\,) if it
satisfies\,: whenever $x\neq y$ are vertices of $X$ from a common
$G$-orbit, then $xy$ is a non-edge, and there is $g\in G$ such
that $g(x) = y$ and $g$ fixes all the vertices in ${\rm lk}(x)
\cap {\rm lk}(y)$.}
\end{defn}

\begin{coro}$\!\!\!${\bf .} \label{L2.6}
Let $G$ be a group of simplicial automorphisms of a simplicial
complex $X$. If $G$ is good then the action of $G$ is pure and
hence $|X/G| = |X|/G$.
\end{coro}

\noindent {\bf Proof.} Clearly $G$ satisfies condition (a). We
need to show that $G$ also satisfies condition $(b)$ for purity.
So, let $\alpha$, $\theta$ be as in $(b)$ and $x\neq y$ be
vertices in $\theta \cap {\rm lk}(\alpha)$. Then the element $g$
of $G$ (given by Definition \ref{good}) fixes all the vertices of
$\alpha$, and hence $g \in G_{\alpha}$. \hfill $\Box$

\begin{remark}$\!\!\!${\bf .} \label{remark2}
{\rm Goodness of $G$ may not be necessary for the purity of its
action. However, goodness has the virtue of simplicity\,: to
verify it one need only examine the action of $G$ on the
1-skeleton of $X$. Notice also that purity of $G$ implies that the
stabiliser of any face fixes the face vertex-wise, so that the
definition of a good action is not as stringent as it may appear.
Note that if $G$ is a simplicial automorphism group of $X$ such
that any two distinct vertices from each $G$-orbit are at distance
at least 3 in the 1-skeleton of $X$, then $G$ is good. This strong
requirement is often employed in the literature to ensure that the
quotient of a triangulated manifold $X$ under such a group $G$ is
again a triangulated manifold.}
\end{remark}

In the following, we shall say that an isometry group $G$ of a
geometric simplicial complex is {\em good} (or {\em pure}) if the
induced action of $G$ on its abstract scheme is good
(respectively, pure).

\section{Proofs}

By a slight abuse of notation, we let $S^{\,2}_4$ denote the
boundary complex of the regular tetrahedron in $\RR^3$ whose
vertices are $x_1 = (1, -1, -1)$, $x_2 = (-1, 1, -1)$, $x_3 = (-1,
-1, 1)$ and $x_4 = (1, 1, 1)$. We use the notation $x_{ijk}$ to
denote the vertex $(x_i, x_j, x_k)$ of $S^{\,2}_4 \times S^{\,2}_4
\times S^{\,2}_4$ in $\RR^9$. So, the vertex-set of $S^{\,2}_4
\times S^{\,2}_4 \times S^{\,2}_4$ is $\{x_{ijk} ~ : 1\leq i, j,k
\leq 4\}$. (With our notation, for instance $x_{214} = (-1, 1, -1,
1, -1, -1, 1, 1, 1)$.) Observe that $S_3$ acts on these vertices
by permuting the positions of the three subscripts and this action
induces the same action on $|S^{\,2}_4| \times |S^{\,2}_4| \times
|S^{\,2}_4|$ as was used in Lemma \ref{L2.4}. $S_3$ is an isometry
group of $S^{\,2}_4 \times S^{\,2}_4 \times S^{\,2}_4$.

\bigskip

\noindent {\bf A sketch of the proof of Theorem \ref{T1}.} Let $W$
denote the product polytopal complex $S^{\,2}_4 \times S^{\,2}_4
\times S^{\,2}_4$ in $\RR^{\,9}$. For $0 \leq k \leq 6$, let $W^k$
be the $k$-skeleton of $W$. That is, $W^k$ is the subcomplex of
$W$ consisting of all its polytopes of dimension at most $k$. For
$0 \leq k \leq 6$, we construct a simplicial subdivision $X^k$ of
$W^k$. The construction is by finite induction on $k$, with the
following constraints\,: (i) $X^0 = W^0$, $X^1 = W^1$, (ii)
$X^{k-1} \subseteq X^k$ for $k\geq 1$, (iii) $X^k$ is stable under
the isometry group $S_3$ of $W$, (iv) the $S_3$-action on $X^k$ is
good. It turns out that there are three choices for $X^2$ subject
to these constraints. We choose the most symmetrical of these
three choices, with an isometry group $S_3 \times A_4$.

To keep the construction under control, we also require that (v)
for $k \geq  3$, $X^k$ inherits the isometry group $S_3 \times A_4$
from $X^{k-1}$.

To get $X^3$, we need to add 60 new vertices, namely $u_{ij}$,
$1\leq i, j \leq 6$, $v_{pqr}$, $1\leq p, q, r \leq 4$ ($p, q, r$
distinct). So, the vertex-set of $X^3$ is
\begin{eqnarray}
V_{124} := \{x_{ijk} : 1\leq i, j, k \leq 4\} \cup \{u_{ij} :
1\leq i, j \leq 6\}  \cup \{v_{pqr} : 1\leq p, q, r \leq 4 \,
\mbox{are distinct}\}.
\end{eqnarray}

At the end of this induction, we obtain the simplicial subdivision
$X = X^6$ of $W$ with the vertex-set $V_{124}$ and the following
set of edges.
\begin{eqnarray}
E(X) \!\!\! & = \!\!\!& (S_3\times A_4)(\{x_{111}y \, : \, x_{111}
\neq y \in V_{124}\} \cup \{u_{11}u_{22},
u_{11}v_{312}, v_{321}x_{112}, v_{321}x_{113}\} \cup \nonumber \\
&& \hspace{8mm} \{u_{11}x_{114}, u_{11}x_{221}, u_{11}x_{223},
u_{11}x_{133}, u_{11}x_{121}, u_{11}x_{131},  u_{11}x_{124}\}
 \cup \{x_{112}x_{122},   \nonumber \\
&& \hspace{12mm}
x_{112}x_{133}, x_{112}x_{144}, x_{112}x_{322}, x_{112}x_{422},
x_{112}x_{113}, x_{112}x_{114},  x_{112}x_{332},
\nonumber \\
&& \hspace{16mm}
x_{112}x_{442}, x_{112}x_{132}, x_{112}x_{142}, x_{112}x_{134},
x_{112}x_{342},
x_{112}x_{334}, x_{112}x_{344}\}). \label{edges}
\end{eqnarray}
The abstract scheme of $X$ is the 124-vertex simplicial complex
$(S^{\,2}\times S^{\,2} \times S^{\,2})_{124}$, which was promised
in Theorem \ref{T1}.

Now let $x \neq y$ be vertices of $X$ from a common $S_3$-orbit.
Since $A_4$ commutes with $S_3$, if $\{x, y\}$ satisfies the
requirements for the goodness of $S_3$, then so does $\{g(x),
g(y)\}$ for all $g \in S_3 \times A_4$. Thus, up to the action of
$S_3 \times A_4$, we can assume that $\{x, y\} \in \{\{x_{123},
x_{213}\}$, $\{x_{123}, x_{312}\}$, $\{x_{112}, x_{121}\}$,
$\{u_{11}, u_{21}\}$, $\{u_{11}, u_{41}\}$, $\{u_{11}, u_{51}\}$,
$\{v_{123}, v_{213}\}$, $\{v_{123}, v_{312}\}\}$. From
(\ref{edges}), one sees that none of these eight pairs forms an
edge. This verifies the first condition in the definition of a
good action.

From (\ref{edges}), we get $V({\rm lk}_X(x_{123})) = \{x_{111}$,
$x_{222}$, $x_{333}$, $x_{444}$, $x_{113}$, $x_{223}$, $x_{121}$,
$x_{323}$, $x_{122}$, $x_{133}$, $x_{144}$, $x_{424}$, $x_{443}$,
$x_{114}$, $x_{422}$, $x_{343}, u_{66}, u_{45}, u_{23}\}$.
Therefore, $V({\rm lk}_X(x_{123})) \cap V({\rm lk}_X(x_{213})) =
\{x_{111}$, $x_{222}$, $x_{333}$, $x_{444}$, $x_{113}$, $x_{223}$,
$x_{443}$, $x_{114}\}$. Then $(1,2)(x_{123}) = x_{213}$ and
$(1,2)$ fixes all the vertices in  $V({\rm lk}_X(x_{123})) \cap
V({\rm lk}_X(x_{213}))$.
Again, $V({\rm lk}_X(x_{123})) \cap V({\rm lk}_X(x_{312})) =
\{x_{111}$, $x_{222}$, $x_{333}$, $x_{444}\}$. Thus, $(1,2,3)
(x_{123}) = x_{312}$ and $(1,2,3)$ fixes all the vertices in
$V({\rm lk}_X(x_{123})) \cap V({\rm lk}_X(x_{312}))$.

From (\ref{edges}), $V({\rm lk}_X(x_{112})) = \{x_{111}$,
$x_{222}$, $x_{333}$, $x_{444}$, $x_{122}$, $x_{212}$, $x_{113}$,
$x_{133}$, $x_{313}$, $x_{114}$, $x_{144}$, $x_{414}$, $x_{232}$,
$x_{322}$, $x_{332}$, $x_{242}$, $x_{422}$, $x_{442}$, $x_{334}$,
$x_{344}$, $x_{434}$, $x_{132}$, $x_{312}$, $x_{142}$, $x_{412}$,
$x_{314}$, $x_{134}$, $x_{432}$, $x_{342}$, $u_{21}$, $u_{51}$,
$u_{32}$, $u_{42}$, $u_{33}$, $u_{43}$, $u_{14}$, $u_{64}$,
$u_{25}$, $u_{55}$, $u_{36}$, $u_{46}$, $v_{231}$, $v_{321}$,
$v_{241}$, $v_{421}\}$. Therefore, $V({\rm lk}_X(x_{112})) \cap
V({\rm lk}_X(x_{121})) = \{x_{111}$, $x_{222}$, $x_{333}$,
$x_{444}$, $x_{122}$, $x_{133}$, $x_{144}$, $x_{322}$, $x_{422},
x_{344}\}$. Then $(2,3)(x_{112}) = x_{121}$ and $(2,3)$ fixes all
the vertices in $V({\rm lk}_X(x_{112}))\cap V({\rm
lk}_X(x_{121}))$.

From (\ref{edges}), $V({\rm lk}_X(u_{11})) =  \{x_{111}$,
$x_{222}$, $x_{333}$, $x_{444}$, $x_{114}$, $x_{221}$, $x_{223}$,
$x_{133}$, $x_{121}$, $x_{131}$, $x_{144}$, $x_{433}$, $x_{233}$,
$x_{224}$, $x_{434}$, $x_{424}$, $x_{124}$, $x_{134}$, $u_{22}$,
$u_{53}$, $u_{34}$, $u_{65}$, $v_{312}$, $v_{342}\}$. So, $V({\rm
lk}_X(u_{11})) \cap V({\rm lk}_X(u_{41})) = \{x_{111}$, $x_{222}$,
$x_{333}$, $x_{444}$, $x_{121}$, $x_{131}$, $x_{434}$,
$x_{424}\}$. Then $(1,3)(u_{11}) = u_{41}$ and $(1,3)$ fixes all
the vertices in  $V({\rm lk}_X(u_{11})) \cap V({\rm
lk}_X(u_{41}))$. Similarly, $(2,3)(u_{11}) = u_{51}$ and $(2,3)$
fixes all the vertices in $V({\rm lk}_X(u_{11}))\cap V({\rm
lk}_X(u_{51}))$. Observe that $V({\rm lk}_X(u_{11}))\cap V({\rm
lk}_X(u_{21})) = \{x_{111}$, $x_{222}$, $x_{333}$, $x_{444}\}$.
Thus, $(1,2,3)(u_{11}) = u_{21}$ and $(1,2,3)$ fixes all the
vertices in $V({\rm lk}_X(u_{11}))\cap V({\rm lk}_X(u_{21}))$.

Finally, from (\ref{edges}), $V({\rm lk}_X(v_{123})) = \{x_{111}$,
$x_{222}$, $x_{333}$, $x_{444}$, $x_{121}$, $x_{131}$, $x_{122}$,
$x_{331}$, $x_{332}$, $x_{322}$, $u_{41}$, $u_{62}$, $u_{54}\}$.
So, $V({\rm lk}_X(v_{123})) \cap V({\rm lk}_X(v_{213})) =
\{x_{111}$, $x_{222}$, $x_{333}$, $x_{444}$, $x_{331}$,
$x_{332}\}$. Then $(1,2)(v_{123}) = v_{213}$ and $(1,2)$ fixes all
the vertices in  $V({\rm lk}_X(v_{123})) \cap V({\rm
lk}_X(v_{213}))$. Since $V({\rm lk}_X(v_{123})) \cap V({\rm
lk}_X(v_{312}))= \{x_{111}$, $x_{222}$, $x_{333}$, $x_{444}\}$,
$(1,2,3)(v_{123}) = v_{312}$ and $(1,2,3)$ fixes all the vertices
in  $V({\rm lk}_X(v_{123})) \cap V({\rm lk}_X(v_{312}))$. So, the
$S_3$-action on $X$ is good. Therefore, the $S_3$-action on the
abstract scheme $(S^{\,2}\times S^{\,2} \times S^{\,2})_{124}$ of
$X$ is good. Now, Lemma \ref{L2.4} and Corollary \ref{L2.6} imply
that the
quotient $\CC P^{\,3}_{30} := (S^{\,2}\times S^{\,2} \times
S^{\,2})_{124}/S_3$ triangulates $\CC P^{\,3}$. \hfill $\Box$

\bigskip

The details of this proof follow.

\begin{lemma}$\!\!\!${\bf .} \label{R1}
Consider the isometry group $\ZZ_2$ $($generated by $\alpha : (x,
y) \mapsto (y, x))$ of the polytopal complex $S^{\,2}_4 \times
S^{\,2}_4$. There are exactly two $\ZZ_2$-stable $16$-vertex
simplicial subdivisions, say $X_0$ and $X_1$, of the $2$-skeleton
$U^2$ of $S^{\,2}_4 \times S^{\,2}_4$, on which the $\ZZ_2$-action
is pure. These two complexes are isometric.
\end{lemma}

\noindent {\bf Proof.} We use the notation $x_{ij}$ to denote the
vertex $(x_i, x_j)$ of $S^{\,2}_4 \times S^{\,2}_4$. So, the
vertex-set of $S^{\,2}_4 \times S^{\,2}_4$ is $\{x_{ij} ~ : 1\leq
i, j \leq 4\}$.

The 2-polytopes in $S^{\,2}_4 \times S^{\,2}_4$ are of the form
$\alpha \times \beta$, where one of $\alpha$, $\beta$ is a
2-simplex and other one is a vertex or both are edges. In the
first case, $\alpha \times \beta$ is a simplex and in the
second case, $\alpha \times \beta$ is not a simplex.

Let $X$ be a simplicial subdivisions of $U^2$ without adding any
new vertex such that the above $\ZZ_2$-action on the vertex-set
induces a good action on $X$. Then $X$ contains simplicial
subdivisions of all the 2-polytopes of $U^2$ which are
non-simplices.

For $i\neq j$, consider the 2-polytope
$x_ix_j \times x_ix_j$. Since the $\ZZ_2$-action is pure,
by condition $(a)$, $x_{ij} x_{ji}$ can not be an edge in $X$.
This implies that $x_{ii}x_{jj}$ is an edge in $X$ and
$x_{ii}x_{jj}\ast S^{\,0}_2(x_{ji}, x_{ij})$ is the
subdivision of $x_ix_j \times x_ix_j$ in $X$ (cf. Figure 1 (a)).

For $i, j, k$ distinct, consider the 2-polytope $x_ix_j \times
x_ix_k$. Since $x_{ik}x_{ij}$ is an edge in the polytopal complex
$S^{\,2}_4 \times S^{\,2}_4$, it is an edge in $X$. Therefore, by
condition $(b)$, $x_{ik} x_{ji}$ can not be an edge in $X$. This
implies that $x_{ii}x_{jk}$ is an edge in $X$ and $X[x_ix_j\times
x_ix_k] = x_{ii}x_{jk}\ast S^{\,0}_2(x_{ji}, x_{ik})$  (cf. Figure
1 (b)). Thus, $x_{ii}x_{jk}$ are edges in $X$ for all $i, j, k$.

\bigskip

\setlength{\unitlength}{3mm}

\begin{picture}(48.5,14)(2,-3.5)

\thicklines

\put(3,4){\line(1,0){5}} \put(3,4){\line(0,1){5}}
\put(3,9){\line(1,0){5}} \put(8,4){\line(0,1){5}}

\thinlines

\put(1,-3.5){\line(1,0){49}} \put(1,11){\line(1,0){49}}

\put(3,4){\line(1,1){5}}

\put(2,3){\mbox{$x_{ii}$}} \put(7.4,3){\mbox{$x_{ji}$}}
\put(2,9.7){\mbox{$x_{ij}$}} \put(7.4,9.7){\mbox{$x_{jj}$}}
\put(2.2,1.3){\mbox{$x_{i}x_j\times x_ix_j$}}


\thicklines

\put(13,4){\line(1,0){5}} \put(13,4){\line(0,1){5}}
\put(13,9){\line(1,0){5}} \put(18,4){\line(0,1){5}}

\thinlines

\put(13,4){\line(1,1){5}}

\put(12,3){\mbox{$x_{ii}$}} \put(17.4,3){\mbox{$x_{ji}$}}
\put(12,9.7){\mbox{$x_{ik}$}} \put(17.4,9.7){\mbox{$x_{jk}$}}
\put(12.2,1.3){\mbox{$x_ix_j\times x_ix_k$}}


\thicklines

\put(23,4){\line(1,0){5}} \put(23,4){\line(0,1){5}}
\put(23,9){\line(1,0){5}} \put(28,4){\line(0,1){5}}

\thinlines

\put(23,4){\line(1,1){5}}

\put(22,3){\mbox{$x_{12}$}} \put(27.4,3){\mbox{$x_{32}$}}
\put(22,9.7){\mbox{$x_{14}$}} \put(27.4,9.7){\mbox{$x_{34}$}}
\put(22.2,1.3){\mbox{$x_{1}x_3 \times x_2x_4$}}


\thicklines

\put(33,4){\line(1,0){5}} \put(33,4){\line(0,1){5}}
\put(33,9){\line(1,0){5}} \put(38,4){\line(0,1){5}}

\thinlines

\put(38,4){\line(-1,1){5}}

\put(32,3){\mbox{$x_{21}$}} \put(37.4,3){\mbox{$x_{31}$}}
\put(32,9.7){\mbox{$x_{24}$}} \put(37.4,9.7){\mbox{$x_{34}$}}
\put(32.2,1.3){\mbox{$x_{2}x_3\times x_1x_4$}}


\thicklines

\put(43,4){\line(1,0){5}} \put(43,4){\line(0,1){5}}
\put(43,9){\line(1,0){5}} \put(48,4){\line(0,1){5}}

\thinlines

\put(48,4){\line(-1,1){5}}

\put(42,3){\mbox{$x_{13}$}} \put(47.4,3){\mbox{$x_{23}$}}
\put(42,9.7){\mbox{$x_{14}$}} \put(47.4,9.7){\mbox{$x_{24}$}}
\put(42.2,1.3){\mbox{$x_{1}x_2\times x_3x_4$}}


\put(4,-0.6){\boldmath $(a)$} \put(14,-0.6){\boldmath $(b)$}
\put(24,-0.6){\boldmath $(c)$} \put(34,-0.6){\boldmath $(d)$}
\put(44,-0.6){\boldmath $(e)$}

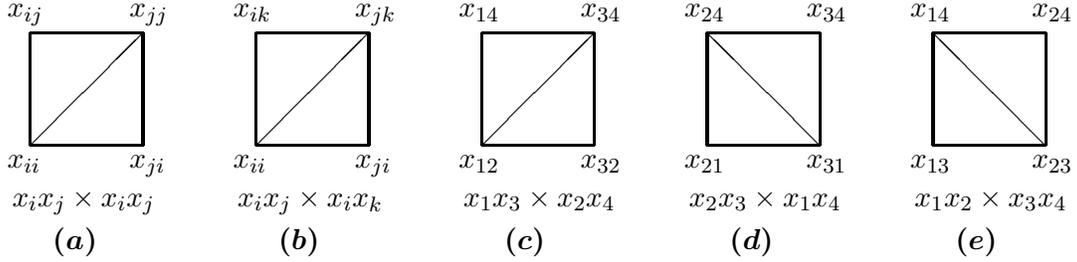
\put(2,-2.5){\mbox{\bf Figure 1\,: Simplicial subdivisions of
rectangular 2-cells of {\boldmath $S^{\,2}_4 \times S^{\,2}_4$}}}

\end{picture}

\medskip

Consider the 2-polytope $x_1x_3 \times x_2x_4$. Clearly, $X[x_1x_3
\times x_2x_4] = x_{12}x_{34}\ast S^{\,0}_2(x_{32}, x_{14})$ or $=
x_{32}x_{14} \ast S^{\,0}_2(x_{12}, x_{34})$.

\medskip

\noindent {\bf Case 1.} $X[x_1x_3 \times x_2x_4] = x_{12}x_{34}
\ast S^{\,0}_2(x_{32}, x_{14})$ (cf. Figure 1 (c)). Therefore
$x_{12}x_{34} \in X$ and hence $x_{21}x_{43}\in X$. Then, by
condition $(b)$, $x_{12} x_{43}$, $x_{21}x_{34} \not\in X$. This
implies that $X[x_2x_3\times x_1x_4] = x_{31}x_{24} \ast
S^0_2(x_{21}, x_{34})$ (cf. Figure 1 (d)). So, $x_{31}x_{24}\in X$
and hence $x_{13}x_{42} \in X$. Then, by condition $(b)$,
$x_{13}x_{24}$, $x_{31}x_{42} \not\in X$. This implies that
$X[x_1x_2\times x_3x_4] = x_{23}x_{14} \ast S^{\,0}_2(x_{13},
x_{24})$ (cf. Figure 1 (e)). So, $x_{14} x_{23} \in X$ and hence
$x_{41}x_{32}\in X$. Then, by condition $(b)$, $x_{14}x_{32}$,
$x_{41}x_{23} \not\in X$. Observe that, for $i, j, k, l$ distinct,
$x_{ij}x_{kl}$ is an edge if and only if $ijkl$ is an even
permutation of $1234$. This gives all the simplices in $X$. Denote
this $X$ by $X_0$. Observe that the set $E(X_0)$ of edges in $X_0$
is the following\,:
\begin{eqnarray}
E(X_0)=  & \hspace{-2mm} \{x_{ij}x_{ik} \,
: \, 1\leq i, j \neq k\leq 4\}\cup
\{x_{ii}x_{jk} \, : \, 1\leq j\neq i\neq k\leq 4\}
\nonumber \\
& \cup
\{x_{ij}x_{kl} \, : \,
ijkl \mbox{ is an even permutation of } 1234\}. \label{EX0}
\end{eqnarray}

\noindent {\bf Case 2.} $X[x_1x_3 \times x_2x_4] = x_{14}x_{32}
\ast S^{\,0}_2(x_{12}, x_{34})$. By the same argument we get the
subdivision of all the 2-polytopes. Denote this $X$ by $X_1$. The
set $E(X_1)$ of edges in $X_1$ is the following\,:
\begin{eqnarray}
E(X_1)=  & \hspace{-2mm} \{x_{ij}x_{ik} \,
: \, 1\leq i, j \neq k\leq 4\}\cup
\{x_{ii}x_{jk} \, : \, 1\leq j\neq i\neq k\leq 4\}
\nonumber \\
& \cup
\{x_{ij}x_{kl} \, : \,
ijkl \mbox{ is an odd permutation of } 1234\}. \nonumber
\end{eqnarray}
The complex $X_1$ is isomorphic to $X_0$ via the map $f$ given by
the transposition $(1, 2) \in S_4$ on the suffices. It is easy to
check that $X_i$ is $\ZZ_2$-invariant and the $\ZZ_2$-action is
good and hence pure for $0 \leq i\leq 1$. \hfill $\Box$

\begin{remark}$\!\!\!${\bf .} \label{remark3}
{\rm Recall that a simplicial complex is said to be a {\em clique
complex} if it consists of all the cliques (sets of mutually
adjacent vertices) of its edge graph. Each of $X_0$ and $X_1$ is a
clique complex, hence is determined uniquely by its edge graph. }
\end{remark}

\begin{lemma}$\!\!\!${\bf .} \label{R2}
Up to isometry, there are exactly three $S_3$-invariant
$64$-vertex simplicial subdivisions of the $2$-skeleton $W^2$ of
the polytopal complex $S^{\,2}_4 \times S^{\,2}_4 \times
S^{\,2}_4$, say $X^2_i$, $0\leq i \leq 2$ such that the induced
$S_3$ action on each of them is pure.
\end{lemma}

\noindent {\bf Proof.} Let $X^2$ be a simplicial subdivision of
$W^2$ without adding any new vertex such that the $S_3$-action on
the vertex-set induces a pure action on $X^2$. Then the vertex-set
of $X^2$ is $\{x_{ijk} ~ : 1\leq i, j,k \leq 4\}$ and $X^2$
contains simplicial subdivisions of all the 2-polytopes of
$S^{\,2}_4 \times S^{\,2}_4 \times S^{\,2}_4$ which are
non-simplices.

Notice that $W^2$ is the union of twelve 2-dimensional polyhedral
complexes, namely, $W^2 = \cup\{S_3(x_h \times U^2) \, : \, U^2$
is the 2-skeleton of $S^{\,2}_4 \times S^{\,2}_4$, $1\leq h \leq
4\}$. Therefore, the $S_3$-stable 64-vertex subdivision $X^2$ of
$W^2$ is determined by its four subcomplexes $X^2[x_h \times
U^2]$, $1\leq h \leq 4$. On each of them, the action of $\ZZ_2
\subseteq S_3$ generated by $(2, 3)$ must be pure. Therefore, by
Lemma \ref{R1}, for each $h \in \{1, 2, 3, 4\}$, we must have
$X^2[x_h \times U^2] = x_h \times X_0$ or $x_h \times X_1$. Since
there is a group $S_4$ of isometries  of $S^{\,2}_4 \times
S^{\,2}_4 \times S^{\,2}_4$ acting naturally on the indices $1, 2,
3, 4$ (i.e., $S_4 \ni \pi : x_{ijk} \mapsto x_{\pi(i)\pi(j)
\pi(k)}$) which commutes with $S_3$, it is now obvious that there
are exactly three non-isomorphic choices $X^2_i$ ($0 \leq i \leq
2$), where $X^2_i[x_h \times U^2] = x_h \times X_1$ for exactly
$i$ values of $h$. (We have $X^2_i \cong X^2_{4-i}$ via any odd
element of $S_4$. Hence the restriction $i \leq 2$.) \hfill $\Box$

\begin{lemma}$\!\!\!${\bf .} \label{R3}
There does not exist any $64$-vertex $S_3$-invariant simplicial
subdivision of the $3$-skeleton $W^3$ of the polytopal complex
$S^{\,2}_4 \times S^{\,2}_4 \times S^{\,2}_4$ on which the induced
$S_3$ action is pure.
\end{lemma}

\setlength{\unitlength}{3.9mm}

\begin{picture}(36,18.7)(0,-2.4)

\thicklines

\put(2,3){\line(1,0){8}} \put(2,3){\line(0,1){8}}
\put(2,11){\line(1,0){8}} \put(10,3){\line(0,1){8}}
\put(6,15){\line(1,0){8}} \put(14,7){\line(0,1){8}}
\put(6,7){\line(1,0){3.5}} \put(10.5,7){\line(1,0){3.5}}
\put(6,7){\line(0,1){3.5}} \put(6,11.5){\line(0,1){3.5}}

\put(2,3){\line(1,1){4}} \put(2,11){\line(1,1){4}}
\put(10,3){\line(1,1){4}} \put(10,11){\line(1,1){4}}

\put(18,6){\line(1,0){8}} \put(18,6){\line(3,5){6}}

\put(22,10){\line(1,3){2}} \put(24,0){\line(1,3){2}}
\put(24,0){\line(3,5){6}}  \put(22,10){\line(1,0){2.7}}
\put(25.6,10){\line(1,0){0.8}} \put(30,10){\line(-1,0){2.8}}

\put(24,0){\line(-2,5){2}} \put(18,6){\line(4,-1){4}}
\put(22,10){\line(0,-1){3.7}}  \put(22,5){\line(0,1){0.6}}
\put(26,6){\line(1,5){1}} \put(30,10){\line(-3,1){3}}
\put(27,11){\line(-3,5){3}}

\thinlines

\put(-1,-1.9){\line(1,0){38}} \put(-1,16.5){\line(1,0){38}}

\put(18,6){\line(1,1){4}} \put(24,0){\line(-1,1){6}}
\put(24,0){\line(-1,5){1.1}} \put(22,10){\line(1,-5){0.7}}
 \put(26,6){\line(1,1){4}}  \put(26,6){\line(-1,5){2}}
\put(30,10){\line(-1,1){6}}

\put(3.3,3.6){\mbox{$x_{112}$}} \put(8,3.6){\mbox{$x_{212}$}}
\put(6.5,7.6){\mbox{$x_{132}$}} \put(12,7.7){\mbox{$x_{232}$}}
\put(2.4,10){\mbox{$x_{113}$}} \put(8,10){\mbox{$x_{213}$}}
\put(6.5,14.2){\mbox{$x_{133}$}} \put(11,14.2){\mbox{$x_{233}$}}

\put(25,15.5){\mbox{$x_{233}$}} \put(25,0.2){\mbox{$x_{112}$}}
\put(22.5,10.6){\mbox{$x_{232}$}} \put(23.7,6.6){\mbox{$x_{113}$}}
\put(30.2,9.2){\mbox{$x_{212}$}} \put(16.3,6.7){\mbox{$x_{133}$}}
\put(20,4.4){\mbox{$x_{132}$}} \put(25.1,10.8){\mbox{$x_{213}$}}


\put(2,1){\mbox{\bf (a): \boldmath{$C= x_1x_2\times x_1x_3
\times x_2x_3 = C^3_5$}}}

\put(18,-1.3){\mbox{\bf (b): \boldmath{$X^2_i[C]$}}}

\put(4,-1.3){\mbox{\bf Figure 2}}

\end{picture}


\noindent {\bf Proof.} Assume, on the contrary, that there is such
a subdivision $X^3$ of $W^3$. Consider the 3-polytope $C = x_1x_2
\times x_1x_3 \times x_2x_3$ in $W^3$. By Lemma \ref{R2},
$X^3[\partial \overline{C}] \cong X^2_i[\partial \overline{C}]$
for some $i$ ($0\leq i \leq 2$). But, for each choice of $i$,
$x_{112}x_{113}x_{133}$ is a triangle of $X^2_i$ and
$x_{133}x_{123}$, $x_{133}x_{122}$, $x_{112}x_{332}$,
$x_{113}x_{123}$, $x_{113}x_{223} \in X^2_i$. Therefore, by
condition $(b)$ of purity, $x_{133}x_{213}$, $x_{133}x_{212}$,
$x_{112}x_{233}$, $x_{113}x_{132}$, $x_{113} x_{232} \not\in
X^2_i$. So, the same is true for $X^3$. Then $X^3$ must contain a
tetrahedron $\sigma \supseteq x_{112}x_{113} x_{133}$ with
$V(\sigma) \subseteq V(C)$, and there is no valid choice of
$\sigma$ because of these five non-edges (cf. Fig. 2). \hfill
$\Box$

\begin{remark}$\!\!\!${\bf .} \label{remark4}
{\rm The proof of Lemma \ref{R3} shows that we need at least 24
more vertices to get simplicial subdivisions of the 24 3-polytopes
$x_ix_j \times x_ix_k \times x_jx_k$, where $i, j, k$ are
distinct. We will see (in the proof of the next lemma) that there
are some more 3-polytopes in $S^{\,2}_4 \times S^{\,2}_4 \times
S^{\,2}_4$ which do not allow simplicial subdivision without
adding new vertices, while maintaining the purity of the
$S_3$-action.}
\end{remark}

\begin{lemma}$\!\!\!${\bf .} \label{R4}
There exists an $124$-vertex simplicial subdivision $X^3$ of the
$3$-skeleton $W^3$ of  $S^{\,2}_4 \times S^{\,2}_4 \times
S^{\,2}_4$, such that $X^3$ admits an isometry group $S_3 \times
A_4$. The action of $S_3 \times A_4$ on $|W^3|$ induced by $X^3$
is the same as the action induced by its action on $W^3$.
Further, the action of $S_3$ on $X^3$ is good.
\end{lemma}

\noindent {\bf Proof.} We take $X^2_0$ (defined in the proof of
Lemma \ref{R2}) as the simplicial subdivision $X^2$ of the
2-skeleton $W^2$ of $W^3$. The set $\{x_{ijk} ~ : 1\leq i, j, k
\leq 4\}$ is the vertex-set of both $X^2$ and $S^{\,2}_4 \times
S^{\,2}_4 \times S^{\,2}_4$. From (\ref{EX0}) it follows that
$E(X^2)$ (the set of edges in $X^2$) is the following\,:

\begin{eqnarray}
E(X^2) =  &\hspace{-2mm} S_3(\{x_{hij}x_{hik} \, : \, 1\leq h, i,
j \neq k\leq 4\})\cup S_3(\{x_{hii}x_{hjk} \, : \, 1\leq h, j\neq
i\neq k\leq 4\})
\nonumber \\
& \cup \, S_3(\{x_{hij}x_{hkl} \, : \, 1\leq h\leq 4, \, ijkl
\mbox{ is an even permutation of } 1234\}). \label{EX2}
\end{eqnarray}

Recall that $S_3$ acts on these vertices by permuting the
positions of the three subscripts. Now, consider the action of
$S_4$ on these vertices by permuting the values of these
subscripts (which are elements of $\{1, 2, 3, 4\}$). It is easy to
see that $S_3 \times A_4 \subseteq S_3 \times S_4$ is a group of
isometries of $X^2$. We try to construct $X^3$ retaining the group
$S_3 \times A_4$ and with the provision that the $S_3$-action is
good.

We have to subdivide all the 3-polytopes. A 3-polytope in $W^3$ is
either the product of one 0-simplex, one 1-simplex and one
2-simplex or the product of three 1-simplices.

In the first case, modulo the $S_3$-action, the 3-polytopes are of
the form $x_ix_jx_k\times x_ix_j \times x_h$ for $i, j, k$
distinct or of the form $x_ix_jx_k \times x_ix_l \times x_h$ for
$i, j, k, l$ distinct. Since $x_{iih}x_{jjh}$, $x_{iih}x_{kjh}$
and $x_{jjh}x_{kih}$ are edges, by Lemma \ref{L2.1}, $x_ix_jx_k
\times x_ix_j \times x_h = x_{iih}x_{ijh}x_{kjh}x_{jjh} \cup
x_{iih}x_{kih}x_{kjh}x_{jjh} \cup x_{iih}x_{kih}x_{jih}x_{jjh}$ is
the unique subdivision of $x_ix_jx_k \times x_ix_j \times x_h$
(see Figure 3). For $i, j, k, l$ distinct, consider the 3-polytope
$x_ix_jx_k \times x_ix_l\times x_h$. Here $x_{iih}x_{jlh}$ and
$x_{iih}x_{klh}$ are edges. By interchanging $j$ and $k$ (if
required) we may assume that $ijkl$ (and hence also $kijl$) is an
even permutation of $1234$. Then $x_{kih} x_{jlh}$ is an edge and
hence, by Lemma \ref{L2.1}, $x_ix_jx_k \times x_ix_l \times x_h =
x_{iih}x_{ilh}x_{klh}x_{j\,lh} \cup x_{iih}x_{kih}x_{klh}x_{jlh}
\cup x_{iih}x_{kih}x_{jih}x_{jlh}$ is the unique subdivision of
$x_ix_jx_k \times x_ix_l\times x_h$.

\medskip

\setlength{\unitlength}{4mm}

\begin{picture}(36.5,14)(3,-1)


\thicklines

\put(6,5){\line(-3,1){3}} \put(6,5){\line(3,1){3}}
\put(6,5){\line(0,1){5}} \put(3,6){\line(1,0){0.75}}
\put(4.25,6){\line(1,0){0.75}} \put(8.25,6){\line(1,0){0.75}}
\put(7,6){\line(1,0){0.75}} \put(3,6){\line(0,1){5}}
\put(9,6){\line(0,1){5}}

\put(6,10){\line(-3,1){3}} \put(6,10){\line(3,1){3}}
\put(3,11){\line(1,0){6}}

\thinlines

\put(2,-1){\line(1,0){37}} \put(2,13){\line(1,0){37}}

\put(6,5){\line(1,2){3}} \put(6,5){\line(-1,2){3}}

\put(9,6){\line(-6,5){1.5}} \put(3,11){\line(6,-5){1.1}}
\put(4.5,9.75){\line(6,-5){1.1}}
\put(6.15,8.375){\line(6,-5){0.88}} 

\put(6.3,4.3){\mbox{$x_{iih}$}} \put(6.3,9.4){\mbox{$x_{ijh}$}}
\put(3,5){\mbox{$x_{jih}$}} \put(3,11.5){\mbox{$x_{jjh}$}}
\put(8,5){\mbox{$x_{kih}$}} \put(8,11.5){\mbox{$x_{kjh}$}}

\put(3,3){\mbox{${x_ix_jx_k\!\times\!x_ix_j\!\times\!x_h}$}}

\thicklines

\put(15,5){\line(-3,1){3}} \put(15,5){\line(3,1){3}}
\put(15,5){\line(0,1){5}} \put(12,6){\line(1,0){0.75}}
\put(13.25,6){\line(1,0){0.75}} \put(17.25,6){\line(1,0){0.75}}
\put(16,6){\line(1,0){0.75}} \put(12,6){\line(0,1){5}}
\put(18,6){\line(0,1){5}}

\put(15,10){\line(-3,1){3}} \put(15,10){\line(3,1){3}}
\put(12,11){\line(1,0){6}}

\thinlines

\put(15,5){\line(1,2){3}} \put(15,5){\line(-1,2){3}}

\put(18,6){\line(-6,5){1.5}} \put(12,11){\line(6,-5){1.1}}
\put(13.5,9.75){\line(6,-5){1.1}}
\put(15.15,8.375){\line(6,-5){0.88}} 

\put(15.3,4.3){\mbox{$x_{iih}$}} \put(15.3,9.4){\mbox{$x_{ilh}$}}
\put(12,5){\mbox{$x_{jih}$}} \put(12,11.5){\mbox{$x_{j\,lh}$}}
\put(17,5){\mbox{$x_{kih}$}} \put(17,11.5){\mbox{$x_{klh}$}}

\put(12,3){\mbox{${x_ix_jx_k\!\times\!x_ix_l\!\times\!x_h}$}}

\put(4,1.4){($ijkl$ an even permutation of $1234$)}

\put(4,0){\mbox{\bf Figure 3\,: Simplicial subdivision of a prism}}

\end{picture}

\medskip

Consider the 3-polytopes which are products of three 1-simplices.
Modulo the action of the group $S_3 \times A_4$, such 3-polytopes
are the following\,: $C^3_1 = x_{1}x_{2} \times x_{1}x_{2} \times
x_{1}x_{2}$, $C^3_2 = x_{1}x_{2} \times x_{1}x_{2} \times
x_{1}x_{3}$, $C^3_3 = x_{1}x_{2} \times x_{1}x_{3} \times
x_{1}x_{4}$, $C^3_4 = x_{1}x_{2} \times x_{1}x_{2} \times
x_{4}x_{3}$, $C^3_5 = x_{1}x_{2} \times x_{1}x_{3} \times
x_{2}x_{3}$, $C^3_6 = x_{1}x_{2} \times x_{2}x_{3} \times x_3x_4$.

Observe that $X^2[C^3_1] = S^{\,0}_2(x_{111}, x_{222}) \ast
S^1_6(x_{121}, x_{221}, x_{211}, x_{212}, x_{112}, x_{122})$ is
the subdivision of the boundary $\partial \overline{C^3_1}$ (see
Figure 4). Since $x_{111}x_{222}$ is the only possible extra edge
with vertices in $C^3_1$ such that its introduction does not
destroy the goodness of the $S_3$-action, it follows that
$$
\widehat{C^3_1} : = x_{111}x_{222} \ast
\begin{picture}(8,1)(-1,1)
\setlength{\unitlength}{2.5mm}
\put(1.7,0){$_{\bullet}$}
\put(1.7,3.85){$_{\bullet}$}
\put(4.8,0){$_{\bullet}$}
\put(4.8,3.85){$_{\bullet}$}
\put(-0.1,1.9){$_{\bullet}$}
\put(6.7,1.9){$_{\bullet}$}
\thicklines
\put(2,0){\line(1,0){3}}
\put(2,0){\line(-1,1){2}}
\put(5,0){\line(1,1){2}}
\put(5,4){\line(1,-1){2}}
\put(2,4){\line(1,0){3}}
\put(0,2){\line(1,1){2}}
\put(-1.5,3.9){\small {$x_{122}$}}
\put(1,1.8){\small {$x_{121}$}}
\put(-1.5,0){\small {$x_{221}$}}
\put(6.7,0){\small {$x_{211}$}}
\put(7.7,1.8){\small {$x_{212}$}}
\put(6.2,3.9){\small {$x_{112}$}}
\end{picture}
$$

\medskip

\noindent is the unique subdivision of $C^3_1$ without adding any
new vertex. (Since ${\rm ast}_{X^2[C^3_1]}(x_{111})$ is a
simplicial subdivision of ${\rm Ast}_{C^3_1}(x_{111})$ and
$\widehat{C^3_1} = x_{111}\ast {\rm ast}_{X^2[C^3_1]}(x_{111})$,
by Lemma \ref{L2.2}, $\widehat{C^3_1}$ triangulates  $C^3_1$.)

\medskip

\setlength{\unitlength}{3.9mm}

\begin{picture}(36,18)(0,-2)

\thicklines

\put(2,3){\line(1,0){8}} \put(2,3){\line(0,1){8}}
\put(2,11){\line(1,0){8}} \put(10,3){\line(0,1){8}}
\put(6,15){\line(1,0){8}} \put(14,7){\line(0,1){8}}
\put(6,7){\line(1,0){3.5}} \put(10.5,7){\line(1,0){3.5}}
\put(6,7){\line(0,1){3.5}}  \put(6,11.5){\line(0,1){3.5}}

\put(2,3){\line(1,1){4}} \put(2,11){\line(1,1){4}}
\put(10,3){\line(1,1){4}} \put(10,11){\line(1,1){4}}

\put(18,8){\line(2,-1){4}} \put(22,6){\line(1,0){6}}
\put(28,6){\line(2,1){4}} \put(18,8){\line(2,1){4}}
 \put(28,10){\line(2,-1){4}}
 \put(22,6){\line(1,3){3}}  \put(32,8){\line(-1,1){7}}
\put(22,10){\line(3,5){3}} \put(28,10){\line(-3,5){3}}

\put(22,10){\line(1,0){1}} \put(28,10){\line(-1,0){1}}
\put(23.8,10){\line(1,0){2.4}}

\put(25,1){\line(-1,1){7}} \put(25,1){\line(-3,5){3}}
\put(25,1){\line(3,5){3}}

\put(25,1){\line(1,3){1.55}}
\put(28,10){\line(-1,-3){0.5}} \put(27.2,7.6){\line(-1,-3){0.4}}

\thinlines

\put(-1,-2){\line(1,0){38}} \put(-1,16){\line(1,0){38}}

\put(18,8){\line(1,1){7}} \put(28,6){\line(-1,3){3}}
\put(25,1){\line(1,1){7}} \put(25,1){\line(-1,3){1.55}}
\put(22,10){\line(1,-3){0.5}} \put(22.8,7.6){\line(1,-3){0.4}}


\put(3.2,3.5){\mbox{$x_{111}$}} \put(8,3.5){\mbox{$x_{211}$}}
\put(6.5,7.6){\mbox{$x_{121}$}} \put(12,7.7){\mbox{$x_{221}$}}
\put(2.4,10){\mbox{$x_{112}$}} \put(8,10){\mbox{$x_{212}$}}
\put(6.5,14.2){\mbox{$x_{122}$}} \put(11,14.2){\mbox{$x_{222}$}}



\put(26,14.7){\mbox{$x_{222}$}} \put(26,1.3){\mbox{$x_{111}$}}
\put(20.3,8.5){\mbox{$x_{122}$}} \put(27.8,7.2){\mbox{$x_{211}$}}
\put(31.5,9){\mbox{$x_{212}$}} \put(16.3,6.7){\mbox{$x_{121}$}}
\put(20.2,7.2){\mbox{$x_{221}$}} \put(27.8,8.5){\mbox{$x_{112}$}}


\put(2,1){\mbox{\bf (a)\,: \boldmath{$C^3_1 = x_1x_2\times x_1x_2
\times x_1x_2$}}}

\put(18,-1){\mbox{\bf (b)\,: Boundary of \boldmath{$\widehat{C^3_1} = X^2[C^3_1]$}}}

\put(5,-1){\mbox{\bf Figure 4}}

\end{picture}

\medskip

Since $X^2[C^3_2] = S^{\,0}_2(x_{111}, x_{223}) \ast
S^1_6(x_{113}, x_{213}, x_{211}, x_{221}, x_{121}, x_{123})$ is
the subdivision of $\partial \overline{C^3_2}$, as before,
$$
\widehat{C^3_2} :=
x_{111}x_{223} \ast
\begin{picture}(8,1.5)(-1,1)
\setlength{\unitlength}{2.5mm}
\put(1.7,0){$_{\bullet}$}
\put(1.7,3.85){$_{\bullet}$}
\put(4.8,0){$_{\bullet}$}
\put(4.8,3.85){$_{\bullet}$}
\put(-0.1,1.9){$_{\bullet}$}
\put(6.7,1.9){$_{\bullet}$}
\thicklines
\put(2,0){\line(1,0){3}}
\put(2,0){\line(-1,1){2}}
\put(5,0){\line(1,1){2}}
\put(5,4){\line(1,-1){2}}
\put(2,4){\line(1,0){3}}
\put(0,2){\line(1,1){2}}
\put(-1.5,3.9){\small {$x_{123}$}}
\put(1,1.8){\small {$x_{113}$}}
\put(-1.5,0){\small {$x_{213}$}}
\put(6.7,0){\small {$x_{211}$}}
\put(7.7,1.8){\small {$x_{221}$}}
\put(6.2,3.9){\small {$x_{121}$}}
\end{picture}
$$
\smallskip

\noindent is the unique subdivision of $C^3_2$ without adding any
new vertex. (Since ${\rm ast}_{X^2[C^3_2]}(x_{111})$ is a
simplicial subdivision of ${\rm Ast}_{C^3_2}(x_{111})$ and
$\widehat{C^3_2} = x_{111}\ast {\rm ast}_{X^2[C^3_2]}(x_{111})$,
by Lemma \ref{L2.2}, $\widehat{C^3_2}$ triangulates $C^3_2$.)

\smallskip

Since $X^2[C^3_3] = S^{\,0}_2(x_{111}, x_{234}) \ast
S^1_6(x_{131}, x_{231}, x_{211}, x_{214}, x_{114}, x_{134})$ is
the subdivision of $\partial \overline{C^3_3}$ and
$x_{111}x_{234}$ is the only possible new edge whose vertices are
in $C^3_3$, it follows that
$$
\widehat{C^3_3} :=
x_{111}x_{234} \ast
\begin{picture}(8,1.8)(-1,1)
\setlength{\unitlength}{2.5mm}
\put(1.7,0){$_{\bullet}$}
\put(1.7,3.85){$_{\bullet}$}
\put(4.8,0){$_{\bullet}$}
\put(4.8,3.85){$_{\bullet}$}
\put(-0.1,1.9){$_{\bullet}$}
\put(6.7,1.9){$_{\bullet}$}
\thicklines
\put(2,0){\line(1,0){3}}
\put(2,0){\line(-1,1){2}}
\put(5,0){\line(1,1){2}}
\put(5,4){\line(1,-1){2}}
\put(2,4){\line(1,0){3}}
\put(0,2){\line(1,1){2}}
\put(-1.5,3.9){\small {$x_{134}$}}
\put(1,1.8){\small {$x_{131}$}}
\put(-1.5,0){\small {$x_{231}$}}
\put(6.7,0){\small {$x_{211}$}}
\put(7.7,1.8){\small {$x_{214}$}}
\put(6.2,3.9){\small {$x_{114}$}}
\end{picture}
$$

\medskip

\noindent is the unique subdivision of $C^3_3$ without adding any
new vertex. (Since ${\rm ast}_{X^2[C^3_3]}(x_{111})$ is a simplicial
subdivision of ${\rm Ast}_{C^3_3}(x_{111})$ and $\widehat{C^3_3} =
x_{111}\ast {\rm ast}_{X^2[C^3_3]}(x_{111})$, by Lemma \ref{L2.2},
$\widehat{C^3_3}$ triangulates $C^3_3$.)

Consider the polytope $C^3_4$. Observe that $X^2_i[x_1x_2 \times
x_1x_2 \times x_k] = x_{11k}x_{22k} \ast S^{\,0}_2(x_{21k},
x_{12k})$ for $k = 3, 4$ and for all $i$ and $X^2[x_h \times
x_1x_2 \times x_4x_3] = x_{h14}x_{h23}\ast S^{\,0}_2(x_{h13},
x_{h24})$, $X^2[x_1x_2 \times x_h\times x_4x_3]=
x_{1h4}x_{2h3}\ast S^{\,0}_2(x_{1h3}, x_{2h4})$ for $h = 1, 2$.
Thus, $X^2[C^3_4] = S^{\,0}_2(x_{114}, x_{223}) \ast
S^1_6(x_{113}, x_{123}, x_{124}, x_{224}, x_{214}, x_{213})$ (cf.
Fig. 5). Since $X^2 \subseteq X^3$, by condition $(b)$,
$x_{114}x_{223}$ and $x_{113}x_{224}$ are the only possible new
edges whose vertices are in $C^3_4$. Thus,
$$
\widehat{C^3_4} : = x_{114}x_{223} \ast
\begin{picture}(8,1.2)(-1,1)
\setlength{\unitlength}{2.5mm}
\put(1.7,0){$_{\bullet}$}
\put(1.7,3.85){$_{\bullet}$}
\put(4.8,0){$_{\bullet}$}
\put(4.8,3.85){$_{\bullet}$}
\put(-0.1,1.9){$_{\bullet}$}
\put(6.7,1.9){$_{\bullet}$}
\thicklines
\put(2,0){\line(1,0){3}}
\put(2,0){\line(-1,1){2}}
\put(5,0){\line(1,1){2}}
\put(5,4){\line(1,-1){2}}
\put(2,4){\line(1,0){3}}
\put(0,2){\line(1,1){2}}
\put(-1.5,3.9){\small {$x_{213}$}}
\put(1,1.8){\small {$x_{113}$}}
\put(-1.5,0){\small {$x_{123}$}}
\put(6.7,0){\small {$x_{124}$}}
\put(7.7,1.8){\small {$x_{224}$}}
\put(6.2,3.9){\small {$x_{214}$}}
\end{picture}
$$

\medskip

\noindent is the unique subdivision of $C^3_4$ without adding any
new vertex. (Since ${\rm ast}_{X^2[C^3_4]}(x_{114})$ is a
simplicial subdivision of ${\rm Ast}_{C^3_4}(x_{114})$ and
$\widehat{C^3_4} = x_{114}\ast {\rm ast}_{X^2[C^3_4]}(x_{111})$,
by Lemma \ref{L2.2}, $\widehat{C^3_4}$ triangulates $C^3_4$.)
Here, we have the new edge $x_{114}x_{223}$. In general, for the
3-polytope $x_ix_k \times x_ix_k \times x_jx_l$, we can assume, by
interchanging $j$ and $l$ (if necessary), that $ijkl$ is an even
permutation of $1234$. We subdivide it by adding the new edge
$x_{iij}x_{kkl}$.

\medskip

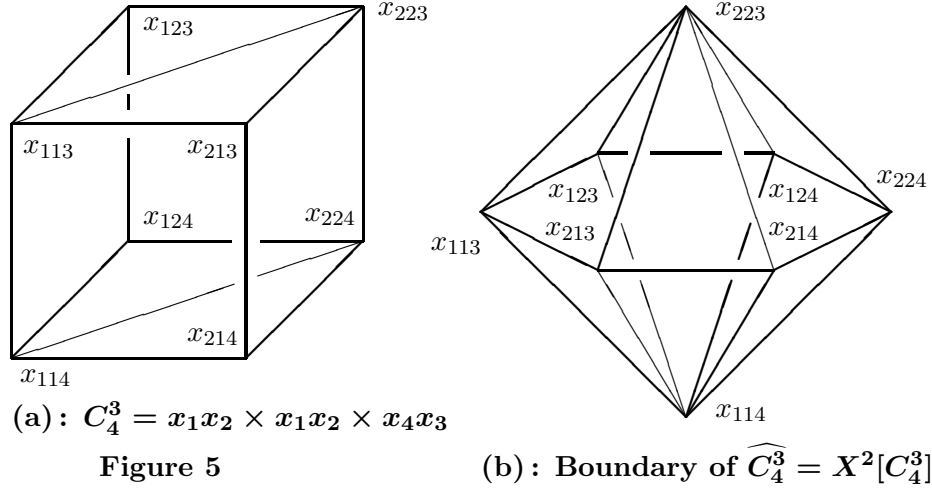
\setlength{\unitlength}{3.9mm}

\begin{picture}(36,18)(0,-2)

\thicklines

\put(2,3){\line(1,0){8}} \put(2,3){\line(0,1){8}}
\put(2,11){\line(1,0){8}} \put(10,3){\line(0,1){8}}
\put(6,15){\line(1,0){8}} \put(14,7){\line(0,1){8}}
\put(6,7){\line(1,0){3.5}} \put(10.5,7){\line(1,0){3.5}}
\put(6,7){\line(0,1){3.5}}  \put(6,11.5){\line(0,1){0.5}}
\put(6,12.7){\line(0,1){2.3}}

\put(2,3){\line(1,1){4}} \put(2,11){\line(1,1){4}}
\put(10,3){\line(1,1){4}} \put(10,11){\line(1,1){4}}

\put(18,8){\line(2,-1){4}} \put(22,6){\line(1,0){6}}
\put(28,6){\line(2,1){4}} \put(18,8){\line(2,1){4}}
 \put(28,10){\line(2,-1){4}}
 \put(22,6){\line(1,3){3}}  \put(32,8){\line(-1,1){7}}
\put(22,10){\line(3,5){3}} %

\put(22,10){\line(1,0){1}} \put(28,10){\line(-1,0){1}}
\put(23.8,10){\line(1,0){2.4}}

\put(25,1){\line(-1,1){7}} %
\put(25,1){\line(3,5){3}}

\put(25,1){\line(1,3){1.55}}
\put(28,10){\line(-1,-3){0.5}} \put(27.2,7.6){\line(-1,-3){0.4}}

\put(25,1){\line(1,1){7}} \put(18,8){\line(1,1){7}}

\thinlines

\put(-1,-2){\line(1,0){38}} \put(-1,16){\line(1,0){38}}

\put(28,6){\line(-1,3){3}} \put(25,1){\line(-3,5){3}}
 \put(25,1){\line(-1,3){1.55}} \put(28,10){\line(-3,5){3}}
\put(22,10){\line(1,-3){0.5}} \put(22.8,7.6){\line(1,-3){0.4}}

\put(2,11){\line(3,1){12}} \put(2,3){\line(3,1){7.6}}
\put(14,7){\line(-3,-1){3.6}}

\put(2.3,2.2){\mbox{$x_{114}$}} \put(8,3.6){\mbox{$x_{214}$}}
\put(6.5,7.6){\mbox{$x_{124}$}} \put(12,7.7){\mbox{$x_{224}$}}
\put(2.4,10){\mbox{$x_{113}$}} \put(8,10){\mbox{$x_{213}$}}
\put(6.5,14.2){\mbox{$x_{123}$}} \put(14.5,14.7){\mbox{$x_{223}$}}



\put(26,14.7){\mbox{$x_{223}$}} \put(26,1){\mbox{$x_{114}$}}
\put(16.3,6.7){\mbox{$x_{113}$}} \put(20.3,8.5){\mbox{$x_{123}$}}
\put(27.8,8.5){\mbox{$x_{124}$}} \put(31.5,9){\mbox{$x_{224}$}}
\put(27.8,7.2){\mbox{$x_{214}$}} \put(20.2,7.2){\mbox{$x_{213}$}}


\put(2,0.7){\mbox{\bf (a)\,: \boldmath{$C^3_4 = x_1x_2\times x_1x_2
\times x_4x_3$}}}

\put(18,-1){\mbox{\bf (b)\,: Boundary of \boldmath{$\widehat{C^3_4} = X^2[C^3_4]$}}}

\put(5,-1){\mbox{\bf Figure 5}}

\end{picture}

\medskip

The vertices of $C^3_5$ are $x_{112} = (1,-1,-1,1,-1,-1,-1,1,-1)$,
$x_{212} = (-1,1,-1,1,-1,-1$, $-1,1,-1)$, $x_{132} = (1,-1,-1,-1,
-1, 1,-1,1,-1)$, $x_{232} = (-1,1,-1,-1,-1,1,-1,1,-1)$, $x_{113} =
(1,-1,-1,1,-1,-1,-1,-1,1)$, $x_{213} = (-1,1,-1,1,-1,-1,-1,-1,1)$,
$x_{133} = (1,-1,$ $-1,-1,-1,1,-1,-1, 1)$, $x_{233} = (-1,1,-1,-1,
-1,1,-1,-1,1)$. The simplicial subdivision $X^2[C^3_5]$ of
$\partial \overline{C^3_5}$ induced by $X^2$ is given in Figure 2.
We know, from the proof of Lemma \ref{R3}, that $C^3_5$ has no
simplicial subdivision without adding any new vertex. We add the
new vertex $v_{321} = (0, 0, -1, 0, -1, 0, -1, 0, 0)$ in the
interior of $C^3_5$. Thus $v_{321}$ is the barycenter of the
polytope $C^3_5$. By Lemma \ref{L2.2}, $K := v_{321}\ast
X^2[C^3_5]$ is a simplicial subdivision of $C^3_5$. Since
$\frac{1}{3}x_{132}+ \frac{2}{3} v_{321} = \frac{1}{3}(1, -1, -3,
-1, -3, 1, -3, 1,-1)= \frac{1}{3} x_{112} + \frac{1}{3}x_{133} +
\frac{1}{3}x_{232}$, it follows that the edge $v_{321}x_{132}$ and
the 2-simplex $x_{112}x_{133} x_{232}$ intersect at an interior
points of both. This implies that the 3-polytope
$\langle\{v_{321},  x_{132}, x_{112}, x_{133}, x_{232} \}\rangle =
v_{321}x_{132}x_{112}x_{133} \cup v_{321}x_{132}x_{112}x_{232}\cup
v_{321}x_{132}x_{133}x_{232} = v_{321}x_{112}x_{133}x_{232} \cup
x_{132}x_{112}x_{133}x_{232}$. So, we can replace
$v_{321}x_{132}\ast S^1_3(\{x_{112}, x_{133}, x_{232}\})$ by
$x_{112}x_{133}x_{232}\ast S^{\,0}_2(v_{321}, x_{132})$ in $K$ and
get a new subdivision $L$ of $C^3_5$. Similarly, we can replace
$v_{321}x_{213}\ast S^1_3(\{x_{233}, x_{212}, x_{113}\})$ by
$x_{233}x_{212}x_{113} \ast S^{\,0}_2(v_{321}, x_{213})$ in $L$
and get a new simplicial subdivision $\widehat{C^3_5}$ of $C^3_5$.
So,
$$
\widehat{C^3_5} =
v_{321} \ast S^{0}_2(x_{112}, x_{233}) \ast
S^{0}_2(x_{133}, x_{212}) \ast S^{0}_2(x_{232}, x_{113})\cup
\{x_{132}x_{112}x_{133}x_{232}, x_{213}x_{233}x_{212}x_{113}\}.
$$
Using the action of $S_3\times A_4$, we subdivide each polytope
$x_jx_k \times x_ix_k \times x_ix_j$ by adding the vertex
$v_{ijk}$, where $i, j, k$ are distinct. (Observe that $x_jx_k
\times x_ix_k \times x_ix_j$ and its subdivision are invariant
under the order three automorphism $(1, 2, 3)(i, j, k) \in S_3
\times A_4$.) Thus $v_{ijk}$ is the barycenter of the polytope
$x_jx_k \times x_ix_k \times x_ix_j$, and the action of $S_3
\times A_4$ on these polytopes induces the natural action (in
terms of the subscripts) on these new vertices.

The vertices of $C^3_6$ are
$x_{123} = (1, -1,-1,-1,1,-1,-1,-1,1)$,
$x_{124} = (1, -1,-1,-1,1,-1,$ $1,1,1)$,
$x_{133} = (1, -1,-1,-1,-1,1,-1,-1,1)$,
$x_{134} = (1, -1,-1$, $-1,-1,1,1,1,1)$,
$x_{223} = (-1, 1,-1,-1,1,-1,-1,-1,1)$,
$x_{224} = (-1, 1,-1,-1,1,-1,1,1,1)$,
$x_{233} = (-1, 1,-1,-1,-1,$ $1,-1,-1,1)$,
$x_{234} = (-1, 1,-1,-1,-1,1,1,1,1)$.

Observe that $X^2_i[x_h\times x_2x_3 \times x_4x_3]=
x_{h24}x_{h33} \ast S^{\,0}_2(x_{h23}, x_{h34})$ and $X^2_i[x_1x_2
\times x_2x_3 \times x_k] = x_{13k}x_{22k}\ast S^{\,0}_2(x_{12k},
x_{23k})$ for $1\leq h \leq 2$, $3\leq k \leq 4$, $1\leq i\leq 4$
and $X^2[x_1x_2 \times x_l \times x_4x_3]= x_{1l4}x_{2l3}\ast
S^{\,0}_2(x_{1l3}, x_{2l4})$ for $2\leq l \leq 3$. Since
$x_{224}x_{331}$ is an edge (in $x_2x_3\times x_2x_3 \times
x_4x_1$), by condition $(b)$, $x_{224}x_{133}$ is not an edge. So,
no new edge (with vertices in $C^3_6$) can be added to triangulate
$C^3_6$. We introduce the new vertex $u = (0, 0, -1, -1, 0, 0, 0,
0, 1)$ in the interior of $C^3_6$ and by the same argument as for
$C^3_5$,
$$
\widehat {C^3_6} := u\ast S^{\,0}_2(x_{133}, x_{224}) \ast
S^{\,0}_2(x_{124}, x_{233}) \ast S^{\,0}_2(x_{223}, x_{134}) \cup
\{x_{123}x_{133}x_{124}x_{223}, x_{234}x_{224} x_{233} x_{134}\}
$$
is a subdivision of $C^3_6$.

\setlength{\unitlength}{3.9mm}

\begin{picture}(36,19.5)(0,-2.2)

\thicklines

\put(2,3){\line(1,0){8}} \put(2,3){\line(0,1){8}}
\put(2,11){\line(1,0){8}} \put(10,3){\line(0,1){8}}
\put(6,15){\line(1,0){8}} \put(14,7){\line(0,1){8}}
\put(6,7){\line(1,0){3.5}} \put(10.5,7){\line(1,0){3.5}}
\put(6,7){\line(0,1){3.5}} \put(6,11.5){\line(0,1){3.5}}

\put(2,3){\line(1,1){4}} \put(2,11){\line(1,1){4}}
\put(10,3){\line(1,1){4}} \put(10,11){\line(1,1){4}}

\put(18,6){\line(1,0){8}} \put(18,6){\line(3,5){6}}

\put(22,10){\line(1,3){2}} \put(24,0){\line(1,3){2}}
\put(24,0){\line(3,5){6}}  \put(22,10){\line(1,0){2.7}}
\put(25.6,10){\line(1,0){0.8}} \put(30,10){\line(-1,0){2.8}}

\put(24,0){\line(-2,5){2}} \put(18,6){\line(4,-1){4}}
\put(22,10){\line(0,-1){3.7}}  \put(22,5){\line(0,1){0.6}}
\put(26,6){\line(1,5){1}} \put(30,10){\line(-3,1){3}}
\put(27,11){\line(-3,5){3}}

\thinlines

\put(-1,-2.2){\line(1,0){38}} \put(-1,16.5){\line(1,0){38}}

\put(18,6){\line(1,1){4}} \put(24,0){\line(-1,1){6}}
\put(24,0){\line(-1,5){1.1}} \put(22,10){\line(1,-5){0.7}}
 \put(26,6){\line(1,1){4}}  \put(26,6){\line(-1,5){2}}
\put(30,10){\line(-1,1){6}}

\put(6,15){\line(1,-1){4}} \put(6,7){\line(1,-1){4}}
\put(6,7){\line(-1,1){4}} \put(14,7){\line(-1,1){4}}

\put(3.3,3.6){\mbox{$x_{123}$}} \put(11,3){\mbox{$x_{223}$}}
\put(6.5,7.6){\mbox{$x_{133}$}} \put(11,7.6){\mbox{$x_{233}$}}
\put(3.2,11.5){\mbox{$x_{124}$}} \put(8,10){\mbox{$x_{224}$}}
\put(3.5,14.7){\mbox{$x_{134}$}} \put(11,14.2){\mbox{$x_{234}$}}

\put(25,15.5){\mbox{$x_{224}$}} \put(25,0.2){\mbox{$x_{133}$}}
\put(22.5,10.6){\mbox{$x_{124}$}} \put(23.7,6.6){\mbox{$x_{233}$}}
\put(30.2,9.2){\mbox{$x_{134}$}} \put(16.3,6.7){\mbox{$x_{223}$}}
\put(20,4.4){\mbox{$x_{123}$}} \put(25.1,10.8){\mbox{$x_{234}$}}

\put(2,0.7){\mbox{\bf (a)\,: \boldmath{$C^3_6 = x_1x_2\times
x_2x_3 \times x_3x_4$}}}

\put(18,-1.2){\mbox{\bf (b)\,: Boundary of
\boldmath{$\widehat{C^3_6} = X^2[C^3_6]$}}}

\put(5,-1){\mbox{\bf Figure 6}}

\end{picture}

\bigskip

By the action of $S_3 \times A_4$, we get the subdivisions of all
the 3-polytopes of the type $x_ix_j \times x_jx_k \times x_kx_l$,
where $i, j, k, l$ are distinct. Since $(1,3)\times (1,4)(2,3) \in
S_3 \times A_4$ is the stabilizer of $C^3_6$ and $\widehat
{C^3_6}$, any 3-polytope of this type can be written as $(\alpha,
\beta)(C^3_6)$, where $(\alpha, \beta) \in A_3 \times A_4$.
Accordingly, there are 36 new vertices, namely,
$u^{\alpha}_{\beta} := (\alpha, \beta)(u)$, where $(\alpha, \beta)
\in A_3 \times A_4$.

Alternatively, let us label the elements of $S_3$ and $A_4$ as
follows. $S_3 = \{\alpha_1 = {\rm identity}, \alpha_2 =(123),
\alpha_3 = (132), \alpha_4= (13), \alpha_5=(23), \alpha_6=(12)\}$
and $A_4 = \{\beta_1 = {\rm identity}, \beta_{2} = (123),
\beta_{3} = (124), \beta_{4} = (132), \beta_{5} = (134), \beta_{6}
= (12)(34), \beta_{7} = (14)(23), \beta_{8} = (142), \beta_{9} =
(234), \beta_{10} = (143), \beta_{11} = (243), \beta_{12} =
(13)(24)\}$. Let $u_{ij} = u^{\alpha_i}_{\beta_j}$. Check that
$u^{\alpha_i}_{\beta_j} = u^{\alpha_{i+3}}_{\beta_{j -6}}$, for $7
\leq j \leq 12$ (summation in the subscripts of $\alpha$ is modulo
6). Thus, $\{u^{\alpha}_{\beta} \, : \, (\alpha, \beta) \in S_3
\times A_4\} = \{u^{\alpha}_{\beta} \, : \, (\alpha, \beta) \in
A_3 \times A_4\} = \{u_{ij} \, : \, 1 \leq i, j\leq 6\}$.

This gives a simplicial subdivision $X^3$ of $W^3$ by
adding $24 + 36 = 60$ new vertices. Let $V(X^3)$ and $E(X^3)$ be
vertex-set and edge-set of $X^3$ respectively. Then
\begin{eqnarray}
V(X^3) \!\!\! & = \!\!\!& \{x_{ijk} ~ :
1\leq i, j, k \leq 4\} \cup \{u_{ij} ~ : ~ 1\leq i, j
\leq 6\} \nonumber \\
&&   \cup \, \{v_{pqr} ~ : ~ 1\leq p, q, r \leq 4,
\, p, q, r \mbox{ are distinct}\}, \nonumber \\
E(X^3) \!\!\! & =  \!\!\! &E(X^2) \,
\cup \, (S_3\times A_4)(\{x_{111}x_{222}, x_{111}x_{223},
x_{111}x_{234}, x_{112}x_{334}\} \cup   \nonumber \\
&& \hspace{10mm} \{v_{321}x_{112}, v_{321}x_{113}\}
\cup \{u_{11}x_{124}, u_{11}x_{133},
u_{11}x_{223}\}) \subseteq E(X), \label{EX3}
\end{eqnarray}
where $E(X^2)$ and $E(X)$ are as in equations (\ref{EX2}) and
(\ref{edges}) respectively.

Since $E(X^3) \subseteq E(X)$, by the same argument as in the
sketch of the proof of Theorem \ref{T1}, the $S_3$-action on $X^3$
is good. \hfill $\Box$

\begin{lemma}$\!\!\!${\bf .} \label{R5}
Let $K$ be the simplicial complex whose
facets are the following $24$ $\,4$-simplices\,: \vspace{-6mm}
\begin{eqnarray*}
&& u_{64}v_{321}x_{112}x_{113}x_{133},
u_{64}v_{321}x_{112}x_{113}x_{212},
u_{64}v_{321}x_{112}x_{232}x_{133},
u_{64}v_{321}x_{112}x_{232}x_{212}, \\
&& u_{64}v_{321}x_{233}x_{113}x_{133},
u_{64}v_{321}x_{233}x_{113}x_{212},
u_{64}v_{321}x_{233}x_{232}x_{133},
u_{64}v_{321}x_{233}x_{232}x_{212}, \\
&& u_{64}x_{113}x_{114}x_{133}x_{233},
u_{64}x_{114}x_{133}x_{134}x_{233},
u_{64}x_{114}x_{134}x_{233}x_{234},
u_{64}x_{112}x_{113}x_{114}x_{212}, \\
&& u_{64}x_{133}x_{134}x_{232}x_{233},
u_{64}x_{134}x_{232}x_{233}x_{234},
u_{64}x_{112}x_{113}x_{114}x_{133},
u_{64}x_{112}x_{114}x_{133}x_{134}, \\
&& u_{64}x_{212}x_{232}x_{233}x_{234},
u_{64}x_{112}x_{133}x_{134}x_{232},
u_{64}x_{113}x_{114}x_{212}x_{233},
u_{64}x_{114}x_{212}x_{233}x_{234}, \\
&& x_{132}x_{112}x_{133}x_{134}x_{232},
x_{113}x_{114}x_{212}x_{213}x_{233},
x_{114}x_{212}x_{213}x_{214}x_{233},
x_{114}x_{212}x_{214}x_{233}x_{234}.
\end{eqnarray*}
Then $K$ is a simplicial subdivision of the $4$-polytope $C^4 =
x_1x_2 \times x_1x_3 \times x_2x_3x_4$, the boundary $\partial
\overline{K}$ of $K$ is a subcomplex of $X^3$, $K$ contains four
more edges than those in $\partial \overline{K}$, namely,
$u_{64}x_{113}, u_{64}x_{133}$, $u_{64}x_{233}$,  $u_{64}v_{321}$
and there is no vertex adjacent to both $x_{132}$ and $x_{213}$ in
$K$.
\end{lemma}

\noindent {\bf Proof.} The 3-dimensional faces of $C^4$ are
 $F_1 = x_1x_2 \times x_1x_3 \times x_2x_4$,
 $F_2 = x_1x_2 \times x_1x_3 \times x_2x_3$,
 $F_3 = x_1x_2 \times x_1x_3 \times x_3x_4$,
 $F_4 = x_1x_2 \times x_1 \times x_2x_3x_4$,
 $F_5 = x_1x_2 \times x_3 \times x_2x_3x_4$,
 $F_6 = x_1 \times x_1x_3 \times x_2x_3x_4$ and
 $F_7 = x_2 \times x_1x_3 \times x_2x_3x_4$.
For $1\leq i\leq 7$, let $X^3[F_i]$ be the triangulation of $F_i$
induced by $X^3$. Then the interior point $u_{64}$ of $F_1$ is a
vertex of $X^3[F_1]$ and the interior point $v_{321}$ of $F_2$ is
a vertex of $X^3[F_2]$. By Lemma \ref{L2.2}, the 4-polytopes
$u_{64}\ast F_2, \dots, u_{64}\ast F_7$ give a polytopal
subdivision $\widetilde{C^4}$ of $C^4$ and $K_1 := \{u_{64}\ast
\alpha \, : \, \alpha\in X^3[F_2] \cup \cdots \cup X^3[F_7]\}$ is
a simplicial subdivision of $C^4$. Now, from $X^3$, we have\,:
\begin{eqnarray*}
u_{64}\ast X^3[F_2] & = & u_{64}v_{321}\ast S^0_2(x_{112},
x_{233}) \ast S^0_2(x_{113}, x_{232})\ast S^0_2(x_{133}, x_{212}) \\
&& ~ \cup\{u_{64}x_{132}x_{112}x_{232}x_{133}, u_{64}x_{213}
x_{233} x_{113}x_{212}\}, \\ u_{64}\ast X^3[F_3] & = & u_{64}
x_{114}x_{233}\ast
\begin{picture}(8,2.5)(-1,1)
\setlength{\unitlength}{2.5mm}
\put(1.7,0){$_{\bullet}$}
\put(1.7,3.85){$_{\bullet}$}
\put(4.8,0){$_{\bullet}$}
\put(4.8,3.85){$_{\bullet}$}
\put(-0.1,1.9){$_{\bullet}$}
\put(6.7,1.9){$_{\bullet}$}
\thicklines
\put(2,0){\line(1,0){3}}
\put(2,0){\line(-1,1){2}}
\put(5,0){\line(1,1){2}}
\put(5,4){\line(1,-1){2}}
\put(2,4){\line(1,0){3}}
\put(0,2){\line(1,1){2}}
\put(-1.5,3.9){\small {$x_{134}$}}
\put(1,1.8){\small {$x_{133}$}}
\put(-1.5,0){\small {$x_{113}$}}
\put(6.7,0){\small {$x_{213}$}}
\put(7.7,1.8){\small {$x_{214}$}}
\put(6.2,3.9){\small {$x_{234}$}}
\end{picture}, \\
&&\\
u_{64}\ast X^3[F_4] & = & u_{64}x_{114}x_{212}\ast
\begin{picture}(13,1)(0,0)
\setlength{\unitlength}{2.5mm}
\put(0.3,0.2){$_{\bullet}$}
\put(5.8,0.2){$_{\bullet}$}
\put(11.3,0.2){$_{\bullet}$}
\put(16.8,0.2){$_{\bullet}$}
\thicklines
\put(0.5,0.3){\line(1,0){16.5}}
\put(0.2,0.9){\small {$x_{112}$}}
\put(5.7,0.9){\small {$x_{113}$}}
\put(11.2,0.9){\small {$x_{213}$}}
\put(16.7,0.9){\small {$x_{214}$}}
\end{picture},
 \\
u_{64}\ast X^3[F_5] & = & u_{64}x_{134}x_{232}\ast
\begin{picture}(13,1)(0,0)
\setlength{\unitlength}{2.5mm}
\put(0.3,0.2){$_{\bullet}$}
\put(5.8,0.2){$_{\bullet}$}
\put(11.3,0.2){$_{\bullet}$}
\put(16.8,0.2){$_{\bullet}$}
\thicklines
\put(0.5,0.3){\line(1,0){16.5}}
\put(0.2,0.9){\small {$x_{132}$}}
\put(5.7,0.9){\small {$x_{133}$}}
\put(11.2,0.9){\small {$x_{233}$}}
\put(16.7,0.9){\small {$x_{234}$}}
\end{picture}, \\
u_{64}\ast X^3[F_6] & = & u_{64}x_{112}x_{133}\ast
\begin{picture}(13,1)(0,0)
\setlength{\unitlength}{2.5mm}
\put(0.3,0.2){$_{\bullet}$}
\put(5.8,0.2){$_{\bullet}$}
\put(11.3,0.2){$_{\bullet}$}
\put(16.8,0.2){$_{\bullet}$}
\thicklines
\put(0.5,0.3){\line(1,0){16.5}}
\put(0.2,0.9){\small {$x_{113}$}}
\put(5.7,0.9){\small {$x_{114}$}}
\put(11.2,0.9){\small {$x_{134}$}}
\put(16.7,0.9){\small {$x_{132}$}}
\end{picture}, \\
u_{64}\ast X^3[F_7] & = & u_{64}x_{212}x_{233}\ast
\begin{picture}(13,1)(0,0)
\setlength{\unitlength}{2.5mm}
\put(0.3,0.2){$_{\bullet}$}
\put(5.8,0.2){$_{\bullet}$}
\put(11.3,0.2){$_{\bullet}$}
\put(16.8,0.2){$_{\bullet}$}
\thicklines
\put(0.5,0.3){\line(1,0){16.5}}
\put(0.2,0.9){\small {$x_{213}$}}
\put(5.7,0.9){\small {$x_{214}$}}
\put(11.2,0.9){\small {$x_{234}$}}
\put(16.7,0.9){\small {$x_{232}$}}
\end{picture}.
\end{eqnarray*}
Observe that the induced complex $K_1[F_1] \neq X^3[F_1]$.

Now, $u_{64} = \frac{1}{2}(x_{114}+x_{232}) = (0, 0, -1, 0, -1, 0,
0, 1, 0)$. Thus, $\frac{2}{3} u_{64} + \frac{1}{3} x_{132} =
\frac{1}{3}x_{112} + \frac{1}{3} x_{134} + \frac{1}{3} x_{232} =
\frac{1}{3}(1, -1, -3, -1, -3, 1, -1, 3,-1)$. Thus $\intalpha \cap
\intbeta \neq \emptyset$, where $\alpha = u_{64}x_{132}$ and
$\beta = x_{112} x_{134}x_{232}$. Since
$x_{133}u_{64}x_{132}x_{112}x_{134}$ is a simplex, by Lemma
\ref{L2.3}, $\overline{x_{133}}\ast \partial \overline{\alpha}
\ast \overline{\beta}$ is a simplicial complex and
$|\overline{x_{133}}\ast \overline{\alpha} \ast
\partial \overline{\beta}| = |\overline{x_{133}}\ast \partial
\overline{\alpha} \ast \overline{\beta}|$. So, we can
replace $\overline{x_{133}}\ast \overline{\alpha} \ast
\partial \overline{\beta}$ ($= \overline{x_{133}u_{64} x_{132}}
\ast S^{1}_3(\{x_{112}, x_{134}, x_{232}\})$) by
$\overline{x_{133}}\ast \partial \overline{\alpha} \ast
\overline{\beta}$ ($=\overline{x_{133}x_{112}x_{134}x_{232}} \ast
S^{\,0}_2(u_{64}, x_{132})$) in $K_1$ to get a new triangulation
$K_2$ of $C^4$.

Again, $\frac{2}{3}u_{64} + \frac{1}{3} x_{213} = \frac{1}{3}
x_{114} + \frac{1}{3} x_{212} + \frac{1}{3} x_{233} =
\frac{1}{3}(-1, 1, -3, 1, -3, -1, -1, 1,1)$ and
$x_{113}u_{64}x_{213}x_{114}x_{212}$, $x_{214}u_{64}x_{213}
x_{114} x_{212}$ are simplices. So, by the same argument as above,
we can replace $\overline{x_{113} u_{64}x_{213}} \ast
S^1_3(\{x_{114}, x_{212},$ $x_{233}\})$ by $\overline{x_{113}}
\ast S^{\,0}_2(u_{64}, x_{213}) \ast   \overline{x_{114}
x_{212}x_{233}}$ and $\overline{x_{214}u_{64}x_{213}} \ast
S^1_3(\{x_{114}, x_{212}, x_{233}\})$ by $\overline{x_{214}} \ast
S^{\,0}_2(u_{64}, x_{213})\ast \overline{x_{114}x_{212}x_{233}}$
in $K_2$ to get a new triangulation $K_3$ of $C^4$. So, $K_3$ is
obtained from $K_2$ by replacing $S^{\,0}_2(x_{214}, x_{113}) \ast
\overline{u_{64}x_{213}} \ast S^1_3(\{x_{114}, x_{212},
x_{233}\})$ by $S^{\,0}_2(x_{113}, x_{214}) \ast S^{\,0}_2(u_{64},
x_{213})\ast \overline{x_{114}x_{212}x_{233}}$.

Finally, $\frac{2}{3}u_{64} + \frac{1}{3} x_{214} =  \frac{1}{3}
x_{114} + \frac{1}{3} x_{212} + \frac{1}{3} x_{234}$ and $x_{233}
u_{64}x_{214}x_{114}x_{212}$ is a simplex in the simplicial
subdivision $K_3$ of $C^4$. So, by the same argument as before, we
can replace $\overline{x_{233}u_{64}x_{214}} \ast S^1_3(\{x_{114},
x_{212}, x_{234}\})$ by $\overline{x_{233}} \ast S^{\,0}_2(u_{64},
x_{214}) \ast \overline{x_{114}x_{212} x_{234}}$ in $K_3$ to get a
new triangulation $K_4$ of $C^4$. Observe that $K_4 = K$.

It is easy to check that $\partial \overline{K} = X^3[F_1] \cup
\cdots  \cup X^3[F_7]$ and the new edges in $K$ are
$u_{64}x_{113}, u_{64}x_{133}$, $u_{64}x_{233}$ and
$u_{64}v_{321}$. This proves the lemma. \hfill $\Box$

\begin{lemma}$\!\!\!${\bf .} \label{R6}
Let $X^3$, $W^j$ $(0\leq j \leq 6)$ and $V_{124}$ be as in  Lemma
$\ref{R4}$. There exists a simplicial complex $X^4$ which satisfy
the following\,: {\rm (i)} $X^4$ is a simplicial subdivision of
$W^4$, {\rm (ii)} the $3$-skeleton of $X^4$ is $X^3$, {\rm (iii)}
$S_3 \times A_4$ acts as an isometry group, where the action of
$S_3$ and $A_4$ on the vertices is as in Lemma $\ref{R4}$, {\rm
(iv)} the induced $S_3$ $($respectively, $A_4)$ action on $|X^4|$
is same as that on $|W^4|$ and {\rm (v)} the $S_3$ action on $X^4$
is good.
\end{lemma}

\noindent {\bf Proof.} A 4-polytope in $W^4$ is the product of two
2-simplices  and one 0-simplex or product of one 2-simplex and two
1-simplices. Up to the action of $S_3\times A_4$, the 4-polytopes
are of the form $C^4_{1,i} = x_i \times x_1x_2x_3 \times
x_1x_2x_3$, $C^4_{2,i} = x_i \times x_1x_2x_3 \times x_1x_2x_4$,
$C^4_3 = x_1x_2 \times x_1x_2 \times x_1x_2x_3$, $C^4_4 = x_1x_2
\times x_1x_2 \times x_1x_3x_4$, $C^4_5 = x_1x_2 \times x_1x_3
\times x_1x_2x_3$, $C^4_6 = x_1x_2 \times x_1x_3 \times
x_1x_3x_4$, $C^4_7 = x_1x_2 \times x_1x_4 \times x_1x_3x_4$,
$C^4_8 = x_1x_2 \times x_1x_2x_3 \times x_3x_4$, $C^4_9 = x_1x_2
\times x_1x_3 \times x_2x_3x_4$, $1\leq i\leq 4$. For $1\leq j\leq
9$, we know the simplicial subdivision $X^3[\partial
\overline{C^4_j}]$ of the boundary complex $\partial
\overline{C^4_j}$ of $C^4_j$.

The 4-polytope $C^4_{1,i}$ has six 3-faces and only two of them,
namely $C_{11} := x_i \times x_2x_3 \times x_1x_2x_3$ and $C_{12}
:= x_i \times x_1x_2x_3 \times x_2x_3$ do not contain the vertex
$x_{i11}$. So, by Lemma \ref{L2.2}, $\widehat{C}^4_{1,i} :=
\{x_{i11} \ast \alpha \, : \, \alpha \in X^3[{C_{11}}] \cup
X^3[{C_{12}}]\}$ is a simplicial subdivision of $C^4_{1,i}$. Since
the anti-star ${\rm ast}_{X^3[C^4_{1,i}]}(x_{i11})$ of $x_{i11}$
in $X^3[\partial \overline{C^4_{1,i}}] = X^3[C^4_{1,i}]$ is same
as $X^3[C_{11}] \cup X^3[C_{12}]$, it follows that
$\widehat{C}^4_{1,i} = x_{i11} \ast {\rm
ast}_{X^3[C^4_{1,i}]}(x_{i11})$. Since $X^3[D] = x_{i11} \ast {\rm
ast}_{X^2[D]}(x_{i11})$, it follows that $X^3[D] =
\widehat{C}^4_{1,i}[D]$ for each 3-face $D$ of $C^4_{1,i}$
containing $x_{i11}$. Thus, $X^3[{C^4_{1,i}}] = \partial
\widehat{C}^4_{1,i}$. We take $\widehat{C}^4_{1,i}$ is the
simplicial subdivision of $C^4_{1,i}$. More explicitly we have\,:

\begin{eqnarray}
\widehat{C}^4_{1,i} & = & \{x_{i11}x_{i22}x_{i33}x_{i12}x_{i13},
x_{i11}x_{i22}x_{i33}x_{i13}x_{i23},
x_{i11}x_{i22}x_{i33}x_{i23}x_{i21}, \nonumber \\
&& x_{i11}x_{i22}x_{i33}x_{i21}x_{i31},
x_{i11}x_{i22}x_{i33}x_{i31}x_{i32},
x_{i11}x_{i22}x_{i33}x_{i32}x_{i12}\}. \label{C41}
\end{eqnarray}

Similarly, if $\widehat{C}^4_{2,i} := x_{i11} \ast {\rm
ast}_{X^3[C^4_{2,i}]}(x_{i11})$ and $\widehat{C}^4_8 := x_{114}
\ast {\rm ast}_{X^3[C^4_8]}(x_{114})$ then $\widehat{C}^4_{2,i}$
and $\widehat{C}^4_8$ are simplicial subdivisions of the polytopes
$C^4_{2,i}$ and $C^4_8$ respectively with the property that
$\partial \widehat{C}^4_{2,i} = X^3[C^4_{2,i}]$ and $\partial
\widehat{C}^4_8 = X^3[C^4_8]$. We take $\widehat{C}^4_{2,i}$ and
$\widehat{C}^4_8$ as the simplicial subdivisions of $C^4_{2,i}$
and $C^4_8$ respectively. Observe that, we are adding the new edge
$u_{11}x_{114}$ in $\widehat{C}^4_8$. Easy to see the following\,:
\begin{eqnarray}
\widehat{C}^4_{2,i} & = & \{x_{i11}x_{i22}x_{i14}x_{i24}x_{i34},
x_{i11}x_{i22}x_{i24}x_{i31}x_{i34},
x_{i11}x_{i22}x_{i21}x_{i24}x_{i31}, \nonumber \\
&& ~~~ x_{i11}x_{i22}x_{i31}x_{i32}x_{i34},
x_{i11}x_{i22}x_{i12}x_{i32}x_{i34},
x_{i11}x_{i22}x_{i12}x_{i14}x_{i34}\}. \\ \label{C42}
\widehat{C}^4_8  &= & \{x_{114}x_{123}x_{124}x_{133}x_{223},
x_{114}x_{134}x_{224}x_{233}x_{234},
x_{114}x_{133}x_{124}x_{134}u_{11}, \nonumber \\
&& ~~~ x_{114}x_{133}x_{124}x_{223}u_{11}, x_{114}x_{133}x_{233}x_{134}u_{11},
x_{114}x_{133}x_{233}x_{223}u_{11}, \nonumber \\
&& ~~~ x_{114}x_{224}x_{124}x_{134}u_{11}, x_{114}x_{224}x_{124}x_{223}u_{11},
x_{114}x_{224}x_{233}x_{134}u_{11}, \nonumber \\
&& ~~~ x_{114}x_{224}x_{233}x_{223}u_{11}, x_{114}x_{113}x_{223}x_{123}x_{133},
x_{114}x_{113}x_{223}x_{133}x_{233},  \nonumber \\
&& ~~~ x_{114}x_{113}x_{223}x_{233}x_{213},  x_{114}x_{214}x_{233}x_{213}x_{223},
x_{114}x_{214}x_{233}x_{223}x_{224},  \nonumber \\
&& ~~~ x_{114}x_{214}x_{233}x_{224}x_{234}\}. \label{C48}
\end{eqnarray}

Similarly, $\widehat{C}^4_j := x_{111} \ast {\rm
ast}_{X^3[C^4_j]}(x_{111})$ is a simplicial subdivision of the
4-polytope $C^4_j$ with $\partial \widehat{C}^4_j = X^3[C^4_j]$,
for $3\leq j \leq 7$. We take these subdivisions (for $C^4_3,
\dots, C^4_7$ respectively) together with the simplicial complex
$K$ in Lemma \ref{R5} as the simplicial subdivision of $C^4_9$.

By the action of $S_3 \times A_4$, we get simplicial subdivisions
of all the 4-polytopes in $W^4$ and get $X^4$. Since $v_{321} \in
C^4_5$, $u_{16} \in C^4_7$ and $x_{233} \in C^4_8$, the edge-set
of $X^4$ is the following\,:
\begin{eqnarray}
E(X^4) \!\!\! & =  \!\!\!& E(X^3) \, \cup \,
(S_3\times A_4)(\{x_{111}u_{16}, x_{111}v_{321}, u_{11}x_{114},
u_{64}x_{113}, \nonumber \\
&& \hspace{37mm} u_{64}x_{133}, u_{64}x_{233}, u_{64}v_{321},
x_{114}x_{233}\}) \nonumber \\
&=  \!\!\!& E(X^3) \, \cup \, (S_3\times A_4)(\{x_{111}u_{12},
x_{111}v_{123}, u_{11}x_{114}, u_{11}x_{221}, \nonumber \\
&& \hspace{37mm} u_{11}x_{121}, u_{11}x_{131}, u_{11}v_{312},
x_{112}x_{344}\}) \subseteq  E(X), \label{EX4}
\end{eqnarray}
where $E(X^3)$ and $E(X)$ are as in equations (\ref{EX3}) and
(\ref{edges}) respectively.

Since $E(X^4) \subseteq E(X)$, by the same argument as in the
sketch of the proof of Theorem \ref{T1}, the $S_3$-action on $X^4$
is good. \hfill $\Box$


\begin{lemma}$\!\!\!${\bf .} \label{R7}
Let $L$ be the simplicial complex whose
facets are the following $75$ $\,5$-simplices\,: \vspace{-6mm}
\begin{eqnarray*}
   x_{114}x_{144}x_{113}x_{112}u_{16}u_{64},
&  x_{114}x_{144}x_{112}v_{421}u_{16}u_{64},
&  x_{114}x_{144}v_{421}x_{244}u_{16}u_{64},   \\
   x_{114}x_{212}x_{113}x_{112}u_{16}u_{64},
&  x_{114}x_{212}x_{112}v_{421}u_{16}u_{64},
&  x_{114}x_{212}v_{421}x_{244}u_{16}u_{64},   \\
   x_{242}x_{144}x_{113}x_{112}u_{16}u_{64},
&  x_{242}x_{144}x_{112}v_{421}u_{16}u_{64},
&  x_{242}x_{144}v_{421}x_{244}u_{16}u_{64},   \\
   x_{242}x_{212}x_{113}x_{112}u_{16}u_{64},
&  x_{242}x_{212}x_{112}v_{421}u_{16}u_{64},
&  x_{242}x_{212}v_{421}x_{244}u_{16}u_{64},   \\
   x_{113}x_{112}x_{242}x_{212}u_{51}u_{64},
&  x_{113}x_{112}x_{212}v_{321}u_{51}u_{64},
&  x_{113}x_{112}v_{321}x_{133}u_{51}u_{64},   \\
   x_{113}x_{233}x_{242}x_{212}u_{51}u_{64},
&  x_{113}x_{233}x_{212}v_{321}u_{51}u_{64},
&  x_{113}x_{233}v_{321}x_{133}u_{51}u_{64},   \\
   x_{232}x_{112}x_{242}x_{212}u_{51}u_{64},
&  x_{232}x_{112}x_{212}v_{321}u_{51}u_{64},
&  x_{232}x_{112}v_{321}x_{133}u_{51}u_{64},   \\
   x_{232}x_{233}x_{242}x_{212}u_{51}u_{64},
&  x_{232}x_{233}x_{212}v_{321}u_{51}u_{64},
&  x_{232}x_{233}v_{321}x_{133}u_{51}u_{64},   \\
   x_{114}x_{144}x_{233}x_{244}u_{16}u_{64},
&  x_{114}x_{144}x_{233}x_{113}u_{16}u_{64},
&  x_{114}x_{144}x_{233}x_{244}x_{234}u_{64},  \\
   x_{114}x_{144}x_{233}x_{234}x_{134}u_{64},
&  x_{114}x_{144}x_{233}x_{134}x_{133}u_{64},
&  x_{114}x_{144}x_{233}x_{133}x_{113}u_{64},  \\
   x_{144}x_{232}x_{233}x_{242}u_{51}u_{64},
&  x_{144}x_{232}x_{233}x_{133}u_{51}u_{64},
 &  x_{144}x_{232}x_{233}x_{242}x_{244}u_{64}, \\
  x_{144}x_{232}x_{233}x_{244}x_{234}u_{64},
&  x_{144}x_{232}x_{233}x_{234}x_{134}u_{64},
&  x_{144}x_{232}x_{233}x_{134}x_{133}u_{64}, \\
  x_{144}x_{114}x_{112}x_{133}x_{134}u_{64}, &
  x_{144}x_{114}x_{112}x_{133}x_{113}u_{64}, &
  x_{244}x_{212}x_{232}x_{233}x_{234}u_{64},\\
  x_{244}x_{212}x_{232}x_{233}x_{242}u_{64}, &
  x_{144}x_{242}x_{243}x_{233}x_{143}x_{113}, &
  x_{144}x_{233}x_{113}x_{242}u_{51}u_{64}, \\
  x_{144}x_{233}x_{113}x_{133}u_{51}u_{64}, &
  x_{144}x_{143}x_{113}x_{242}x_{142}u_{51}, &
  x_{144}x_{143}x_{113}x_{242}x_{233}u_{51}, \\
  x_{144}x_{143}x_{113}x_{133}x_{142}u_{51}, &
  x_{144}x_{143}x_{113}x_{133}x_{233}u_{51}, &
  x_{144}x_{242}x_{233}x_{244}u_{16}u_{64}, \\
  x_{144}x_{242}x_{233}x_{113}u_{16}u_{64}, &
  x_{242}x_{243}x_{233}x_{244}x_{144}u_{16}, &
  x_{242}x_{243}x_{233}x_{244}x_{213}u_{16}, \\
  x_{242}x_{243}x_{233}x_{113}x_{144}u_{16}, &
  x_{242}x_{243}x_{233}x_{113}x_{213}u_{16}, &
  x_{144}x_{232}x_{112}x_{133}x_{134}u_{64},\\
  x_{144}x_{232}x_{112}x_{132}x_{133}x_{134}, &
  x_{144}x_{232}x_{112}x_{132}x_{133}x_{142}, &
  x_{244}x_{114}x_{212}x_{233}x_{234}u_{64}, \\
  x_{244}x_{214}x_{114}x_{212}x_{233}x_{234}, &
  x_{244}x_{214}x_{114}x_{212}x_{233}x_{213}, &
  x_{144}x_{112}x_{232}x_{242}u_{51}u_{64}, \\
  x_{144}x_{112}x_{232}x_{133}u_{51}u_{64}, &
  x_{144}x_{112}x_{113}x_{242}u_{51}u_{64}, &
  x_{144}x_{112}x_{113}x_{133}u_{51}u_{64}, \\
  x_{212}x_{233}x_{114}x_{244}u_{16}u_{64}, &
  x_{212}x_{233}x_{114}x_{113}u_{16}u_{64}, &
  x_{212}x_{233}x_{242}x_{244}u_{16}u_{64}, \\
  x_{212}x_{233}x_{242}x_{113}u_{16}u_{64}, &
  x_{212}x_{213}x_{233}x_{114}x_{244}u_{16}, &
  x_{212}x_{213}x_{233}x_{114}x_{113}u_{16}, \\
  x_{212}x_{213}x_{233}x_{242}x_{244}u_{16}, &
  x_{212}x_{213}x_{233}x_{242}x_{113}u_{16}, &
  x_{144}x_{142}x_{112}x_{232}x_{242}u_{51}, \\
  x_{144}x_{142}x_{112}x_{232}x_{133}u_{51}, &
  x_{144}x_{142}x_{112}x_{113}x_{242}u_{51}, &
  x_{144}x_{142}x_{112}x_{113}x_{133}u_{51}.
\end{eqnarray*}
Then $L$ is a simplicial subdivision of the $5$-polytope $D^5 =
x_1x_2 \times x_1x_3x_4 \times x_2x_3x_4$, the boundary $\partial
L$ of $L$ is the induced subcomplex $X^4[L]$ of $X^4$, $L$
contains two more edges than those in $\partial L$, namely,
$u_{64}u_{16}$ and $u_{64}u_{51}$.
\end{lemma}

\noindent {\bf Proof.} Observe that $D^5$ is invariant under the
automorphism $f = (2, 3) \times (1, 2)(3, 4)\in S_3 \times A_4$.
The polytope $D_5$ has 23 vertices, the vertex $u_{64}$ is fixed
by $f$ and the action of $f$ on other vertices is given by: $f
\equiv (u_{16}, u_{51})(v_{321}, v_{421})(x_{112}, x_{212})
(x_{113}, x_{242})(x_{114}, x_{232}) (x_{132}, x_{214})$
$(x_{133}, x_{244})(x_{134}, x_{234}) (x_{142}, x_{213})(x_{143},
x_{243})(x_{144}, x_{233})$. Here $u_{64} = (0, 0, -1, 0, -1, 0,
0, 1, 0)$ and $u_{51} = \frac{1}{2}(x_{133}+x_{242}) = (0, 0, -1,
0, 0, 1, -1, 0, 0)$.

The 4-dimensional faces of $D^5$ are
 $D_1 = x_1x_2 \times x_1x_3 \times x_2x_3x_4 = C^4_9$,
 $D_2 = x_1x_2 \times x_1x_3x_4 \times x_2x_4  = f(C^4_9)$,
 $D_3 = x_1x_2 \times x_1x_4 \times x_2x_3x_4
= ((1,2)\times (2,4,3))(C^4_9)$,
 $D_4 = x_1x_2 \times x_1x_3x_4 \times x_2x_3
= f(D_3) = ((2,1,3)\times (1,2,3))(C^4_9)$,
 $D_5 = x_1x_2 \times x_1x_3x_4 \times x_3x_4
= ((1,3)\times (1,3)(2,4))(C^4_8)$,
 $D_6 = x_1x_2 \times x_3x_4 \times x_2x_3x_4
= f(D_5) = ((1,2,3)\times (1,4)(2,3))(C^4_8)$,
 $D_7 = x_1 \times x_1x_3x_4 \times x_2x_3x_4 \cong C^4_{2,i}$ and
 $D_8 = x_2 \times x_1x_3x_4 \times x_2x_3x_4 = f(D_7)$,
where $C^4_j$ are as in the proof of Lemma \ref{R6}. For $1\leq
i\leq 8$, let $X^4[D_i]$ be the triangulation of $D_i$ induced by
$X^4$. Then the interior point $u_{64}$ of $D_j$ is a vertex of
$X^4[D_j]$ for $j = 1, 2$. By Lemma \ref{L2.2}, the 5-polytopes
$u_{64}\ast D_3, \dots, u_{64}\ast D_8$ give a polytopal
subdivision $\widetilde{D^5}$ of $D^5$ and $L_1 := \{u_{64}\ast
\alpha \, : \, \alpha\in X^4[D_3] \cup \cdots \cup X^4[D_8]\}$ is
a simplicial subdivision of $D^5$. Now, from $X^4$, we have\,:
\begin{eqnarray*}
u_{64}\ast X^4[D_3] & = & u_{64} \ast ((1,2)\times (2, 4,3))(K), \\
u_{64}\ast X^4[D_4] & = & u_{64} \ast ((2, 1,3)\times (1,2, 3))(K), \\
u_{64}\ast X^4[D_5] & = &
u_{64} \ast ((1,3)\times (1,3)(2, 4))(\widehat{C}^4_8), \\
u_{64}\ast X^4[D_6] & = &
u_{64} \ast ((1,2,3)\times (1,4)(2,3))(\widehat{C}^4_8), \\
u_{64}\ast X^4[D_7] & = & \{u_{64}x_{133}x_{144}x_{113}x_{142}x_{143}\}
\, \cup u_{64}x_{133}x_{144}x_{112}\ast
\begin{picture}(7,1.5)(-1,1)
\setlength{\unitlength}{2.5mm}
\put(1.7,0){$_{\bullet}$}
\put(1.7,3.85){$_{\bullet}$}
\put(4.8,0){$_{\bullet}$}
\put(4.8,3.85){$_{\bullet}$}
\put(-0.1,1.9){$_{\bullet}$}
\thicklines
\put(2,0){\line(1,0){3}}
\put(2,0){\line(-1,1){2}}
\put(5,0){\line(0,1){4}}
\put(2,4){\line(1,0){3}}
\put(0,2){\line(1,1){2}}
\put(1,1.8){\small {$x_{113}$}}
\put(-1.5,0){\small {$x_{142}$}}
\put(6,0){\small {$x_{132}$}}
\put(5.8,3.9){\small {$x_{134}$}}
\put(-1.5,3.9){\small {$x_{114}$}}
\end{picture}, \\
u_{64}\ast X^4[D_8] & = &
\{u_{64}x_{233}x_{244}x_{213}x_{242}x_{243}\} \, \cup
\begin{picture}(6,2.7)(-1,1)
\setlength{\unitlength}{2.5mm}
\put(1.7,0){$_{\bullet}$}
\put(1.7,3.85){$_{\bullet}$}
\put(4.8,0){$_{\bullet}$}
\put(4.8,3.85){$_{\bullet}$}
\put(6.7,1.9){$_{\bullet}$}
\thicklines
\put(2,0){\line(1,0){3}}
\put(2,0){\line(0,1){4}}
\put(5,0){\line(1,1){2}}
\put(5,4){\line(1,-1){2}}
\put(2,4){\line(1,0){3}}
\put(-1.5,0){\small {$x_{242}$}}
\put(6.7,0){\small {$x_{232}$}}
\put(3.5,1.8){\small {$x_{234}$}}
\put(6.2,3.9){\small {$x_{214}$}}
\put(-1.5,3.9){\small {$x_{213}$}}
\end{picture} \ast
u_{64}x_{233}x_{244}x_{212},
\end{eqnarray*}

\medskip

\noindent where $K$ is as in Lemma \ref{R5} and $\widehat{C}^4_8$
is as in equation (\ref{C48}). These give all the $24+24+16+16+6+6
= 92$ facets of $L_1$. Observe that $L_1[D_1] \neq X^4[D_1]$ and
$L_1[D_2] \neq X^4[D_2]$. Also, the eight edges $u_{64}u_{16}$,
$u_{64}u_{51}$, $u_{64}x_{143}$, $u_{64}x_{243}$, $u_{64}x_{132}$,
$u_{64}x_{214}$, $u_{64}x_{142}$ and $u_{64}x_{213}$ are in $L_1$
and not in $X^4[D^5]$. Forty of these ninetytwo facets do not
contain any of the vertices from $\{x_{143}, x_{243}, x_{142},
x_{213}, x_{132}, x_{214}\}$. These are the first forty in the
list given in the lemma.

Observe that ${\rm lk}_{L_1}(u_{64}x_{143}x_{243}) = S^{\,2}_4(\{
x_{113}, x_{144}, x_{233}, x_{242}\})$ and $\frac{1}{2}u_{64} +
\frac{1}{4}x_{143} + \frac{1}{4} x_{243} = (0, 0, -1, \frac{1}{2},
0, \frac{1}{2}, -\frac{1}{2}, 0, \frac{1}{2})= \frac{1}{4}(x_{113}
+ x_{144} + x_{233} + x_{242})$. Thus, by Lemma \ref{L2.3}, $Y_1
:= \overline{u_{64} x_{143} x_{243}} \ast S^{\,2}_4(\{x_{113},
x_{144}, x_{233}, x_{242}\})$ and $Z_1 := S^1_3(\{u_{64}, x_{143},
x_{243}\}) \ast \overline{x_{113}x_{144}x_{233}x_{242}}$
triangulate the 5-polytope $\langle\{u_{64}, x_{143}, x_{243},
x_{113}, x_{144}, x_{233}, x_{242}\}\rangle$. Let $L_2$ be the
simplicial complex obtained from $L_1$ by replacing $Y_1$ by
$Z_1$. Then $L_2$ triangulates $D^5$ and has $92-1 = 91$ facets.

Now, ${\rm lk}_{L_2}(u_{64}x_{143}) = S^1_3(\{u_{51}, x_{113},
x_{144}\}) \ast S^{\,0}_2(x_{133}, x_{242}) \ast
S^{\,0}_2(x_{142}, x_{233})$ and $\frac{1}{2}(u_{64} + x_{143}) =
\frac{1}{2}(1, -1, -2, 1, 0, 1, -1, 0, 1) = \frac{1}{2}u_{51} +
\frac{1}{4}x_{113}  + \frac{1}{4}x_{144}$. Therefore, by Lemma
\ref{L2.3}, $|\overline{\gamma}\ast \overline{u_{64}x_{143}} \ast
S^1_3(\{u_{51}, x_{113}, x_{144}\})| = |\overline{\gamma}\ast
S^{\,0}_2(u_{64}, x_{143}) \ast \overline{u_{51} x_{113}x_{144}}|$
for $\gamma = x_{133}x_{142}$, $x_{133}x_{233}$, $x_{242}x_{142}$,
$x_{242}x_{233}$. Thus, we can replace $Y_2:=
\overline{u_{64}x_{143}} \ast S^1_3(\{u_{51}, x_{113}, x_{144}\})
\ast S^{\,0}_2(x_{133}, x_{242}) \ast S^{\,0}_2(x_{142}, x_{233})$
in $L_2$ by $Z_2 := S^{\,0}_2(u_{64}, x_{143}) \ast
\overline{u_{51}x_{113} x_{144}} \ast S^{\,0}_2(x_{133}, x_{242})
\ast S^{\,0}_2(x_{142}, x_{233})$. Similarly, we can replace
$f(Y_2) = \overline{u_{64}x_{243}} \ast S^1_3(\{u_{16}, x_{242},
x_{233}\}) \ast S^{\,0}_2(x_{244}, x_{113})\ast S^{\,0}_2(x_{213},
x_{144})$ in $L_2$ by $f(Z_2) = S^{\,0}_2(u_{64}, x_{243}) \ast
\overline{u_{16}x_{242} x_{233}} \ast  S^{\,0}_2(x_{244},
x_{113})\ast S^{\,0}_2(x_{213}, x_{144})$. Let $L_3$ be obtained
from $L_2$ by replacing $Y_2$ by $Z_2$ and by replacing $f(Y_2)$
by $f(Z_2)$. Then $u_{64}x_{143}$, $u_{64}x_{243}$ are non-edges
in $L_3$, $L_3$ triangulates $D^5$ and $L_3$ has $91- (2\times 4)
= 83$ facets.

Consider the subcomplex $Y_3 := \{x_{212}x_{244}u_{64}x_{214}
x_{213}x_{114}\} \cup (\overline{x_{244}x_{233}u_{64}x_{214}} \ast
S^{\,0}_2(x_{212},$ $x_{114}) \ast S^{\,0}_2(x_{234}, x_{213}))$
of $L_3$. Then $Y_3 = {\rm star}_{L_3}(u_{64}x_{214})$ and hence
$Y_3$ is a combinatorial 5-ball (i.e., ($|Y_3|$ is pl homeomorphic
to a 5-simplex). Now, $\alpha_1 := x_{233}x_{244}x_{212}u_{64}
x_{234}$ is a 4-simplex (in $Y_3$). If there exist real numbers
$a_1, \dots, a_5$ such that $a_1x_{244} + a_2x_{233} + a_3u_{64} +
a_4x_{212} + a_5x_{234} = x_{114}$, then from first, third, eighth
and ninth coordinates $a_2=0$, $a_3 = 2$ and $a_4=-1$. Now, from
fourth and fifth coordinates, we get $2= a_1-a_5 = 0$, a
contradiction. Thus, $x_{114} \not\in {\rm AH}(\alpha_1)$ and
hence $x_{114} \ast \alpha_1$ is a 5-simplex (inside $D^5$). Again
$\alpha_2 := x_{212}x_{233}x_{244}x_{214}x_{213}$, $\alpha_3 :=
x_{212}x_{233}x_{244}x_{214}x_{234}$ are simplices (in $Y_3$),
$|\alpha_2|$, $|\alpha_3| \subseteq D_8$ and $x_{114} \not \in
D_8$. These imply that $x_{114} \ast \alpha_2$ and $x_{114} \ast
\alpha_3$ are 5-dimensional simplices (with geometric carriers in
$D^5$). So, $Z_3 := \overline{x_{114}x_{212}x_{233}x_{244}} \ast
\begin{picture}(10,1.1)(0,0) \setlength{\unitlength}{2mm}
\put(0.3,0.2){$_{\bullet}$} \put(5.8,0.2){$_{\bullet}$}
\put(11.3,0.2){$_{\bullet}$} \put(16.8,0.2){$_{\bullet}$}
\thicklines \put(0.5,0.3){\line(1,0){16.5}} \put(0.2,0.9){\small
{$u_{64}$}} \put(5.7,0.9){\small {$x_{234}$}}
\put(11.2,0.9){\small {$x_{214}$}} \put(16,0.9){\small
{$x_{213}$}}
\end{picture}
$ is a simplicial complex and $|Z_5| \subseteq D^5$. Clearly,
$Z_3$ is a combinatorial 5-ball.

\medskip

\noindent {\sf Claim 1.} $|Y_3| = |Z_3|$.

\smallskip

Observe that $\partial Y_3 = \{x_{214}x_{244}x_{213}x_{212}
x_{114}\} \cup (\overline{x_{214}x_{233}x_{244}} \ast
S^{\,0}_2(x_{212}, x_{114}) \ast S^{\,0}_2(x_{234}, x_{213}))$
$\cup (\overline{u_{64}x_{233}x_{244}x_{234}} \ast
S^{\,0}_2(x_{212}, x_{114})) \cup M_1 \cup (\overline{u_{64}
x_{214}x_{233}x_{234}} \ast S^{\,0}_2(x_{212}, x_{114})) \cup
(\overline{u_{64}x_{214}x_{244}} \ast S^{1}_3(\{x_{212}, x_{234},
x_{114}\}))$, where $M_1:= S^{\,0}_2(x_{214}, x_{244}) \ast
\overline{u_{64}x_{213}} \ast  S^{1}_3(\{x_{212}, x_{233},
x_{114}\})$. Since $M_1 \subseteq \partial Y_3$ and
$\frac{1}{3}(-1, 1, -3, 1, -3, -1, -1, 1,1)$ is an interior point
of both the simplices $x_{114}x_{212}x_{233}$ and $u_{64}x_{213}$,
by Lemma \ref{L2.3} (iii), $|\overline{\gamma} \ast
\overline{u_{64}x_{213}} \ast S^{1}_3(\{x_{212}, x_{233},
x_{114}\})| = |\overline{\gamma} \ast S^{\,0}_2(u_{64}, x_{213})
\ast \overline{x_{212}x_{233} x_{114}}|$, for $\gamma = x_{214},
x_{244}$. Thus $|M_1| = |N_1|$, where $N_1 := S^{\,0}_2(x_{214},
x_{244}) \ast S^{\,0}_2(u_{64}, x_{213}) \ast
\overline{x_{212}x_{233}x_{114}}$. Let $\widetilde{\partial Y_3}$
be the complex obtained from $\partial Y_3$ by replacing $M_1$ by
$N_1$. Then $|\widetilde{\partial Y_3}| = |\partial Y_3|$.

Similarly, $|M_2| = |N_2|$, where $M_2 :=  S^{\,0}_2(x_{233},
x_{244}) \ast \overline{u_{64}x_{214}} \ast S^{1}_3(\{x_{212},
x_{234}, x_{114}\})$ is a subcomplex of $\widetilde{\partial Y_3}$
and $N_2 := S^{\,0}_2(x_{233}, x_{244}) \ast S^{\,0}_2(u_{64},
x_{214}) \ast \overline{x_{212}x_{234}x_{114}}$. Let
$\widehat{\partial Y_3}$ be the complex obtained from
$\widetilde{\partial Y_3}$ by replacing $M_2$ by $N_2$. Then
$\widehat{\partial Y_3} = \partial Z_3$. So, $|\partial Y_3| =
|\widetilde{\partial Y_3}| = |\widehat{\partial Y_3}| = |\partial
Z_3|$. Now, $|Y_3|$ and $|Z_3|$ are two 5-balls in the
5-dimensional affine space ${\rm AH}(D^5)$ and $\partial |Y_3| =
|\partial Y_3| = |\partial Z_3| = \partial |Z_3|$. These imply
that $|Y_3| = |Z_3|$. This proves the claim.

\smallskip

Similarly,  $|f(Y_3)| = |f(Z_3)|$, where $f(Y_3) =
\{x_{112}x_{133}u_{64}x_{132}x_{142}x_{232}\} \cup
(\overline{x_{133}x_{144}u_{64}x_{132}} \ast S^{\,0}_2(x_{112},
x_{232}) \ast S^{\,0}_2(x_{134},$ $x_{142}))$ and $f(Z_3) :=
\overline{x_{112}x_{133}x_{144}x_{232}} \ast
\begin{picture}(10,1.1)(0,0)
\setlength{\unitlength}{2mm}
\put(0.3,0.2){$_{\bullet}$}
\put(5.8,0.2){$_{\bullet}$}
\put(11.3,0.2){$_{\bullet}$}
\put(16.8,0.2){$_{\bullet}$}
\thicklines
\put(0.5,0.3){\line(1,0){16.5}}
\put(0.2,0.9){\small {$u_{64}$}}
\put(5.7,0.9){\small {$x_{134}$}}
\put(11.2,0.9){\small {$x_{132}$}}
\put(16,0.9){\small {$x_{142}$}}
\end{picture}
$. Let $L_4$ be the simplicial complex obtained from $L_3$ by
replacing $Y_3$ by $Z_3$ and $f(Y_3)$ by $f(Z_3)$. Then, by Claim
1, $|L_4| = |L_3|$, and hence $L_4$ triangulates $D^5$. Observe
that $u_{64}x_{132}$, $u_{64}x_{214}$ are non-edges in $L_4$ and
$L_4$ has $83 -2\times 2 = 79$ facets.


Consider the subcomplex $Y_4 := \overline{u_{64}x_{142}}\ast
S^{\,0}_2(x_{133}, x_{242}) \ast D$ of $L_4$, where $D = \{x_{112}
x_{113}x_{144},$ $u_{51}x_{112}x_{113}, u_{51}x_{112} x_{232},
u_{51}x_{113}x_{144}, u_{51}x_{144}x_{232}\}$. Then $D$ is a
triangulation of the 2-disc and hence $Y_4$ is a combinatorial
5-ball. Let $M:=S^{\,0}_2(x_{133}, x_{242})\ast S^{\,0}_2(x_{113},
x_{232}) \ast \overline{u_{51}x_{112}} \subseteq Y_4$. Then
$\overline{u_{64}} \ast M$, $\overline{x_{142}} \ast M \subseteq
Y_4$. Since $|\overline{x_{142}} \ast M| \subseteq D_4$ and
$x_{144} \not \in D_4$, $\overline{x_{142}} \ast M \ast
\overline{x_{144}}$ is a simplicial complex. Again, from the 7-th
coordinate of $x_{144}$, it is clear that $x_{144} \not\in {\rm
AH}(|\overline{u_{64}} \ast M|)$. Thus, $\overline{u_{64}} \ast M
\ast \overline{x_{144}}$ is a simplicial complex. These imply $Z_4
:= S^{\,0}_2(u_{64}, x_{142})\ast S^{\,0}_2(x_{133}, x_{242}) \ast
S^{\,0}_2(x_{113}, x_{232}) \ast \overline{u_{51}x_{112}x_{144}}$
is a simplicial complex. Clearly, $Z_4$ is a combinatorial 5-ball.

\medskip

\noindent {\sf Claim 2.} $|Y_4| = |Z_4|$.

\smallskip

Observe that $\partial Y_4 := (S^{\,0}_2(u_{64}, x_{142}) \ast
S^{\,0}_2(x_{133}, x_{242}) \ast D) \cup N$, where $N =
\overline{u_{64}x_{142}} \ast S^{\,0}_2(x_{133}, x_{242}) \ast
S^{1}_3(\{x_{112}, x_{144}, x_{232}\})$. Since $\frac{1}{3}(1, -1,
-3, 1, -1, 1, -1, 3, -1)$ is an interior point of both
$u_{64}x_{142}$ and $x_{112}x_{144}x_{232}$, by Lemma \ref{L2.3},
$|\overline{\gamma} \ast \overline{u_{64}x_{142}} \ast
S^{1}_3(\{x_{112}, x_{144}, x_{232}\})|$ $= |\overline{\gamma}
\ast S^{\,0}_2(u_{64}, x_{142}) \ast \overline{x_{112}x_{144}
x_{232}}|$, for $\gamma = x_{133}$, $x_{242}$. This implies $|N| =
|\widetilde{N}|$, where $\widetilde{N} = S^{\,0}_2(u_{64},
x_{142}) \ast S^{\,0}_2(u_{133}, x_{242}) \ast \overline{x_{112}
x_{144}x_{232}}$. Observe that $\partial Z_4$ can be obtained from
$\partial Y_4$ by replacing $N$ by $\widetilde{N}$. Thus
$|\partial Z_4| = |\partial Y_4|$. Now, $|Y_4|$ and $|Z_4|$ are
two 5-balls in the 5-dimensional affine space ${\rm AH}(D^5)$ and
$\partial |Y_4| = |\partial Y_4| = |\partial Z_4| = \partial
|Z_4|$. Therefore, $|Y_4| = |Z_4|$. This proves the claim.

\smallskip

Let $L_5$ be the simplicial complex obtained from $L_4$ by
replacing $Y_4$ by $Z_4$ and $f(Y_4)$ by $f(Z_4)$. Then, by Claim
2,  $|L_5| = |L_4|$, and hence $L_5$ triangulates $D^5$. Clearly,
$u_{64}x_{142}$, $u_{64}x_{213}$ are non-edges in $L_5$ and $L_5$
has $79 -2\times 2 = 75$ facets. Observe that $L_5 = L$ (the
sixteen facets in $Z_4 \cup f(Z_4)$ are the last sixteen
facets in $L$). This completes the proof. \hfill $\Box$


\begin{lemma}$\!\!\!${\bf .} \label{R8}
Let $X^4$, $W^j$ $(0\leq j \leq 6)$ and $V_{124}$ be as in  Lemma
$\ref{R6}$. There exists a simplicial complex $X^5$ which satisfy
the following\,: {\rm (i)} $X^5$ is a simplicial subdivision of
$W^5$, {\rm (ii)} the $4$-skeleton of $X^5$ is $X^4$, {\rm (iii)}
$S_3 \times A_4$ acts as an isometry group, where the action of
$S_3$ and $A_4$ on the vertices is as in Lemma $\ref{R4}$, {\rm (iv)}
the induced $S_3$ $($respectively, $A_4)$ action on $|X^5|$ is the
same as that on $|W^5|$ and {\rm (v)} the $S_3$ action on $X^5$ is
good.
\end{lemma}

\noindent {\bf Proof.} A 5-polytope is the product of two
2-simplices  and one 1-simplex. Up to the action of $S_3\times
A_4$, the 5-polytopes are of the form $C^5_1 = x_1x_3 \times
x_2x_3x_4 \times x_2x_3x_4$, $C^5_2 = x_3x_4 \times x_2x_3x_4
\times x_2x_3x_4$, $C^5_3 = x_3x_4 \times x_1x_3x_4 \times
x_2x_3x_4$, $C^5_4 = x_1x_3 \times x_1x_3x_4 \times x_2x_3x_4$,
$C^5_5 = x_1x_4 \times x_1x_3x_4 \times x_2x_3x_4$ or $C^5_6 =
x_1x_2 \times x_1x_3x_4 \times x_2x_3x_4$. For $1\leq j\leq 6$, we
know (from Lemma \ref{R6}) the simplicial subdivision
$X^4[\partial \overline{C^5_j}]$ of the boundary complex $\partial
\overline{C^5_j}$ of $C^5_j$. For $1 \leq j \leq 4$, let $S_j$ be
the anti-star of the vertex $x_{333}$ in $X^4[\partial
\overline{C^5_j}]$ and let $S_5$ be the anti-star of the vertex
$x_{444}$ in $X^4[\partial \overline{C^5_5}]$. Then
\begin{eqnarray*}
S_1 &=&
X^4[x_1 \times x_2x_3x_4 \times x_2x_3x_4] \cup X^4[x_1x_3 \times
x_2x_4 \times x_2x_3x_4] \cup X^4[x_1x_3 \times x_2x_3x_4 \times
x_2x_4] \\
&=& ({\rm Id}\times (14)(23))(\widehat{C^4_{1,4}}) \cup ((123) \times (124))
(\widehat{C^4_8}) \cup ((13) \times (124))(\widehat{C^4_8}), \\
S_2 &=& X^4[x_4 \times x_2x_3x_4 \times x_2x_3x_4] \cup X^4[x_3x_4
\times x_2x_4 \times x_2x_3x_4] \cup X^4[x_3x_4 \times x_2x_3x_4
\times x_2x_4], \\
&=& ({\rm Id}\times (14)(23))(\widehat{C^4_{1,1}}) \cup ({\rm Id}
\times (14)(23))(\widehat{C^4_5}) \cup ((23) \times (14)(23))
(\widehat{C^4_5}), \\
S_3 &=& X^4[x_4 \times x_1x_3x_4 \times x_2x_3x_4] \cup X^4[x_3x_4
\times x_1x_4 \times x_2x_3x_4] \cup X^4[x_3x_4 \times x_1x_3x_4
\times x_2x_4], \\
&=& ({\rm Id}\times (13)(24))(\widehat{C^4_{2,2}}) \cup ((12) \times (142))
(\widehat{C^4_6}) \cup ((132) \times (143))(\widehat{C^4_7}), \\
S_4 &=& X^4[x_1 \times x_1x_3x_4 \times x_2x_3x_4] \cup X^4[x_1x_3
\times x_1x_4 \times x_2x_3x_4] \cup X^4[x_1x_3 \times x_1x_3x_4
\times x_2x_4] \\
&=& ({\rm Id}\times (13)(24))(\widehat{C^4_{2,3}}) \cup ({\rm Id} \times (234))
(\widehat{C^4_9}) \cup ({\rm Id} \times (234))(\widehat{C^4_8}), \\
S_5 &=& X^4[x_1 \times
x_1x_3x_4 \times x_2x_3x_4] \cup X^4[x_1x_4 \times x_1x_3 \times
x_2x_3x_4] \cup X^4[x_1x_4 \times x_1x_3x_4 \times x_2x_3] \\
&=& ({\rm Id}\times (13)(24))(\widehat{C^4_{2,3}}) \cup ((12)
\times (234))
(\widehat{C^4_9}) \cup ({\rm Id} \times (142))(\widehat{C^4_8}),
\end{eqnarray*}
where $\widehat{C^4_{j,i}}$, $\widehat{C^4_k}$ are as in the
proof of Lemma \ref{R6}. We take $\widehat{C^5_j} :=x_{333}
\ast S_j$ as the triangulation of $C^5_j$ for $1\leq j \leq 4$
and $\widehat{C^5_5} := x_{444} \ast S_5$ as the triangulation of
$C^5_5$. We take the simplicial complex $L$ in Lemma \ref{R7} as
the simplicial subdivision $\widehat{C^5_6}$ of $C^5_6$.
By the action of $S_3 \times A_4$, we get simplicial subdivisions
of all the 5-polytopes in $W^5$ and get $X^5$.
Since $u_{42} \in S_5$, the edge-set of $X^5$ is
the following\,:
\begin{eqnarray}
E(X^5) \!\!\! & =  \!\!\!& E(X^4) \, \cup \, (S_3\times
A_4)(\{x_{333}x_{122}, x_{444}u_{42}, u_{64}u_{16}\}) \nonumber \\
&=  \!\!\!& E(X^4) \, \cup \, (S_3\times A_4)(\{x_{111}x_{233},
x_{111}u_{11},u_{11}u_{22}\})
\subseteq E(X), \label{EX5}
\end{eqnarray}
where $E(X^4)$ and $E(X)$ are as in equations (\ref{EX4}) and
(\ref{edges}) respectively. Since $E(X^5) \subseteq E(X)$, it
follows that the $S_3$-action on $X^5$ is good. \hfill $\Box$

\bigskip

\noindent {\bf Proof of Theorem \ref{T1}.} Let $W$, $X$, $W^k$ and
$X^k$ be as in the sketch of the proof. A facet in the polytopal
complex $S^{\,2}_4 \times  S^{\,2}_4 \times S^{\,2}_4$ is the
product of three 2-simplices. Up to the action of $S_3\times A_4$,
the facets are of the form $C^6_1 = x_1x_2x_4 \times x_1x_2x_3
\times x_1x_3x_4$, $C^6_2 = x_1x_2x_3 \times x_1x_2x_3 \times
x_1x_3x_4$, $C^6_3 = x_1x_2x_3 \times x_1x_2x_3 \times x_1x_2x_3$.
So, the set of facets in $S^{\,2}_4 \times  S^{\,2}_4 \times
S^{\,2}_4$ is $(S_3\times A_4)(\{C^6_1, C^6_2, C^6_3\})$.

By Lemma \ref{L2.2}, $\{x_{111}\ast (x_2x_4 \times x_1x_2x_3
\times x_1x_3x_4), x_{111}\ast (x_1x_2x_4 \times x_2x_3 \times
x_1x_3x_4), x_{111}\ast (x_1x_2x_4 \times x_1x_2x_3 \times
x_3x_4)\}= \{x_{111}\ast ((23)\times (143))(D^5), x_{111}\ast
((123)\times (134))(D^5), x_{111}\ast ((132)\times (14)(23))(D^5)$
gives a polytopal subdivision of $C^6_1$, where $D^5 = x_1x_2
\times x_1x_3x_4 \times x_2x_3x_4$. We choose $L$ in Lemma
\ref{R7} as the subdivision (as in Lemma \ref{R8}) of $D^5$.
Accordingly, $\overline{x_{111}}\ast ((23)\times (143)(L)) \cup
\overline{x_{111}}\ast ((123)\times (134)(L)) \cup
\overline{x_{111}}\ast ((132)\times (14)(23)(L))$ triangulates
$D^6_1$. Thus, $(S_3\times A_4)(\overline{x_{111}} \ast L)$ gives
simplicial subdivisions of all the seventytwo facets $(S_3 \times
A_4)(D^6_1)$. These give $72 \times 75$ simplices which are in the
first seventyfive orbits in the description of $(S^{\,2} \times
S^{\,2}\times S^{\,2})_{124}$. Here we are adding new set $(S_3
\times A_4)(\{x_{111}v_{342}\})$ of edges.

By Lemma \ref{L2.2}, $\{x_{111}\ast (x_2x_3 \times x_1x_2x_3
\times x_1x_3x_4), x_{111}\ast (x_1x_2x_3 \times x_2x_3 \times
x_1x_3x_4), x_{111}\ast (x_1x_2x_3 \times x_1x_2x_3 \times
x_3x_4)\}$ gives a polytopal subdivision of $C^6_2$. Now,
$\{x_{333}\ast (x_2 \times x_1x_2x_3 \times x_1x_3x_4),
x_{333}\ast (x_2x_3 \times x_1x_2 \times x_1x_3x_4), x_{333}\ast
(x_2x_3 \times x_1x_2x_3 \times x_1x_4)= \{({\rm Id} \times
(132))(C^4_{2,3}), ({\rm Id} \times (123))(C^4_{9}), ({\rm Id}
\times (123))(C^4_{8})\}$ is a polytopal subdivision of $x_2x_3
\times x_1x_2x_3 \times x_1x_3x_4= ({\rm Id} \times
(124))(C^5_4)$, where $C^4_{2,3}, C^4_{8}, C^4_{9}$ are as in the
proof of Lemma \ref{R6}. Again, Lemma \ref{L2.2}, $\{x_{333}\ast
(x_1x_2 \times x_1x_2x_3 \times x_3x_4), x_{333}\ast (x_1x_2x_3
\times x_1x_2\times x_3x_4), x_{333}\ast (x_1x_2x_3 \times
x_1x_2x_3 \times x_4)\}$ is a polytopal subdivision of $x_1x_2x_3
\times x_1x_2x_3 \times x_3x_4= (13 \times (142))(C^5_1)$ and
$x_1x_2 \times x_1x_2x_3 \times x_3x_4 = C^4_8$, $x_{123} \times
x_{123} \times x_4 = ((123) \times {\rm Id})(C^4_{1,4})$, where
$C^4_{1,4}, C^4_{8}, C^5_1$ are as in the proofs of Lemmas
\ref{R6} and \ref{R8}. Also, $((123) \times {\rm
Id})(\widehat{C}^4_{1,4}) = \{x_{114}x_{224}x_{334}x_{124}x_{134},
x_{114}x_{224}x_{334}x_{134}x_{234}, x_{114}x_{224}x_{334}x_{234}
x_{214}\} \cup ((12)$ $\times {\rm Id})(\{x_{114}x_{224}x_{334}
x_{124}x_{134}, x_{114}x_{224}x_{334}x_{134}x_{234},
x_{114}x_{224}x_{334}x_{234}x_{214}\})$. These imply that
$\overline{x_{111}x_{333}} \ast [A \cup ((12) \times {\rm
Id})(A)]$ is a simplicial subdivision of the 6-polytope $D^6_2$,
where \linebreak $A = ({\rm Id} \times (132))(\widehat{C}^4_{2,3})
\, \cup \, ({\rm Id} \times (123))(L)\, \cup \, ({\rm Id} \times
(123))(\widehat{C}^4_{8}) \, \cup \, \widehat{C}^4_8 \, \cup \,
\{x_{114}x_{224}x_{334}x_{124}x_{134}$,
$x_{114}x_{224}x_{334}x_{134}x_{234}$,
$x_{114}x_{224}x_{334}x_{234}x_{214}\}$. Thus, $(S_3\times
A_4)(\overline{x_{111}x_{333}} \ast A)$ gives simplicial
subdivisions of all the thirtysix facets $(S_3 \times
A_4)(D^6_2)$. These give $72 \times 65$ simplices which are in
76-th to 140-th orbits in the description of $(S^{\,2} \times
S^{\,2}\times S^{\,2})_{124}$.

By Lemma \ref{L2.2}, $\{x_{111}\ast (x_2x_3 \times x_1x_2x_3
\times x_1x_2x_3), x_{111}\ast (x_1x_2x_3 \times x_2x_3 \times
x_1x_2x_3), x_{111}\ast (x_1x_2x_3 \times x_1x_2x_3 \times
x_2x_3)\}$ gives a polytopal subdivision of $C^6_3$. Again, by
Lemma \ref{L2.2}, $\{x_{222}\ast (x_3 \times x_1x_2x_3 \times
x_1x_2x_3), x_{222}\ast (x_2x_3 \times x_1x_3 \times x_1x_2x_3),
x_{222}\ast (x_2x_3 \times x_1x_2x_3 \times x_1x_3)\}$ gives a
polytopal subdivision of $x_2x_3 \times x_1x_2x_3 \times
x_1x_2x_3$, $\{x_{333}\ast (x_3 \times x_1x_2 \times x_1x_2x_3),
x_{333}\ast (x_3 \times x_1x_2x_3 \times x_1x_2)\}$ gives a
polytopal subdivision of $x_3 \times x_1x_2x_3 \times x_1x_2x_3$
and $\{x_{333}\ast (x_2 \times x_1x_3 \times x_1x_2x_3),
x_{333}\ast (x_2x_3 \times x_1 \times x_1x_2x_3), x_{333}\ast
(x_2x_3 \times x_1x_3 \times x_1x_2)\}$ gives a polytopal
subdivision of $x_2x_3 \times x_1x_3 \times x_1x_2x_3$. So, $(S_3
\times A_3)(\{x_{111}x_{222}x_{333} \ast (x_1 \times x_{1}x_2x_3
\times x_2x_3), x_{111}x_{222}x_{333} \ast (x_1x_2 \times x_1x_3
\times x_2x_3)\})$ is a polytopal subdivision of $C^6_3$. We take
(see Lemma \ref{R4}, Fig. 3) $\{x_{122}x_{112}x_{132}x_{133}$,
$x_{122}x_{112}x_{113}x_{133}$, $x_{122}x_{123}x_{113}x_{133}\}$
is the subdivision of $(x_1 \times x_{1}x_2x_3 \times x_2x_3$ and
$\widehat{C^3_5} = H(\{x_{132}x_{112}x_{133}x_{232}$,
$x_{213}x_{233}x_{212}x_{113}$, $x_{112}x_{133}x_{113}v_{321}$,
$x_{112}x_{212}x_{113}v_{321}$, $x_{112}x_{133}x_{232}v_{321}$,
$x_{233}x_{212}x_{113}v_{321}\})$ in Lemma \ref{R4} is the
subdivision of $x_1x_2 \times x_1x_3 \times x_2x_3$, where $H
\langle (1,2,3)\times (1,2,3)\rangle$ is a subgroup of $S_3 \times
A_3$. Thus, $(S_3\times A_4) (\overline{x_{111}x_{222}x_{333}}
\ast \{x_{122}x_{133}x_{112}x_{132}$,
$x_{122}x_{133}x_{112}x_{113}$, $x_{122}x_{133}x_{113}x_{123}$,
$x_{112}x_{133}x_{113}v_{321}$, $x_{112}x_{212}x_{113}v_{321}$,
$x_{132}x_{112}x_{133}x_{232}$, $x_{213}x_{233}x_{212}x_{113}$,
$x_{112}x_{133}x_{232}v_{321}$, $x_{233}x_{212}x_{113}v_{321}\})$
\, \, gives simplicial subdivisions of all the 4 faces $(S_3\times
A_4)(D^6_3)$. These give $5\times 72 + 4 \times 24$ simplices
which are the last nine orbits in the description of $(S^{\,2} \times
S^{\,2}\times S^{\,2})_{124}$. Observe that the set of edges of
$X$ is $E(X^5) \cup (S_3\times A_4)(\{x_{111}v_{234}\})$ which is
same as $E(X)$ in equation (\ref{edges}), where $E(X^5)$ is as in
equation (\ref{EX5}).

We have shown in the sketch of the proof of Theorem \ref{T1} that
the $S_3$-action on $X$ is good. This implies that the
$S_3$-action on the abstract scheme $(S^{\,2}\times S^{\,2} \times
S^{\,2})_{124}$ of $X$ is good. Now, Lemma \ref{L2.4} and
Corollary \ref{L2.6} imply that the quotient $\CC P^{\,3}_{30} :=
(S^{\,2}\times S^{\,2} \times S^{\,2})_{124}/S_3$ triangulates
$\CC P^{\,3}$. This completes the proof. \hfill $\Box$

\section{Appendix: An 18-vertex triangulation of \boldmath{$\CCC P^3$}
\newline \mbox{} \hspace{20mm} {\rm (by Bhaskar Bagchi, Basudeb Datta,
Nitin Singh)}}

We applied the BISTELLAR program of Lutz (\cite{lu2}) on $\CC
P^{\,3}_{30}$. The BISTELLAR flips were carried out on a PC with
intel Atom processor (1.66 Ghz) for approximately 16 hours. The
final output was  an 18-vertex simplicial complex $\CC
P^{\,3}_{18}$ which is bistellar equivalent to $\CC P^{\,3}_{30}$
and hence triangulates $\CC P^{\,3}$. Its face vector is $(18,
153, 783, 2110, 3021, 2177, 622)$. Its automorphism group is
trivial (verified by the simpcomp program of Effenberger and
Spreer \cite{es}). After an appropriate renaming of the vertices,
the description of $\CC P^{\,3}_{18}$ is as follows. The vertices
are $a_i$, $b_i$, $1\leq i \leq 9$. The induced subcomplex of $\CC
P^{\,3}_{18}$ on the nine vertices $a_i$, $1\leq i \leq 9$, is a
six-dimensional ball (whose facets are the first seven facets in
the following list). The induced subcomplex of $\CC P^{\,3}_{18}$
on the remaining nine vertices $b_i$, $1\leq i \leq 9$, is an
isomorphic copy of K\"{u}hnel's $\CC P^{\,2}_9$. The map $b_i
\mapsto i$ gives an isomorphism of this subcomplex with $\CC
P^{\,2}_9$ as given in \cite{kb2}. This shows that the geometric
carrier of $\CC P^{\,3}_{18}$ is a compactification by $\CC
P^{\,2}$ of an open 6-ball, hence is homeomorphic to $\CC
P^{\,3}$. A computer search shows that $\CC P^{\,3}_{18}$ has a
unique induced 9-vertex 6-ball, namely the one mentioned above.
The complete list of the 622 facets of $\CC P^{\,3}_{18}$ is as
follows\,:

\medskip

$a_1a_2a_3a_4a_5a_7a_8 $, $a_1a_2a_3a_4a_5a_7a_9 $,
$a_1a_2a_3a_4a_5a_8a_9 $, $a_1a_2a_3a_4a_6a_7a_8 $,
$a_1a_2a_3a_4a_6a_7a_9 $,

$a_1a_2a_3a_4a_6a_8a_9 $, $a_1a_2a_3a_5a_6a_7a_9 $,
$a_1a_2a_3a_5a_6a_7b_2 $, $a_1a_2a_3a_5a_6a_9b_1 $,
$a_1a_2a_3a_5a_6b_1b_2 $,

$a_1a_2a_3a_5a_7a_8b_2 $, $a_1a_2a_3a_5a_8a_9b_5 $,
$a_1a_2a_3a_5a_8b_2b_5 $, $a_1a_2a_3a_5a_9b_1b_5 $,
$a_1a_2a_3a_5b_1b_2b_5 $,

$a_1a_2a_3a_6a_7a_8b_2 $, $a_1a_2a_3a_6a_8a_9b_5 $,
$a_1a_2a_3a_6a_8b_2b_5 $, $a_1a_2a_3a_6a_9b_1b_5 $,
$a_1a_2a_3a_6b_1b_2b_5 $,

$a_1a_2a_4a_5a_7a_8b_3 $, $a_1a_2a_4a_5a_7a_9b_1 $,
$a_1a_2a_4a_5a_7b_1b_3 $, $a_1a_2a_4a_5a_8a_9b_5 $,
$a_1a_2a_4a_5a_8b_3b_8 $,

$a_1a_2a_4a_5a_8b_5b_8 $, $a_1a_2a_4a_5a_9b_1b_5 $,
$a_1a_2a_4a_5b_1b_3b_8 $, $a_1a_2a_4a_5b_1b_5b_8 $,
$a_1a_2a_4a_6a_7a_8b_6 $,

$a_1a_2a_4a_6a_7a_9b_5 $, $a_1a_2a_4a_6a_7b_5b_7 $,
$a_1a_2a_4a_6a_7b_6b_7 $, $a_1a_2a_4a_6a_8a_9b_5 $,
$a_1a_2a_4a_6a_8b_5b_7 $,

$a_1a_2a_4a_6a_8b_6b_7 $, $a_1a_2a_4a_7a_8b_3b_6 $,
$a_1a_2a_4a_7a_9b_1b_5 $, $a_1a_2a_4a_7b_1b_3b_6 $,
$a_1a_2a_4a_7b_1b_5b_6 $,

$a_1a_2a_4a_7b_5b_6b_7 $, $a_1a_2a_4a_8b_3b_6b_8 $,
$a_1a_2a_4a_8b_5b_7b_8 $, $a_1a_2a_4a_8b_6b_7b_8 $,
$a_1a_2a_4b_1b_3b_6b_8 $,

$a_1a_2a_4b_1b_5b_6b_8 $, $a_1a_2a_4b_5b_6b_7b_8 $,
$a_1a_2a_5a_6a_7a_9b_1 $, $a_1a_2a_5a_6a_7b_1b_3 $,
$a_1a_2a_5a_6a_7b_2b_3 $,

$a_1a_2a_5a_6b_1b_2b_3 $, $a_1a_2a_5a_7a_8b_2b_3 $,
$a_1a_2a_5a_8b_2b_3b_8 $, $a_1a_2a_5a_8b_2b_5b_8 $,
$a_1a_2a_5b_1b_2b_3b_8 $,

$a_1a_2a_5b_1b_2b_5b_8 $, $a_1a_2a_6a_7a_8b_2b_6 $,
$a_1a_2a_6a_7a_9b_1b_5 $, $a_1a_2a_6a_7b_1b_3b_5 $,
$a_1a_2a_6a_7b_2b_3b_7 $,

$a_1a_2a_6a_7b_2b_6b_7 $, $a_1a_2a_6a_7b_3b_5b_7 $,
$a_1a_2a_6a_8b_2b_5b_9 $, $a_1a_2a_6a_8b_2b_6b_7 $,
$a_1a_2a_6a_8b_2b_7b_9 $,

$a_1a_2a_6a_8b_5b_7b_9 $, $a_1a_2a_6b_1b_2b_3b_9 $,
$a_1a_2a_6b_1b_2b_5b_9 $, $a_1a_2a_6b_1b_3b_5b_9 $,
$a_1a_2a_6b_2b_3b_7b_9 $,

$a_1a_2a_6b_3b_5b_7b_9 $, $a_1a_2a_7a_8b_2b_3b_6 $,
$a_1a_2a_7b_1b_3b_5b_6 $, $a_1a_2a_7b_2b_3b_6b_7 $,
$a_1a_2a_7b_3b_5b_6b_7 $,

$a_1a_2a_8b_2b_3b_4b_6 $, $a_1a_2a_8b_2b_3b_4b_8 $,
$a_1a_2a_8b_2b_4b_6b_7 $, $a_1a_2a_8b_2b_4b_7b_9 $,
$a_1a_2a_8b_2b_4b_8b_9 $,

$a_1a_2a_8b_2b_5b_8b_9 $, $a_1a_2a_8b_3b_4b_6b_8 $,
$a_1a_2a_8b_4b_6b_7b_8 $, $a_1a_2a_8b_4b_7b_8b_9 $,
$a_1a_2a_8b_5b_7b_8b_9 $,

$a_1a_2b_1b_2b_3b_8b_9 $, $a_1a_2b_1b_2b_5b_8b_9 $,
$a_1a_2b_1b_3b_5b_6b_9 $, $a_1a_2b_1b_3b_6b_8b_9 $,
$a_1a_2b_1b_5b_6b_8b_9 $,

$a_1a_2b_2b_3b_4b_6b_9 $, $a_1a_2b_2b_3b_4b_8b_9 $,
$a_1a_2b_2b_3b_6b_7b_9 $, $a_1a_2b_2b_4b_6b_7b_9 $,
$a_1a_2b_3b_4b_6b_8b_9 $,

$a_1a_2b_3b_5b_6b_7b_9 $, $a_1a_2b_4b_6b_7b_8b_9 $,
$a_1a_2b_5b_6b_7b_8b_9 $, $a_1a_3a_4a_5a_7a_8b_9 $,
$a_1a_3a_4a_5a_7a_9b_9 $,

$a_1a_3a_4a_5a_8a_9b_9 $, $a_1a_3a_4a_6a_7a_8b_9 $,
$a_1a_3a_4a_6a_7a_9b_9 $, $a_1a_3a_4a_6a_8a_9b_9 $,
$a_1a_3a_5a_6a_7a_9b_4 $,

$a_1a_3a_5a_6a_7b_2b_8 $, $a_1a_3a_5a_6a_7b_4b_8 $,
$a_1a_3a_5a_6a_9b_1b_4 $, $a_1a_3a_5a_6b_1b_2b_8 $,
$a_1a_3a_5a_6b_1b_4b_8 $,

$a_1a_3a_5a_7a_8b_2b_9 $, $a_1a_3a_5a_7a_9b_4b_9 $,
$a_1a_3a_5a_7b_2b_8b_9 $, $a_1a_3a_5a_7b_4b_8b_9 $,
$a_1a_3a_5a_8a_9b_5b_9 $,

$a_1a_3a_5a_8b_2b_5b_9 $, $a_1a_3a_5a_9b_1b_4b_7 $,
$a_1a_3a_5a_9b_1b_5b_6 $, $a_1a_3a_5a_9b_1b_6b_7 $,
$a_1a_3a_5a_9b_4b_7b_9 $,

$a_1a_3a_5a_9b_5b_6b_7 $, $a_1a_3a_5a_9b_5b_7b_9 $,
$a_1a_3a_5b_1b_2b_5b_8 $, $a_1a_3a_5b_1b_4b_7b_8 $,
$a_1a_3a_5b_1b_5b_6b_8 $,

$a_1a_3a_5b_1b_6b_7b_8 $, $a_1a_3a_5b_2b_5b_8b_9 $,
$a_1a_3a_5b_4b_7b_8b_9 $, $a_1a_3a_5b_5b_6b_7b_8 $,
$a_1a_3a_5b_5b_7b_8b_9 $,

$a_1a_3a_6a_7a_8b_2b_9 $, $a_1a_3a_6a_7a_9b_4b_8 $,
$a_1a_3a_6a_7a_9b_8b_9 $, $a_1a_3a_6a_7b_2b_8b_9 $,
$a_1a_3a_6a_8a_9b_5b_9 $,

$a_1a_3a_6a_8b_2b_5b_9 $, $a_1a_3a_6a_9b_1b_4b_5 $,
$a_1a_3a_6a_9b_3b_4b_5 $, $a_1a_3a_6a_9b_3b_4b_8 $,
$a_1a_3a_6a_9b_3b_5b_9 $,

$a_1a_3a_6a_9b_3b_8b_9 $, $a_1a_3a_6b_1b_2b_5b_9 $,
$a_1a_3a_6b_1b_2b_8b_9 $, $a_1a_3a_6b_1b_3b_4b_5 $,
$a_1a_3a_6b_1b_3b_4b_8 $,

$a_1a_3a_6b_1b_3b_5b_9 $, $a_1a_3a_6b_1b_3b_8b_9 $,
$a_1a_3a_7a_9b_4b_8b_9 $, $a_1a_3a_9b_1b_4b_5b_6 $,
$a_1a_3a_9b_1b_4b_6b_7 $,

$a_1a_3a_9b_3b_4b_5b_6 $, $a_1a_3a_9b_3b_4b_6b_9 $,
$a_1a_3a_9b_3b_4b_8b_9 $, $a_1a_3a_9b_3b_5b_6b_9 $,
$a_1a_3a_9b_4b_6b_7b_9 $,

$a_1a_3a_9b_5b_6b_7b_9 $, $a_1a_3b_1b_2b_5b_8b_9 $,
$a_1a_3b_1b_3b_4b_5b_6 $, $a_1a_3b_1b_3b_4b_6b_8 $,
$a_1a_3b_1b_3b_5b_6b_9 $,

$a_1a_3b_1b_3b_6b_8b_9 $, $a_1a_3b_1b_4b_6b_7b_8 $,
$a_1a_3b_1b_5b_6b_8b_9 $, $a_1a_3b_3b_4b_6b_8b_9 $,
$a_1a_3b_4b_6b_7b_8b_9 $,

$a_1a_3b_5b_6b_7b_8b_9 $, $a_1a_4a_5a_7a_8b_1b_3 $,
$a_1a_4a_5a_7a_8b_1b_4 $, $a_1a_4a_5a_7a_8b_4b_9 $,
$a_1a_4a_5a_7a_9b_1b_4 $,

$a_1a_4a_5a_7a_9b_4b_9 $, $a_1a_4a_5a_8a_9b_5b_7 $,
$a_1a_4a_5a_8a_9b_7b_9 $, $a_1a_4a_5a_8b_1b_3b_8 $,
$a_1a_4a_5a_8b_1b_4b_7 $,

$a_1a_4a_5a_8b_1b_7b_8 $, $a_1a_4a_5a_8b_4b_7b_9 $,
$a_1a_4a_5a_8b_5b_7b_8 $, $a_1a_4a_5a_9b_1b_4b_7 $,
$a_1a_4a_5a_9b_1b_5b_6 $,

$a_1a_4a_5a_9b_1b_6b_7 $, $a_1a_4a_5a_9b_4b_7b_9 $,
$a_1a_4a_5a_9b_5b_6b_7 $, $a_1a_4a_5b_1b_5b_6b_8 $,
$a_1a_4a_5b_1b_6b_7b_8 $,

$a_1a_4a_5b_5b_6b_7b_8 $, $a_1a_4a_6a_7a_8b_2b_6 $,
$a_1a_4a_6a_7a_8b_2b_9 $, $a_1a_4a_6a_7a_9b_2b_7 $,
$a_1a_4a_6a_7a_9b_2b_9 $,

$a_1a_4a_6a_7a_9b_5b_7 $, $a_1a_4a_6a_7b_2b_6b_7 $,
$a_1a_4a_6a_8a_9b_5b_7 $, $a_1a_4a_6a_8a_9b_7b_9 $,
$a_1a_4a_6a_8b_2b_6b_7 $,

$a_1a_4a_6a_8b_2b_7b_9 $, $a_1a_4a_6a_9b_2b_7b_9 $,
$a_1a_4a_7a_8b_1b_3b_6 $, $a_1a_4a_7a_8b_1b_4b_6 $,
$a_1a_4a_7a_8b_2b_4b_6 $,

$a_1a_4a_7a_8b_2b_4b_9 $, $a_1a_4a_7a_9b_1b_4b_6 $,
$a_1a_4a_7a_9b_1b_5b_6 $, $a_1a_4a_7a_9b_2b_4b_6 $,
$a_1a_4a_7a_9b_2b_4b_9 $,

$a_1a_4a_7a_9b_2b_6b_7 $, $a_1a_4a_7a_9b_5b_6b_7 $,
$a_1a_4a_8b_1b_3b_6b_8 $, $a_1a_4a_8b_1b_4b_6b_7 $,
$a_1a_4a_8b_1b_6b_7b_8 $,

$a_1a_4a_8b_2b_4b_6b_7 $, $a_1a_4a_8b_2b_4b_7b_9 $,
$a_1a_4a_9b_1b_4b_6b_7 $, $a_1a_4a_9b_2b_4b_6b_7 $,
$a_1a_4a_9b_2b_4b_7b_9 $,

$a_1a_5a_6a_7a_9b_1b_4 $, $a_1a_5a_6a_7b_1b_3b_4 $,
$a_1a_5a_6a_7b_2b_3b_8 $, $a_1a_5a_6a_7b_3b_4b_8 $,
$a_1a_5a_6b_1b_2b_3b_8 $,

$a_1a_5a_6b_1b_3b_4b_8 $, $a_1a_5a_7a_8b_1b_3b_4 $,
$a_1a_5a_7a_8b_2b_3b_8 $, $a_1a_5a_7a_8b_2b_8b_9 $,
$a_1a_5a_7a_8b_3b_4b_8 $,

$a_1a_5a_7a_8b_4b_8b_9 $, $a_1a_5a_8a_9b_5b_7b_9 $,
$a_1a_5a_8b_1b_3b_4b_8 $, $a_1a_5a_8b_1b_4b_7b_8 $,
$a_1a_5a_8b_2b_5b_8b_9 $,

$a_1a_5a_8b_4b_7b_8b_9 $, $a_1a_5a_8b_5b_7b_8b_9 $,
$a_1a_6a_7a_9b_1b_4b_5 $, $a_1a_6a_7a_9b_2b_3b_7 $,
$a_1a_6a_7a_9b_2b_3b_8 $,

$a_1a_6a_7a_9b_2b_8b_9 $, $a_1a_6a_7a_9b_3b_4b_5 $,
$a_1a_6a_7a_9b_3b_4b_8 $, $a_1a_6a_7a_9b_3b_5b_7 $,
$a_1a_6a_7b_1b_3b_4b_5 $,

$a_1a_6a_8a_9b_5b_7b_9 $, $a_1a_6a_9b_2b_3b_7b_9 $,
$a_1a_6a_9b_2b_3b_8b_9 $, $a_1a_6a_9b_3b_5b_7b_9 $,
$a_1a_6b_1b_2b_3b_8b_9 $,

$a_1a_7a_8b_1b_3b_4b_6 $, $a_1a_7a_8b_2b_3b_4b_6 $,
$a_1a_7a_8b_2b_3b_4b_8 $, $a_1a_7a_8b_2b_4b_8b_9 $,
$a_1a_7a_9b_1b_4b_5b_6 $,

$a_1a_7a_9b_2b_3b_4b_6 $, $a_1a_7a_9b_2b_3b_4b_8 $,
$a_1a_7a_9b_2b_3b_6b_7 $, $a_1a_7a_9b_2b_4b_8b_9 $,
$a_1a_7a_9b_3b_4b_5b_6 $,

$a_1a_7a_9b_3b_5b_6b_7 $, $a_1a_7b_1b_3b_4b_5b_6 $,
$a_1a_8b_1b_3b_4b_6b_8 $, $a_1a_8b_1b_4b_6b_7b_8 $,
$a_1a_9b_2b_3b_4b_6b_9 $,

$a_1a_9b_2b_3b_4b_8b_9 $, $a_1a_9b_2b_3b_6b_7b_9 $,
$a_1a_9b_2b_4b_6b_7b_9 $, $a_1a_9b_3b_5b_6b_7b_9 $,
$a_2a_3a_4a_5a_7a_8b_3$,

$a_2a_3a_4a_5a_7a_9b_9 $, $a_2a_3a_4a_5a_7b_3b_9$,
$a_2a_3a_4a_5a_8a_9b_3 $, $a_2a_3a_4a_5a_9b_3b_9$,
$a_2a_3a_4a_6a_7a_8b_6 $,

$a_2a_3a_4a_6a_7a_9b_8$, $a_2a_3a_4a_6a_7b_6b_7 $,
$a_2a_3a_4a_6a_7b_7b_8$, $a_2a_3a_4a_6a_8a_9b_8 $,
$a_2a_3a_4a_6a_8b_6b_7$,

$a_2a_3a_4a_6a_8b_7b_8 $, $a_2a_3a_4a_7a_8b_3b_6$,
$a_2a_3a_4a_7a_9b_8b_9 $, $a_2a_3a_4a_7b_3b_6b_9$,
$a_2a_3a_4a_7b_6b_7b_8 $,

$a_2a_3a_4a_7b_6b_8b_9$, $a_2a_3a_4a_8a_9b_3b_8 $,
$a_2a_3a_4a_8b_3b_6b_8$, $a_2a_3a_4a_8b_6b_7b_8 $,
$a_2a_3a_4a_9b_3b_8b_9$,

$a_2a_3a_4b_3b_6b_8b_9 $, $a_2a_3a_5a_6a_7a_9b_4$,
$a_2a_3a_5a_6a_7b_2b_7 $, $a_2a_3a_5a_6a_7b_4b_7$,
$a_2a_3a_5a_6a_9b_1b_4 $,

$a_2a_3a_5a_6b_1b_2b_7$, $a_2a_3a_5a_6b_1b_4b_7 $,
$a_2a_3a_5a_7a_8b_2b_3$, $a_2a_3a_5a_7a_9b_4b_9 $,
$a_2a_3a_5a_7b_2b_3b_7$,

$a_2a_3a_5a_7b_3b_7b_9 $, $a_2a_3a_5a_7b_4b_7b_9$,
$a_2a_3a_5a_8a_9b_2b_3 $, $a_2a_3a_5a_8a_9b_2b_5$,
$a_2a_3a_5a_9b_1b_2b_5 $,

$a_2a_3a_5a_9b_1b_2b_7$, $a_2a_3a_5a_9b_1b_4b_7 $,
$a_2a_3a_5a_9b_2b_3b_7$, $a_2a_3a_5a_9b_3b_7b_9 $,
$a_2a_3a_5a_9b_4b_7b_9$,

$a_2a_3a_6a_7a_8b_2b_6 $, $a_2a_3a_6a_7a_9b_4b_8$,
$a_2a_3a_6a_7b_2b_6b_7 $, $a_2a_3a_6a_7b_4b_7b_8$,
$a_2a_3a_6a_8a_9b_4b_5 $,

$a_2a_3a_6a_8a_9b_4b_8$, $a_2a_3a_6a_8b_1b_2b_5 $,
$a_2a_3a_6a_8b_1b_2b_7$, $a_2a_3a_6a_8b_1b_4b_5 $,
$a_2a_3a_6a_8b_1b_4b_7$,

$a_2a_3a_6a_8b_2b_6b_7 $, $a_2a_3a_6a_8b_4b_7b_8$,
$a_2a_3a_6a_9b_1b_4b_5 $, $a_2a_3a_7a_8b_2b_3b_6$,
$a_2a_3a_7a_9b_4b_8b_9 $,

$a_2a_3a_7b_2b_3b_6b_7$, $a_2a_3a_7b_3b_6b_7b_9 $,
$a_2a_3a_7b_4b_7b_8b_9$, $a_2a_3a_7b_6b_7b_8b_9 $,
$a_2a_3a_8a_9b_2b_3b_4$,

$a_2a_3a_8a_9b_2b_4b_5 $, $a_2a_3a_8a_9b_3b_4b_8$,
$a_2a_3a_8b_1b_2b_4b_5 $, $a_2a_3a_8b_1b_2b_4b_7$,
$a_2a_3a_8b_2b_3b_4b_6 $,

$a_2a_3a_8b_2b_4b_6b_7$, $a_2a_3a_8b_3b_4b_6b_8 $,
$a_2a_3a_8b_4b_6b_7b_8$, $a_2a_3a_9b_1b_2b_4b_5 $,
$a_2a_3a_9b_1b_2b_4b_7$,

$a_2a_3a_9b_2b_3b_4b_6 $, $a_2a_3a_9b_2b_3b_6b_7$,
$a_2a_3a_9b_2b_4b_6b_7 $, $a_2a_3a_9b_3b_4b_6b_9$,
$a_2a_3a_9b_3b_4b_8b_9 $,

$a_2a_3a_9b_3b_6b_7b_9$, $a_2a_3a_9b_4b_6b_7b_9 $,
$a_2a_3b_3b_4b_6b_8b_9$, $a_2a_3b_4b_6b_7b_8b_9 $,
$a_2a_4a_5a_7a_9b_1b_4$,

$a_2a_4a_5a_7a_9b_4b_9 $, $a_2a_4a_5a_7b_1b_3b_9$,
$a_2a_4a_5a_7b_1b_4b_9 $, $a_2a_4a_5a_8a_9b_3b_8$,
$a_2a_4a_5a_8a_9b_5b_8 $,

$a_2a_4a_5a_9b_1b_2b_3$, $a_2a_4a_5a_9b_1b_2b_5 $,
$a_2a_4a_5a_9b_1b_3b_9$, $a_2a_4a_5a_9b_1b_4b_9 $,
$a_2a_4a_5a_9b_2b_3b_8$,

$a_2a_4a_5a_9b_2b_5b_8 $, $a_2a_4a_5b_1b_2b_3b_8$,
$a_2a_4a_5b_1b_2b_5b_8 $, $a_2a_4a_6a_7a_9b_5b_8$,
$a_2a_4a_6a_7b_5b_7b_8 $,

$a_2a_4a_6a_8a_9b_5b_8$, $a_2a_4a_6a_8b_5b_7b_8 $,
$a_2a_4a_7a_9b_1b_4b_5$, $a_2a_4a_7a_9b_4b_5b_8 $,
$a_2a_4a_7a_9b_4b_8b_9$,

$a_2a_4a_7b_1b_3b_6b_9 $, $a_2a_4a_7b_1b_4b_5b_9$,
$a_2a_4a_7b_1b_5b_6b_8 $, $a_2a_4a_7b_1b_5b_8b_9$,
$a_2a_4a_7b_1b_6b_8b_9 $,

$a_2a_4a_7b_4b_5b_8b_9$, $a_2a_4a_7b_5b_6b_7b_8 $,
$a_2a_4a_9b_1b_2b_3b_9$, $a_2a_4a_9b_1b_2b_4b_5 $,
$a_2a_4a_9b_1b_2b_4b_9$,

$a_2a_4a_9b_2b_3b_8b_9 $, $a_2a_4a_9b_2b_4b_5b_8$,
$a_2a_4a_9b_2b_4b_8b_9 $, $a_2a_4b_1b_2b_3b_8b_9$,
$a_2a_4b_1b_2b_4b_5b_9 $,

$a_2a_4b_1b_2b_5b_8b_9$, $a_2a_4b_1b_3b_6b_8b_9 $,
$a_2a_4b_2b_4b_5b_8b_9$, $a_2a_5a_6a_7a_9b_1b_4 $,
$a_2a_5a_6a_7b_1b_3b_7$,

$a_2a_5a_6a_7b_1b_4b_7 $, $a_2a_5a_6a_7b_2b_3b_7$,
$a_2a_5a_6b_1b_2b_3b_7 $, $a_2a_5a_7b_1b_3b_7b_9$,
$a_2a_5a_7b_1b_4b_7b_9 $,

$a_2a_5a_8a_9b_2b_3b_8$, $a_2a_5a_8a_9b_2b_5b_8 $,
$a_2a_5a_9b_1b_2b_3b_7$, $a_2a_5a_9b_1b_3b_7b_9 $,
$a_2a_5a_9b_1b_4b_7b_9$,

$a_2a_6a_7a_9b_1b_4b_5 $, $a_2a_6a_7a_9b_4b_5b_8$,
$a_2a_6a_7b_1b_3b_5b_7 $, $a_2a_6a_7b_1b_4b_5b_7$,
$a_2a_6a_7b_4b_5b_7b_8 $,

$a_2a_6a_8a_9b_4b_5b_8$, $a_2a_6a_8b_1b_2b_5b_9 $,
$a_2a_6a_8b_1b_2b_7b_9$, $a_2a_6a_8b_1b_4b_5b_7 $,
$a_2a_6a_8b_1b_5b_7b_9$,

$a_2a_6a_8b_4b_5b_7b_8 $, $a_2a_6b_1b_2b_3b_7b_9$,
$a_2a_6b_1b_3b_5b_7b_9 $, $a_2a_7b_1b_3b_5b_6b_9$,
$a_2a_7b_1b_3b_5b_7b_9 $,

$a_2a_7b_1b_4b_5b_7b_9$, $a_2a_7b_1b_5b_6b_8b_9 $,
$a_2a_7b_3b_5b_6b_7b_9$, $a_2a_7b_4b_5b_7b_8b_9 $,
$a_2a_7b_5b_6b_7b_8b_9$,

$a_2a_8a_9b_2b_3b_4b_8 $, $a_2a_8a_9b_2b_4b_5b_8$,
$a_2a_8b_1b_2b_4b_5b_9 $, $a_2a_8b_1b_2b_4b_7b_9$,
$a_2a_8b_1b_4b_5b_7b_9 $,

$a_2a_8b_2b_4b_5b_8b_9$, $a_2a_8b_4b_5b_7b_8b_9 $,
$a_2a_9b_1b_2b_3b_7b_9$, $a_2a_9b_1b_2b_4b_7b_9 $,
$a_2a_9b_2b_3b_4b_6b_9$,

$a_2a_9b_2b_3b_4b_8b_9 $, $a_2a_9b_2b_3b_6b_7b_9$,
$a_2a_9b_2b_4b_6b_7b_9 $, $a_3a_4a_5a_7a_8b_3b_9$,
$a_3a_4a_5a_8a_9b_3b_9 $,

$a_3a_4a_6a_7a_8b_2b_6$, $a_3a_4a_6a_7a_8b_2b_9 $,
$a_3a_4a_6a_7a_9b_8b_9$, $a_3a_4a_6a_7b_2b_6b_7 $,
$a_3a_4a_6a_7b_2b_7b_8$,

$a_3a_4a_6a_7b_2b_8b_9 $, $a_3a_4a_6a_8a_9b_3b_8$,
$a_3a_4a_6a_8a_9b_3b_9 $, $a_3a_4a_6a_8b_1b_2b_7$,
$a_3a_4a_6a_8b_1b_2b_9 $,

$a_3a_4a_6a_8b_1b_3b_8$, $a_3a_4a_6a_8b_1b_3b_9 $,
$a_3a_4a_6a_8b_1b_7b_8$, $a_3a_4a_6a_8b_2b_6b_7 $,
$a_3a_4a_6a_9b_3b_8b_9$,

$a_3a_4a_6b_1b_2b_7b_8 $, $a_3a_4a_6b_1b_2b_8b_9$,
$a_3a_4a_6b_1b_3b_8b_9 $, $a_3a_4a_7a_8b_1b_2b_6$,
$a_3a_4a_7a_8b_1b_2b_9 $,

$a_3a_4a_7a_8b_1b_3b_6$, $a_3a_4a_7a_8b_1b_3b_9 $,
$a_3a_4a_7b_1b_2b_6b_8$, $a_3a_4a_7b_1b_2b_8b_9 $,
$a_3a_4a_7b_1b_3b_6b_9$,

$a_3a_4a_7b_1b_6b_8b_9 $, $a_3a_4a_7b_2b_6b_7b_8$,
$a_3a_4a_8b_1b_2b_6b_7 $, $a_3a_4a_8b_1b_3b_6b_8$,
$a_3a_4a_8b_1b_6b_7b_8 $,

$a_3a_4b_1b_2b_6b_7b_8$, $a_3a_4b_1b_3b_6b_8b_9 $,
$a_3a_5a_6a_7b_2b_7b_8$, $a_3a_5a_6a_7b_4b_7b_8 $,
$a_3a_5a_6b_1b_2b_7b_8$,

$a_3a_5a_6b_1b_4b_7b_8 $, $a_3a_5a_7a_8b_2b_3b_5$,
$a_3a_5a_7a_8b_2b_5b_9 $, $a_3a_5a_7a_8b_3b_5b_9$,
$a_3a_5a_7b_2b_3b_5b_7 $,

$a_3a_5a_7b_2b_5b_7b_8$, $a_3a_5a_7b_2b_5b_8b_9 $,
$a_3a_5a_7b_3b_5b_7b_9$, $a_3a_5a_7b_4b_7b_8b_9 $,
$a_3a_5a_7b_5b_7b_8b_9$,

$a_3a_5a_8a_9b_2b_3b_5 $, $a_3a_5a_8a_9b_3b_5b_9$,
$a_3a_5a_9b_1b_2b_5b_6 $, $a_3a_5a_9b_1b_2b_6b_7$,
$a_3a_5a_9b_2b_3b_5b_7 $,

$a_3a_5a_9b_2b_5b_6b_7$, $a_3a_5a_9b_3b_5b_7b_9 $,
$a_3a_5b_1b_2b_5b_6b_8$, $a_3a_5b_1b_2b_6b_7b_8 $,
$a_3a_5b_2b_5b_6b_7b_8$,

$a_3a_6a_8a_9b_3b_4b_5 $, $a_3a_6a_8a_9b_3b_4b_8$,
$a_3a_6a_8a_9b_3b_5b_9 $, $a_3a_6a_8b_1b_2b_5b_9$,
$a_3a_6a_8b_1b_3b_4b_5 $,

$a_3a_6a_8b_1b_3b_4b_8$, $a_3a_6a_8b_1b_3b_5b_9 $,
$a_3a_6a_8b_1b_4b_7b_8$, $a_3a_7a_8b_1b_2b_5b_6 $,
$a_3a_7a_8b_1b_2b_5b_9$,

$a_3a_7a_8b_1b_3b_5b_6 $, $a_3a_7a_8b_1b_3b_5b_9$,
$a_3a_7a_8b_2b_3b_5b_6 $, $a_3a_7b_1b_2b_5b_6b_8$,
$a_3a_7b_1b_2b_5b_8b_9 $,

$a_3a_7b_1b_3b_5b_6b_9$, $a_3a_7b_1b_5b_6b_8b_9 $,
$a_3a_7b_2b_3b_5b_6b_7$, $a_3a_7b_2b_5b_6b_7b_8 $,
$a_3a_7b_3b_5b_6b_7b_9$,

$a_3a_7b_5b_6b_7b_8b_9 $, $a_3a_8a_9b_2b_3b_4b_5$,
$a_3a_8b_1b_2b_4b_5b_6 $, $a_3a_8b_1b_2b_4b_6b_7$,
$a_3a_8b_1b_3b_4b_5b_6 $,

$a_3a_8b_1b_3b_4b_6b_8$, $a_3a_8b_1b_4b_6b_7b_8 $,
$a_3a_8b_2b_3b_4b_5b_6$, $a_3a_9b_1b_2b_4b_5b_6 $,
$a_3a_9b_1b_2b_4b_6b_7$,

$a_3a_9b_2b_3b_4b_5b_6 $, $a_3a_9b_2b_3b_5b_6b_7$,
$a_3a_9b_3b_5b_6b_7b_9 $, $a_4a_5a_7a_8b_1b_3b_9$,
$a_4a_5a_7a_8b_1b_4b_9 $,

$a_4a_5a_8a_9b_3b_7b_8$, $a_4a_5a_8a_9b_3b_7b_9 $,
$a_4a_5a_8a_9b_5b_7b_8$, $a_4a_5a_8b_1b_3b_7b_8 $,
$a_4a_5a_8b_1b_3b_7b_9$,

$a_4a_5a_8b_1b_4b_7b_9 $, $a_4a_5a_9b_1b_2b_3b_7$,
$a_4a_5a_9b_1b_2b_5b_6 $, $a_4a_5a_9b_1b_2b_6b_7$,
$a_4a_5a_9b_1b_3b_7b_9 $,

$a_4a_5a_9b_1b_4b_7b_9$, $a_4a_5a_9b_2b_3b_7b_8 $,
$a_4a_5a_9b_2b_5b_6b_7$, $a_4a_5a_9b_2b_5b_7b_8 $,
$a_4a_5b_1b_2b_3b_7b_8$,

$a_4a_5b_1b_2b_5b_6b_8 $, $a_4a_5b_1b_2b_6b_7b_8$,
$a_4a_5b_2b_5b_6b_7b_8 $, $a_4a_6a_7a_9b_2b_7b_8$,
$a_4a_6a_7a_9b_2b_8b_9 $,

$a_4a_6a_7a_9b_5b_7b_8$, $a_4a_6a_8a_9b_3b_7b_8 $,
$a_4a_6a_8a_9b_3b_7b_9$, $a_4a_6a_8a_9b_5b_7b_8 $,
$a_4a_6a_8b_1b_2b_7b_9$,

$a_4a_6a_8b_1b_3b_7b_8 $, $a_4a_6a_8b_1b_3b_7b_9$,
$a_4a_6a_9b_2b_3b_7b_8 $, $a_4a_6a_9b_2b_3b_7b_9$,
$a_4a_6a_9b_2b_3b_8b_9 $,

$a_4a_6b_1b_2b_3b_7b_8$, $a_4a_6b_1b_2b_3b_7b_9 $,
$a_4a_6b_1b_2b_3b_8b_9$, $a_4a_7a_8b_1b_2b_4b_6 $,
$a_4a_7a_8b_1b_2b_4b_9$,

$a_4a_7a_9b_1b_2b_4b_5 $, $a_4a_7a_9b_1b_2b_4b_6$,
$a_4a_7a_9b_1b_2b_5b_6 $, $a_4a_7a_9b_2b_4b_5b_8$,
$a_4a_7a_9b_2b_4b_8b_9 $,

$a_4a_7a_9b_2b_5b_6b_7$, $a_4a_7a_9b_2b_5b_7b_8 $,
$a_4a_7b_1b_2b_4b_5b_9$, $a_4a_7b_1b_2b_5b_6b_8 $,
$a_4a_7b_1b_2b_5b_8b_9$,

$a_4a_7b_2b_4b_5b_8b_9 $, $a_4a_7b_2b_5b_6b_7b_8$,
$a_4a_8b_1b_2b_4b_6b_7 $, $a_4a_8b_1b_2b_4b_7b_9$,
$a_4a_9b_1b_2b_3b_7b_9 $,

$a_4a_9b_1b_2b_4b_6b_7$, $a_4a_9b_1b_2b_4b_7b_9 $,
$a_5a_6a_7b_1b_3b_4b_7$, $a_5a_6a_7b_2b_3b_7b_8 $,
$a_5a_6a_7b_3b_4b_7b_8$,

$a_5a_6b_1b_2b_3b_7b_8 $, $a_5a_6b_1b_3b_4b_7b_8$,
$a_5a_7a_8b_1b_3b_4b_7 $, $a_5a_7a_8b_1b_3b_7b_9$,
$a_5a_7a_8b_1b_4b_7b_9 $,

$a_5a_7a_8b_2b_3b_5b_8$, $a_5a_7a_8b_2b_5b_8b_9 $,
$a_5a_7a_8b_3b_4b_7b_8$, $a_5a_7a_8b_3b_5b_7b_8 $,
$a_5a_7a_8b_3b_5b_7b_9$,

$a_5a_7a_8b_4b_7b_8b_9 $, $a_5a_7a_8b_5b_7b_8b_9$,
$a_5a_7b_2b_3b_5b_7b_8 $, $a_5a_8a_9b_2b_3b_5b_8$,
$a_5a_8a_9b_3b_5b_7b_8 $,

$a_5a_8a_9b_3b_5b_7b_9$, $a_5a_8b_1b_3b_4b_7b_8 $,
$a_5a_9b_2b_3b_5b_7b_8$, $a_6a_7a_9b_2b_3b_7b_8 $,
$a_6a_7a_9b_3b_4b_5b_8$,

$a_6a_7a_9b_3b_5b_7b_8 $, $a_6a_7b_1b_3b_4b_5b_7$,
$a_6a_7b_3b_4b_5b_7b_8 $, $a_6a_8a_9b_3b_4b_5b_8$,
$a_6a_8a_9b_3b_5b_7b_8 $,

$a_6a_8a_9b_3b_5b_7b_9$, $a_6a_8b_1b_3b_4b_5b_7 $,
$a_6a_8b_1b_3b_4b_7b_8$, $a_6a_8b_1b_3b_5b_7b_9 $,
$a_6a_8b_3b_4b_5b_7b_8$,

$a_7a_8b_1b_2b_4b_5b_6 $, $a_7a_8b_1b_2b_4b_5b_9$,
$a_7a_8b_1b_3b_4b_5b_6 $, $a_7a_8b_1b_3b_4b_5b_7$,
$a_7a_8b_1b_3b_5b_7b_9 $,

$a_7a_8b_1b_4b_5b_7b_9$, $a_7a_8b_2b_3b_4b_5b_6 $,
$a_7a_8b_2b_3b_4b_5b_8$, $a_7a_8b_2b_4b_5b_8b_9 $,
$a_7a_8b_3b_4b_5b_7b_8$,

$a_7a_8b_4b_5b_7b_8b_9 $, $a_7a_9b_1b_2b_4b_5b_6$,
$a_7a_9b_2b_3b_4b_5b_6 $, $a_7a_9b_2b_3b_4b_5b_8$,
$a_7a_9b_2b_3b_5b_6b_7 $,

$a_7a_9b_2b_3b_5b_7b_8$, $a_8a_9b_2b_3b_4b_5b_8$.

\end{document}